%
\input hyperbasics 
\def\unredoffs{} 

%
%
\newbox\leftpage \newdimen\fullhsize \newdimen\hstitle \newdimen\hsbody
\tolerance=1000\hfuzz=2pt
\catcode`\@=11 
\ifx\hyperdef\UNd@FiNeD\def\hyperdef#1#2#3#4{#4}\def\hyperref#1#2#3#4{#4}\fi
\def\bigans{b }
\let\answ\bigans
\magnification=1100\unredoffs\baselineskip=16pt plus 2pt minus 1pt
\hsbody=\hsize \hstitle=\hsize 
%
%
\newcount\yearltd\yearltd=\year\advance\yearltd by -1900

\def\Title#1#2{\nopagenumbers\abstractfont\hsize=\hstitle\rightline{#1}%
\vskip 1in\centerline{\titlefont #2}\abstractfont\vskip .5in\pageno=0}
\def\Date#1{\vfill\leftline{#1}\tenpoint\supereject\global\hsize=\hsbody%
\footline={\hss\tenrm\hyperdef\hypernoname{page}\folio\folio\hss}}%
%

\def\draftmode{\message{ DRAFTMODE }\def\draftdate{{\rm preliminary draft:
\number\month/\number\day/\number\yearltd\ \ \hourmin}}%
\headline={\hfil\draftdate}\writelabels\baselineskip=20pt plus 2pt minus 2pt
 {\count255=\time\divide\count255 by 60 \xdef\hourmin{\number\count255}
  \multiply\count255 by-60\advance\count255 by\time
  \xdef\hourmin{\hourmin:\ifnum\count255<10 0\fi\the\count255}}}
\def\nolabels{\def\wrlabeL##1{}\def\eqlabeL##1{}\def\reflabeL##1{}}
\def\writelabels{\def\wrlabeL##1{\leavevmode\vadjust{\rlap{\smash%
{\line{{\escapechar=` \hfill\rlap{\sevenrm\hskip.03in\string##1}}}}}}}%
\def\eqlabeL##1{{\escapechar-1\rlap{\sevenrm\hskip.05in\string##1}}}%
\def\reflabeL##1{\noexpand\llap{\noexpand\sevenrm\string\string\string##1}}}
\nolabels
%
\global\newcount\secno \global\secno=0
\global\newcount\meqno \global\meqno=1
\def\newsec#1{\global\advance\secno by1\message{(\the\secno. #1)}
\global\subsecno=0\eqnres@t\teores@t\noindent
{\bf\hyperdef\hypernoname{section}{\the\secno}{\the\secno.} \noexpand#1}%
\writetoca{{\string\hyperref{}{section}{\the\secno}{\secsym}} {#1}}%
\par\nobreak\medskip\nobreak}
\def\eqnres@t{\xdef\secsym{\the\secno.}\global\meqno=1\bigbreak\bigskip}
\def\teores@t{\xdef\secsym{\the\secno.}\global\theno=1}
\def\sequentialequations{\def\eqnres@t{\bigbreak}}\xdef\secsym{}
\global\newcount\subsecno \global\subsecno=0
\def\subsec#1{\global\advance\subsecno by1\message{(\secsym\the\subsecno. #1)}
\ifnum\lastpenalty>9000\else\bigbreak\fi
\noindent{\it\hyperdef\hypernoname{subsection}{\secsym\the\subsecno}%
{\secsym\the\subsecno.} #1}\writetoca{\string\quad
{\string\hyperref{}{subsection}{\secsym\the\subsecno}{\secsym\the\subsecno.}}
{#1}}\par\nobreak\medskip\nobreak}
\def\appendix#1#2{\global\meqno=1\global\subsecno=0\xdef\secsym{\hbox{#1.}}%
\bigbreak\bigskip\noindent{\bf Appendix \hyperdef\hypernoname{appendix}{#1}%
{#1.} #2}\message{(#1. #2)}\xdef\appsym{#1}%
\writetoca{\string\hyperref{}{appendix}{#1}{Appendix {#1.}} {#2}}%
\par\nobreak\medskip\nobreak}
%
%
\def\checkm@de#1#2{\ifmmode{\def\f@rst##1{##1}\hyperdef\hypernoname{equation}%
{#1}{#2}}\else\hyperref{}{equation}{#1}{#2}\fi}
\def\eqnn#1{\DefWarn#1\xdef #1{(\noexpand\relax\noexpand\checkm@de%
{\secsym\the\meqno}{\secsym\the\meqno})}%
\wrlabeL#1\writedef{#1\leftbracket#1}\global\advance\meqno by1}
\def\f@rst#1{\c@t#1a\em@ark}\def\c@t#1#2\em@ark{#1}
\def\eqna#1{\DefWarn#1\wrlabeL{#1$\{\}$}%
\xdef #1##1{(\noexpand\relax\noexpand\checkm@de%
{\secsym\the\meqno\noexpand\f@rst{##1}}{\hbox{$\secsym\the\meqno##1$}})}
\writedef{#1\numbersign1\leftbracket#1{\numbersign1}}\global\advance\meqno by1}
\def\eqn#1#2{\DefWarn#1%
\xdef #1{(\noexpand\hyperref{}{equation}{\secsym\the\meqno}%
{\secsym\the\meqno})}$$#2\eqno(\hyperdef\hypernoname{equation}%
{\secsym\the\meqno}{\secsym\the\meqno})\eqlabeL#1$$%
\writedef{#1\leftbracket#1}\global\advance\meqno by1}
\def\xeqn{\expandafter\xe@n}\def\xe@n(#1){#1}
\def\xeqna#1{\expandafter\xe@n#1}
\def\eqns#1{(\e@ns #1{\hbox{}})}
\def\e@ns#1{\ifx\UNd@FiNeD#1\message{eqnlabel \string#1 is undefined.}%
\xdef#1{(?.?)}\fi{\let\hyperref=\relax\xdef\next{#1}}%
\ifx\next\em@rk\def\next{}\else%
\ifx\next#1\xeqn#1\else\def\n@xt{#1}\ifx\n@xt\next#1\else\xeqna#1\fi
\fi\let\next=\e@ns\fi\next}

\def\DefWarn#1{\ifx\UNd@FiNeD#1\else
\immediate\write16{*** WARNING: the label \string#1 is already defined ***}\fi}
%
\newskip\footskip\footskip14pt plus 1pt minus 1pt 
\def\footnotefont{\ninepoint}\def\f@t#1{\footnotefont #1\@foot}
\def\f@@t{\baselineskip\footskip\bgroup\footnotefont\aftergroup\@foot\let\next}
\setbox\strutbox=\hbox{\vrule height9.5pt depth4.5pt width0pt}
\global\newcount\ftno \global\ftno=0
\def\foot{\global\advance\ftno by1\def\foot@rg{\hyperref{}{footnote}%
{\the\ftno}{\the\ftno}\xdef\foot@rg{\noexpand\hyperdef\noexpand\hypernoname%
{footnote}{\the\ftno}{\the\ftno}}}\footnote{$^{\foot@rg}$}}
%
\newwrite\ftfile
\def\footend{\def\foot{\global\advance\ftno by1\chardef\wfile=\ftfile
\hyperref{}{footnote}{\the\ftno}{$^{\the\ftno}$}%
\ifnum\ftno=1\immediate\openout\ftfile=\jobname.fts\fi%
\immediate\write\ftfile{\noexpand\smallskip%
\noexpand\item{\noexpand\hyperdef\noexpand\hypernoname{footnote}
{\the\ftno}{f\the\ftno}:\ }\pctsign}\findarg}%
\def\footatend{\vfill\eject\immediate\closeout\ftfile{\parindent=20pt
\centerline{\bf Footnotes}\nobreak\bigskip\input \jobname.fts }}}
\def\footatend{}
%
%
\global\newcount\refno \global\refno=1
\newwrite\rfile
\def\ref{[\hyperref{}{reference}{\the\refno}{\the\refno}]\nref}
\def\nref#1{\DefWarn#1%
\xdef#1{[\noexpand\hyperref{}{reference}{\the\refno}{\the\refno}]}%
\writedef{#1\leftbracket#1}%
\ifnum\refno=1\immediate\openout\rfile=\jobname.refs\fi
\chardef\wfile=\rfile\immediate\write\rfile{\noexpand\item{[\noexpand\hyperdef%
\noexpand\hypernoname{reference}{\the\refno}{\the\refno}]\ }%
\reflabeL{#1\hskip.31in}\pctsign}\global\advance\refno by1\findarg}
\def\findarg#1#{\begingroup\obeylines\newlinechar=`\^^M\pass@rg}
{\obeylines\gdef\pass@rg#1{\writ@line\relax #1^^M\hbox{}^^M}%
\gdef\writ@line#1^^M{\expandafter\toks0\expandafter{\striprel@x #1}%
\edef\next{\the\toks0}\ifx\next\em@rk\let\next=\endgroup\else\ifx\next\empty%
\else\immediate\write\wfile{\the\toks0}\fi\let\next=\writ@line\fi\next\relax}}
\def\striprel@x#1{} \def\em@rk{\hbox{}}
\def\lref{\begingroup\obeylines\lr@f}
\def\lr@f#1#2{\DefWarn#1\gdef#1{\let#1=\UNd@FiNeD\ref#1{#2}}\endgroup\unskip}

\def\addref#1{\immediate\write\rfile{\noexpand\item{}#1}} 
\def\listrefs{\footatend\vfill\supereject\immediate\closeout\rfile\writestoppt
\baselineskip=\footskip\centerline{{\bf References}}\bigskip{\parindent=20pt%
\frenchspacing\escapechar=` \input \jobname.refs\vfill\eject}\nonfrenchspacing}
\def\startrefs#1{\immediate\openout\rfile=\jobname.refs\refno=#1}
\def\xref{\expandafter\xr@f}\def\xr@f[#1]{#1}
\def\refs#1{\count255=1[\r@fs #1{\hbox{}}]}
\def\r@fs#1{\ifx\UNd@FiNeD#1\message{reflabel \string#1 is undefined.}%
\nref#1{need to supply reference \string#1.}\fi%
\vphantom{\hphantom{#1}}{\let\hyperref=\relax\xdef\next{#1}}%
\ifx\next\em@rk\def\next{}%
\else\ifx\next#1\ifodd\count255\relax\xref#1\count255=0\fi%
\else#1\count255=1\fi\let\next=\r@fs\fi\next}
%

%
\newwrite\ffile\global\newcount\figno \global\figno=1
\def\fig{fig.~\hyperref{}{figure}{\the\figno}{\the\figno}\nfig}
\def\nfig#1{\DefWarn#1%
\xdef#1{fig.~\noexpand\hyperref{}{figure}{\the\figno}{\the\figno}}%
\writedef{#1\leftbracket fig.\noexpand~\xfig#1}%
\ifnum\figno=1\immediate\openout\ffile=\jobname.figs\fi\chardef\wfile=\ffile%
{\let\hyperref=\relax
\immediate\write\ffile{\noexpand\medskip\noexpand\item{Fig.\ %
\noexpand\hyperdef\noexpand\hypernoname{figure}{\the\figno}{\the\figno}. }
\reflabeL{#1\hskip.55in}\pctsign}}\global\advance\figno by1\findarg}
\def\listfigs{\vfill\eject\immediate\closeout\ffile{\parindent40pt
\baselineskip14pt\centerline{{\bf Figure Captions}}\nobreak\medskip
\escapechar=` \input \jobname.figs\vfill\eject}}
\def\xfig{\expandafter\xf@g}\def\xf@g fig.\penalty\@M\ {}
\def\figs#1{figs.~\f@gs #1{\hbox{}}}
\def\f@gs#1{{\let\hyperref=\relax\xdef\next{#1}}\ifx\next\em@rk\def\next{}\else
\ifx\next#1\xfig #1\else#1\fi\let\next=\f@gs\fi\next}
\def\figin{\epsfcheck\figin}\def\figins{\epsfcheck\figins}
\def\epsfcheck{\ifx\epsfbox\UNd@FiNeD
\message{(NO epsf.tex, FIGURES WILL BE IGNORED)}
\gdef\figin##1{\vskip2in}\gdef\figins##1{\hskip.5in}
\else\message{(FIGURES WILL BE INCLUDED)}%
\gdef\figin##1{##1}\gdef\figins##1{##1}\fi}
\def\DefWarn#1{}
\def\figinsert{\goodbreak\midinsert}
\def\ifig#1#2#3{\DefWarn#1\xdef#1{fig.~\noexpand\hyperref{}{figure}%
{\the\figno}{\the\figno}}\writedef{#1\leftbracket fig.\noexpand~\xfig#1}%
\figinsert\figin{\centerline{#3}}\medskip\centerline{\vbox{\baselineskip12pt
\advance\hsize by -1truein\noindent\wrlabeL{#1=#1}\footnotefont%
{\bf Fig.~\hyperdef\hypernoname{figure}{\the\figno}{\the\figno}:} #2}}
\bigskip\endinsert\global\advance\figno by1}
\newwrite\lfile
{\escapechar-1\xdef\pctsign{\string\%}\xdef\leftbracket{\string\{}
\xdef\rightbracket{\string\}}\xdef\numbersign{\string\#}}
\def\writedefs{\immediate\openout\lfile=\jobname.defs \def\writedef##1{%
{\let\hyperref=\relax\let\hyperdef=\relax\let\hypernoname=\relax
 \immediate\write\lfile{\string\def\string##1\rightbracket}}}}%
\def\writestop{\def\writestoppt{\immediate\write\lfile{\string\pageno%
\the\pageno\string\startrefs\leftbracket\the\refno\rightbracket%
\string\def\string\secsym\leftbracket\secsym\rightbracket%
\string\secno\the\secno\string\meqno\the\meqno}\immediate\closeout\lfile}}
\def\writestoppt{}\def\writedef#1{}
\def\seclab#1{\DefWarn#1%
\xdef #1{\noexpand\hyperref{}{section}{\the\secno}{\the\secno}}%
\writedef{#1\leftbracket#1}\wrlabeL{#1=#1}}
\def\subseclab#1{\DefWarn#1%
\xdef #1{\noexpand\hyperref{}{subsection}{\secsym\the\subsecno}%
{\secsym\the\subsecno}}\writedef{#1\leftbracket#1}\wrlabeL{#1=#1}}
\def\applab#1{\DefWarn#1%
\xdef #1{\noexpand\hyperref{}{appendix}{\appsym}{\appsym}}%
\writedef{#1\leftbracket#1}\wrlabeL{#1=#1}}
\newwrite\tfile \def\writetoca#1{}
\def\leaderfill{\leaders\hbox to 1em{\hss.\hss}\hfill}
\def\writetoc{\immediate\openout\tfile=\jobname.toc
   \def\writetoca##1{{\edef\next{\write\tfile{\noindent ##1
   \string\leaderfill {\string\hyperref{}{page}{\noexpand\number\pageno}%
                       {\noexpand\number\pageno}} \par}}\next}}}
\newread\ch@ckfile
\def\listtoc{\immediate\closeout\tfile\immediate\openin\ch@ckfile=\jobname.toc
\ifeof\ch@ckfile\message{no file \jobname.toc, no table of contents this pass}%
\else\closein\ch@ckfile\centerline{\bf Contents}\nobreak\medskip%
{\baselineskip=12pt\footnotefont\parskip=0pt\catcode`\@=11\input\jobname.toc
\catcode`\@=12\bigbreak\bigskip}\fi}
\catcode`\@=12 
%
\edef\tfontsize{
 scaled\magstep3

 }
\font\titlerm=cmr10 \tfontsize \font\titlerms=cmr7 \tfontsize
\font\titlermss=cmr5 \tfontsize \font\titlei=cmmi10 \tfontsize
\font\titleis=cmmi7 \tfontsize \font\titleiss=cmmi5 \tfontsize
\font\titlesy=cmsy10 \tfontsize \font\titlesys=cmsy7 \tfontsize
\font\titlesyss=cmsy5 \tfontsize \font\titleit=cmti10 \tfontsize
\skewchar\titlei='177 \skewchar\titleis='177 \skewchar\titleiss='177
\skewchar\titlesy='60 \skewchar\titlesys='60 \skewchar\titlesyss='60
\def\titlefont{\def\rm{\fam0\titlerm}
\textfont0=\titlerm \scriptfont0=\titlerms \scriptscriptfont0=\titlermss
\textfont1=\titlei \scriptfont1=\titleis \scriptscriptfont1=\titleiss
\textfont2=\titlesy \scriptfont2=\titlesys \scriptscriptfont2=\titlesyss
\textfont\itfam=\titleit \def\it{\fam\itfam\titleit}\rm}
 \ifx\answ\bigans\else scaled\magstep1\fi
\def\abstractfont{\tenpoint}
\def\tenpoint{\def\rm{\fam0\tenrm}
\textfont0=\tenrm \scriptfont0=\sevenrm \scriptscriptfont0=\fiverm
\textfont1=\teni  \scriptfont1=\seveni  \scriptscriptfont1=\fivei
\textfont2=\tensy \scriptfont2=\sevensy \scriptscriptfont2=\fivesy
\textfont\itfam=\tenit \def\it{\fam\itfam\tenit}\def\footnotefont{\ninepoint}%
\textfont\bffam=\tenbf \def\bf{\fam\bffam\tenbf}\def\sl{\fam\slfam\tensl}\rm}
\font\ninerm=cmr9 \font\sixrm=cmr6 \font\ninei=cmmi9 \font\sixi=cmmi6
\font\ninesy=cmsy9 \font\sixsy=cmsy6 \font\ninebf=cmbx9
\font\nineit=cmti9 \font\ninesl=cmsl9 \skewchar\ninei='177
\skewchar\sixi='177 \skewchar\ninesy='60 \skewchar\sixsy='60
\def\ninepoint{\def\rm{\fam0\ninerm}
\textfont0=\ninerm \scriptfont0=\sixrm \scriptscriptfont0=\fiverm
\textfont1=\ninei \scriptfont1=\sixi \scriptscriptfont1=\fivei
\textfont2=\ninesy \scriptfont2=\sixsy \scriptscriptfont2=\fivesy
\textfont\itfam=\ninei \def\it{\fam\itfam\nineit}\def\sl{\fam\slfam\ninesl}%
\textfont\bffam=\ninebf \def\bf{\fam\bffam\ninebf}\rm}
%
%
\def\noblackbox{\overfullrule=0pt}
\hyphenation{anom-aly anom-alies coun-ter-term coun-ter-terms}
\def\inv{^{\raise.15ex\hbox{${\scriptscriptstyle -}$}\kern-.05em 1}}

\def\Dsl{\,\raise.15ex\hbox{/}\mkern-13.5mu D} 
\def\dsl{\raise.15ex\hbox{/}\kern-.57em\partial}

 \def\Tr{{\rm Tr}}
\def\lspace{\ifx\answ\bigans{}\else\qquad\fi}
\def\lbspace{\ifx\answ\bigans{}\else\hskip-.2in\fi} 
\def\boxeqn#1{\vcenter{\vbox{\hrule\hbox{\vrule\kern3pt\vbox{\kern3pt
    \hbox{${\displaystyle #1}$}\kern3pt}\kern3pt\vrule}\hrule}}}
\def\mbox#1#2{\vcenter{\hrule \hbox{\vrule height#2in
        \kern#1in \vrule} \hrule}}  
%
 \def\CC{{\cal C}}

\def\A{{\cal A}} \def\C{{\cal C}}  
 \def\H{{\cal H}}  
\def\B{{\cal B}} \def\R{{\cal R}} \def\D{{\cal D}} \def\T{{\cal T}}

\def\darr#1{\raise1.5ex\hbox{$\leftrightarrow$}\mkern-16.5mu #1}

\def\roughly#1{\raise.3ex\hbox{$#1$\kern-.75em\lower1ex\hbox{$\sim$}}}

\input amssym.def
\input amssym.tex
\noblackbox
\def\ZZ{{\Bbb Z}}
\def\RR{{\Bbb R}}
\def\CC{{\Bbb C}}

\def\NN{{\Bbb N}}

\def\PP{{\Bbb P}}
\
\def\EE{{\Bbb E}}
\def\A{{\cal A}}
\def\B{{\cal B}}
\def\D{{\cal D}}
\def\N{{\cal N}}

\def\H{{\frak H}}
\def\h{{\cal H}}

\def\J{{\cal J}}

\def\M{{\cal M}}

\def\K{{\cal K}}
\def\O{{\cal O}}
\def\T{{\cal T}}
\def\R{{\cal R}}

\def\a{{\sl a}}
\def\b{{\sl b}}
\def\c{{\sl c}}
\def\d{{\sl d}}
\def\kk{{\sl k}}
\def\Xmin#1#2{|X^{#1}|\matrix{#2}} 
\def\Ymin#1#2{|Y^{#1}|\matrix{#2}} 
\def\ss{\scriptscriptstyle}
\def\sprod#1#2{#1\cdot #2}    
\def\spisom#1{\sprod{\tilde#1}{\tilde#1}} 
\def\spbase#1{(#1\cdot #1)} 

\def\sopra#1#2{\vbox{\baselineskip -10pt\hsize=20pt \centerline{$ #2 $}\centerline{$ #1 $}}}
\def\perm{{\cal P}}  
\def\tp#1{\,{}^{\ss t\!}#1} 
\def\smsum{\mathop{\textstyle{\sum}}\limits} 
\def\smprod{\mathop{\textstyle{\prod}}\limits} 
\def\im{\mathop{\rm Im}}

\def\sgn{\epsilon} 

\def\deltadiv{\theta_\Delta}

\def\1{\frak 1}
\def\2{\frak 2}
\def\3{\frak 3}
\def\n{\frak n}
\def\m{\frak m}
\def\imtau{\tau_2}
\def\kuno{\kappa}
\def\kenne{\kappa}
\def\kdue{\kappa}
\def\supp{\mathop{\rm supp}}
\def\vel{\vee}
\def\Sym{{\rm Sym}}
\def\diag{\rm diag}
\def\abvar{A} 
\def\ppavmod{{\cal A}} 

\def\futuref#1{\csname#1\endcsname}
\def\futurof#1{\if\csname#1\endcsname\relax ??\message{
Undefined reference \string#1
}
\else\csname#1\endcsname\fi}

\def\prendieq#1#2#3#4#5#6#7{#6}
\def\solonum#1{(\expandafter\prendieq#1)}
\def\prendith#1#2#3#4#5{#5}
\def\solonumth#1{\expandafter\prendith#1}
\global\meqno=1 \global\newcount\theno \global\theno=1
\def\newth#1#2#3{\DefWarn#1\vskip-\lastskip\medskip%
\xdef #1{\noexpand\hyperref{}{theorem}{\secsym\the\theno}{\secsym\the\theno}}%
\writedef{#1\leftbracket#1} \noindent{\bf #2
\hyperdef\hypernoname{theorem}
{\secsym\the\theno}{\secsym\the\theno}.} \global\advance\theno
by1{\sl #3}}
\def\newrem#1#2#3{\DefWarn#1\vskip-\lastskip\medskip%
\xdef #1{\noexpand\hyperref{}{theorem}{\secsym\the\theno}{\secsym\the\theno}}
\writedef{#1\leftbracket#1}%
\noindent{\bf #2 \hyperdef\hypernoname{theorem}%
{\secsym\the\theno}{\secsym\the\theno}. }%
\global\advance\theno by1{#3}}

\Title{\vbox{\baselineskip11pt }} {\vbox{ \centerline{The Singular Locus of the Theta Divisor and
}\vskip 0.5cm \centerline{Quadrics through a Canonical Curve} \vskip 1pt }}
\smallskip
\centerline{Marco Matone and Roberto Volpato}
\bigskip
\centerline{Dipartimento di Fisica ``G. Galilei'' and Istituto
Nazionale di Fisica Nucleare} \centerline{Universit\`a di Padova,
Via Marzolo, 8 -- 35131 Padova, Italy}

\vskip 0.5cm

\bigskip
\vskip 0.5cm

\noindent A section $K$ on a genus $g$ canonical curve $C$ is
identified as the key tool to prove new results on the geometry of
the singular locus $\Theta_s$ of the theta divisor. The $K$ divisor
is characterized by the condition of linear dependence of a set of
quadrics containing $C$ and naturally associated to a degree $g$
effective divisor on $C$. $K$ counts the number of intersections of
special varieties on the Jacobian torus defined in terms of
$\Theta_s$. It also identifies sections of line bundles on the
moduli space of algebraic curves, closely related to the Mumford
isomorphism, whose zero loci characterize special varieties in the
framework of the Andreotti-Mayer approach to the Schottky problem, a
result which also reproduces the only previously known case $g=4$.

This new approach, based on the combinatorics of determinantal
relations for two-fold products of holomorphic abelian
differentials, sheds light on basic structures, and leads to the
explicit expressions, in terms of theta functions, of the canonical
basis of the abelian holomorphic differentials and of the constant
defining the Mumford form. Furthermore, the metric on the moduli
space of canonical curves, induced by the Siegel metric, which is
shown to be equivalent to the Kodaira-Spencer map of the square of
the Bergman reproducing kernel, is explicitly expressed in terms of
the Riemann period matrix only, a result previously known for the
trivial cases $g=2$ and $g=3$. Finally, the induced Siegel volume
form is expressed in terms of the Mumford form.

\Date{}
%
\baselineskip13pt

\newread\instream \openin\instream=\jobname.defs
\ifeof\instream
\else \closein\instream \input \jobname.defs
\fi

\listtoc\writetoc
\writedefs
\vfill\eject

\newsec{Introduction}

In this paper we introduce a new approach leading to the characterization
of the geometry of the singular locus of the theta divisor. The starting point concerns
the determinantal relations satisfied by the two-fold products of holomorphic abelian differentials of combinatorial nature
introduced in 
\ref\MatoneBB{
  M.~Matone and R.~Volpato,
  ``Determinantal characterization of canonical curves and combinatorial theta
  identities,''
  arXiv:math.AG/0605734.
}  and here fully developed. The approach also leads to several
results concerning canonical curves, we will summarize shortly. The
key object in the investigation is the section $K$ of a suitable
line bundle on a canonical curve, which is naturally defined in
terms of the determinantal relations and encodes the geometry of the
singular locus  $\Theta_s$ of the theta divisor. It satisfies
remarkable properties and identifies sections on the moduli space of
canonical curves characterizing the geometry of $\Theta_s$. For
$g=4$ it leads to the Hessian of the theta function, evaluated at
$\Theta_s$, whose vanishing characterizes the Jacobian locus. While
this was the only previously known case of the explicit
characterization of the Jacobian, it turns out that $K$ has similar
properties which extend to arbitrary genus. In other words, the
section $K$ is the key building block to construct the higher genus
generalization of the sections on the moduli space of algebraic
curves characterizing the geometry of $\Theta_s$ that for $g=4$
reduces to the Hessian of the theta function.

In the present investigation a central role is played by Petri's
work \ref\petriuno{K.~Petri, $\rm\ddot{U}$ber die invariante
darstellung algebraischer funktionen einer
ver$\rm\ddot{a}$nderlichen, {\it Math.\ Ann.\ } {\bf 88} (1922),
242-289.} and the Andreotti-Mayer paper
\ref\andreottimayer{A.~Andreotti and A.~L.~Mayer, On period
relations for Abelian integrals on algebraic curves, {\it Ann.
Scuola Norm. Sup. Pisa} {\bf 21} (3) (1967), 189–-238.}.  Petri's
Theorem determines the ideal of canonical curves of genus $g\geq4$
by means of relations among holomorphic differentials (see also
\ref\ArbSern{ E.~Arbarello and E.~Sernesi, Petri's approach to the
study of the ideal associated to a special divisor, {\it Invent.
Math. } {\bf 49} (1978), 99-119.}\ref\ArbHarr{E.~Arbarello and
J.~Harris, Canonical curves and quadrics of rank $4$, {\it Compos.
Math.} {\bf 43} (1981), 145-179.}). As emphasized by Mumford,
Petri's relations are fundamental and should have basic applications
(pg.241 of \ref\mumfordd{D. Mumford, {\it The Red Book of Varieties
and Schemes}, Springer Lecture Notes 1358 (1999).}). Petri's
construction relies on the definition of a basis of holomorphic
abelian differentials which is, in some sense, ``dual'' to a
$g$-tuple of distinct points on the curve. Schematically, for a
canonical curve $C$ of genus $g\ge 3$,
$$\hbox{distinct points }p_1,\ldots,p_g\in C\qquad\longleftrightarrow\qquad\hbox{basis
}\sigma_1,\ldots,\sigma_g\hbox{ of }H^0(K_C)\ ,$$
where $p_1,\ldots,p_g\in C$ and $\{\sigma_i\}_{1\le i\le g}$ is
basis of $H^0(K_C)$, with $K_C$ the canonical line bundle on $C$,
such that $\sigma_i(p_j)=0$ if and only if $i\neq j$, $1\le i,j\le
g$. Such a condition determines the basis $\{\sigma_i\}_{1\le i\le
g}$ up to a non-singular diagonal transformation.

 Max Noether's Theorem
assures that, for non-hyperelliptic Riemann surfaces, the space of
holomorphic $n$-differentials, with $n>1$, is generated by $n$-fold
products of abelian differentials. Under general conditions on the
dual $g$-tuple of points, Petri's basis can be used to explicitly
construct bases of holomorphic $n$-differentials. Such a procedure
leads to the Enriques-Babbage-Petri Theorem, according to which the
ideal of a (non-singular) canonical curve, with few exceptions, is
generated by quadrics. The exceptions are the trigonal curves and
the smooth plane quintics, which can be obtained as complete
intersections of quadrics and cubics.

The idea of constructing bases of holomorphic $n$-differentials in terms of $n$-fold products of
abelian differentials is very powerful; in some sense, it may be used to mimic analogous constructions in the
hyperelliptic case, which, for example, lead to explicit results for the Mumford form at genus $2$.
In particular, it would be useful to connect {\it determinants} of $n$-differentials to
{\it determinants} of $1$-differentials. This amounts to a non-trivial combinatorial problem which
is of interest on its own, and is solved in section \futuref{detlemmas}. Determinants of holomorphic
$n$-differentials can be used to the definition of some Petri-like bases of $n$-differentials,
with also provides a standard normalization procedure which is lacking in the original Petri
approach. This is crucial in order to ensure the modular invariance of such Petri-like
bases, once one chooses a marking for the Riemann surface and define such bases in terms of Riemann
theta functions.

\vskip 3pt

As the most direct application of the construction, we will
determine a necessary and sufficient condition for suitable sets of
two-fold and three-fold products of holomorphic abelian
differentials to be linear independent generators for the spaces of,
respectively, holomorphic quadratic and cubic differentials. The
analysis of such a condition is developed in several steps in
section \futuref{primecostr}, and results in Theorem
\futuref{quadgen}. In particular, in the case of two-fold products,
such a condition is deeply related to the geometry of the singular
locus of the theta divisor associated to the Riemann surface. This
leads to the definition of two sections of suitable line bundles on
the curve, we will denote by $H$ and $K$, depending on the points
``dual'' to the basis of abelian differentials, which vanish exactly
when the two-fold products fail to be a basis of holomorphic
$2$-differentials. The divisor of such sections (in particular, of
the section $K$) represents the key for a deeper understanding of
the geometry of the singular locus $\Theta_s$ of the theta divisor;
Theorem \futuref{zeriK} is a remarkable application of such an
approach and one of the main results of this paper.

\vskip 3pt

The two-fold products of holomorphic $1$-differentials also
correspond to the generators of the symmetric product
$\Sym^2(H^0(K_C))$, which is a $M$-dimensional vector space, with
$M:=g(g+1)/2$. On the other hand, the space $H^0(K_C^2)$ is
$N$-dimensional, with $N=3(g-1)$, so that the natural homomorphism
$\psi\colon \Sym^2(H^0(K_C))\to H^0(K_C^2)$ has a kernel of
dimension $M-N=(g-2)(g-3)/2$. A basis of $\ker\psi$ can be described
in two different ways:
\smallskip
\item{{\it a})} A set of $M-N$ linearly independent relations among holomorphic quadratic differentials
$$\sum_{i,j=1}^gC^{\eta}_{k,ij}\eta_i\eta_j=0\ ,\qquad k=N+1,\ldots,M\ ,$$
where $\{\eta_i\}_{1\le i\le g}$ is a basis of $H^0(K_C)$ and $\{C_k^\eta\}_{N<k\le M}$ is a set of
linearly independent elements of
$\PP(\Sym^2H^0(K_C)^*)$.
\smallskip
\item{{\it b})} The choice of a basis $\{\eta_i\}_{1\le i\le g}$ of $H^0(K_C)$ determines a canonical
embedding of a non-hyperelliptic Riemann surface as a canonical curve into $\PP_{g-1}$ via $C\ni p\to
(\eta_1(p),\ldots,\eta_2(p))\in \PP_{g-1}$; in other words, each $\eta_i$, $1\le i\le g$, is
identified with a projective coordinate $X_i$ of $\PP_{g-1}$. Under such an identification, a basis
of $\ker\psi$ corresponds to a set of generators
$$\sum_{i,j=1}^gC^{\eta}_{k,ij}X_iX_j=0\ ,\qquad k=N+1,\ldots,M\ ,$$
of $I_2$, the ideal of quadric passing through the
curve $C$.
\smallskip
As mentioned above, by the Enriques-Babbage-Petri Theorem, in most
cases, canonical curves are complete intersections of quadrics in
$\PP_{g-1}$, so that the knowledge of the generators of $I_2$
completely determines the curve.

However, this is not the only motivation for the interest in $I_2$.
A classical result due to Riemann shows that each double point $e$
on the singular locus of the theta divisor connected to the Riemann
surface (with marking) corresponds to a relation among holomorphic
quadratic differentials
\eqn\RiemRel{\sum_{i,j=1}^g\theta_{ij}(e)\omega_i\omega_j=0\ ,} with
respect to the canonical basis $\{\omega_i\}_{1\le i\le g}$ of
$H^0(K_C)$. Equivalently, such a relation represents an element in
$I_2$, which can be easily proved to be a quadric of rank $r\le 4$.
The interest in the study of the ideal $I_2$ was renewed by the work
of Andreotti and Mayer \andreottimayer. The motivation for their analysis is strictly related to
the Schottky problem, which amounts to the definition of necessary
and sufficient conditions for a principally polarized abelian
variety (ppav) to correspond to the Jacobian torus of a Riemann
surface. In their beautiful construction, Andreotti and Mayer
proposed to characterize the Jacobian locus (denoted by $\J_g$ for
the non-hyperelliptic Riemann surfaces and $\h_g$ for hyperelliptic
ones) inside the moduli space $\ppavmod_g$ of ppav's, through the
dimension of the singular locus of the theta divisor. More
precisely, they showed that $\J_g$ (resp., $\h_g$) is a component of
$\N_{g-4}\subset\ppavmod_g$ (resp., $\N_{g-3}\subset \ppavmod_g$),
where $\N_k$ is the locus of the ppav's such that $\dim\Theta_s\ge
k$. As a crucial point in such a construction, they proved that, for
each trigonal curve, the ideal $I_2$ is generated by relations in
the form \RiemRel, as $e$ varies in $\Theta_s$. Such a result has
received many remarkable generalizations, among which at least two
deserve citation: Arbarello and Harris \ArbHarr\ proved that the
relations \RiemRel\ generate $I_2$ for $g\le 6$ and that, for all
$g$, they generate all the quadrics of rank $\le 4$; finally, Green
\ref\Green{M.~Green, Quadrics of rank four in the ideal of a
canonical curve, {\it Invent. Math.} {\bf 75} (1984), no. 1,
85-104.} proved that such relations generate $I_2$ for all genera,
so that, as a consequence, $I_2$ can be generated by quadrics of
rank $\le 4$ only.

In this respect, some fundamental problems are still unsolved. First
of all, a procedure to determine a finite set of (possibly linearly
independent) generators for $I_2$ in the form \RiemRel, is still
lacking. More generally, no explicit expression of the coefficients
$C^{\eta}_{k,ij}$ in terms of the Riemann period matrices, that is,
in terms of theta constants or modular forms, is actually known. In
view of the Andreotti and Mayer construction, such a result may
represent a key step toward an explicit solution of the Schottky
problem, i.e. as a characterization of the Jacobian locus by
algebraic conditions on suitable modular forms, similar to the
Schottky relation at genus $4$ \ref\FarkasRauch{H.~M.~Farkas and H.~E.~Rauch,
Period relations of Schottky type on Riemann surfaces,
{\it Ann. of Math.} (2) {\bf 92} (1970), 434-461.}.

\vskip 6pt

\subsec{Main results}

In this paper we introduce a new approach to the above problems and the definition of some powerful tools for their
analysis. The main results of the present investigation concern:
\smallskip
\item{--} Several alternative expressions for $M-N=(g-2)(g-3)/2$ linearly independent coefficients
$C^{\eta}_{k,ij}$, $k=N+1,\ldots,M$, depending both on
the Riemann period matrix and on $g-2$ points $p_3,\ldots,p_g\in C$ (Theorem \futuref{thleX} and
Corollary \futuref{threla}).
\smallskip
\item{--} The correspondence between each relation and a combinatorial identity among Riemann theta functions
evaluated on points of the Riemann surface (Theorem \futuref{ththetarel}).
\smallskip
\item{--} The expansion of the relation \RiemRel, for all $e\in\Theta_s$, with respect to the
relations given by the coefficients $C^{\eta}_{k,ij}$ (Lemma \futuref{FaKraRel}).
\smallskip
\item{--} The connection between a zero of the section $K$, which corresponds to a point
$e\in\Theta_s$, and the relation among holomorphic quadratic
differentials given by \RiemRel\ at $e$. More precisely, it turns out that $K=0$ is a
sufficient condition for the existence of a linear relation among a
suitable set of $N$ holomorphic $2$-differentials; Theorem
\futuref{voila} identifies such a relation with Eq.\RiemRel, with
respect to the corresponding point $e\in\Theta_s$.
\smallskip
\item{--} A long standing problem in the study of Riemann surfaces, which arises for example in investigating the
Schottky problem or in constructing modular forms, is that some
basic constants, such as the Mumford form, have an expression that
needs the use of points on $C$. On the other hand, such quantities
should be expressed in terms of {\it Thetanullwerte} and, via the
higher genus generalization of the Jacobi identities, of the related
{\it Jacobian Nullwerte}. Here we introduce a new general strategy
which is based on the idea of identifying the divisors with the one
defining spin structures. In doing this one has to consider several
intermediate problems. For example, the Mumford form, as many other
basic quantities, is expressed in terms of determinants of
holomorphic differentials and, in particular, of the determinant of
the canonical basis $\{\omega_i\}_{i\in I_g}$ of $H^0(K_C)$. Looking
for expressions where only the {\it Thetanullwerte} and the {\it
Jacobian Nullwerte} appear, requires expressing $\det\omega_i(p_j)$
in terms of theta functions only. Furthermore, another problem
concerning $\det\omega_i(p_j)$ is to consider the $g$-points
$p_1,\ldots,p_g$ as defining spin structures associated to divisors
of degree $g-1$. As we will see, this question is strictly related
to the problem of expressing $\det\omega_i(p_j)$ without the use of
the $g/2$-differential $\sigma$ and of the constant $\kuno[\omega]$.
It turns out that there exists a natural solution
leading to the explicit expression of basic quantities in terms of
divisors defining spin structures. In particular, we first
explicitly express the abelian holomorphic differentials in terms of
theta functions only, a result previously known in the trivial
elliptic case. This also implies the expression for the basic
constant $\kuno[\omega]$ (that, as we will see, needs an important
refinement leading to our $\kuno_\nu[\omega]$, a problem also related to
the definition of the $g/2$-differential $\sigma$, the carrier of
the holomorphic and gravitational Liouville anomalies, well-known in
string theories) usually defined in terms of $\det\omega_i(p_j)$, corresponding to the
main building block of the Mumford form. We then
will obtain the explicit expression for products of
$\kuno_{\nu_k}[\omega]$ in terms of theta functions with spin
structures whose arguments involve the difference of points
belonging to the divisors of such spin structures. The present
approach seems having interesting consequences which extend to
branches related to the theory of Riemann surfaces, including the
theory of Siegel modular forms.
\smallskip
\item{--}
The metric $ds^2_{|\hat\M_g}$ on the moduli space $\hat\M_g$ of
genus $g$ canonical curves induced by the Siegel metric, is
expressed in terms of the Riemann period matrix, a result previously
known for the trivial cases $g=2$ and $g=3$. It turns out that
such a metric is equivalent to the Kodaira-Spencer map of the square
of the Bergman reproducing kernel. Furthermore, the induced Siegel
volume form is expressed in terms of the Mumford form.
\smallskip
\item{--} The combinatorial Lemma \futuref{thcombi}, regarding the expansion of a determinant of holomorphic $2$-differentials in terms of
determinants holomorphic abelian differentials, is related to the Mumford
isomorphism $\lambda_2\cong\lambda_1^{13}$, where $\lambda_i$, $i\ge 1$, is the determinant line bundle on the
moduli space $\M_g$, with fiber $\bigwedge^{top}H^0(K^i_C)$ at the point representing $C$. In particular, the section $K$ plays again a key role.
For $g=3$, $K$
is in fact a constant on $C$, and is proportional to $\Psi_9$, a square root of the modular form
$\Psi_{18}$ of weight $18$. For $g=4$, by taking a suitable product of the sections $K$'s, the
dependence on the points of $C$ can be eliminated, and the resulting quantity is proportional to the
Hessian of the theta function $\det_{ij}\theta_{ij}(e)$, where $e$ is one of the two (generally
distinct) points of $\Theta_s$.
\smallskip
\item{--} Remarkably, a conjecture by H.~M.~Farkas \ref\Farkas{H.~M.~Farkas,
Vanishing theta nulls and Jacobians, {\it The geometry of Riemann
surfaces and abelian varieties}, 37-53, Contemp. Math. 397, 2006.}
has been recently proved \ref\grushsalv{S.~Grushevsky and R.~Salvati
Manni, Jacobians with a vanishing theta-null in genus 4, arXiv:
math.AG/0605160.}, according to which, the vanishing of such a
Hessian characterizes the elements in $\J_4$ among the ppav's with a
vanishing theta-null. Note that all such ppav's are elements of
$\N_0$ (i.e. the locus of ppav's with non-empty $\Theta_s$);
according to \ref\Beauville{A.~Beauville, Prym varieties and the
Schottky problem, {\it Invent. Math.}  {\bf 41}  (1977), no. 2,
149-196.} it turns out that for $g=4$, $\N_0$ is given by the union
of such a locus and $\J_4$ and in this case such a result, together
with the Andreotti-Mayer criterium, completely characterizes $\J_4$.
Sections with well-defined properties under modular transformations
can similarly be constructed in terms of the section $K$ for $g\ge
5$; in this sense, $K$ appears as a powerful tool for the study of
$\Theta_s$ and, more generally, of the geometry of the moduli space.
We provide explicit examples for $g=5$ and for any even genus.

The above construction was initiated in 
\ref\MatoneBX{
  M.~Matone and R.~Volpato,
  Linear relations among holomorphic quadratic differentials and induced
  Siegel's metric on $\M_g$,
  arXiv:math.AG/0506550.
}\MatoneBB, in relation to the investigation
in \ref\MatoneVM{
  M.~Matone and R.~Volpato,
  Higher genus superstring amplitudes from the geometry of moduli spaces,
{\it  Nucl.\ Phys.\ } B {\bf 732} (2006), 321-340
  [arXiv:hep-th/0506231].
  } concerning the $g>2$ generalization of the remarkable
D'Hoker and Phong formula 
\ref\DHokerJC{E.~D'Hoker and D.~H.~Phong,
  Two-loop superstrings. VI: Non-renormalization theorems and the 4-point
  function,
{\it Nucl.\ Phys.\ } B {\bf 715} (2005), 3-90
  [arXiv:hep-th/0501197].
} for the four point superstring amplitude.

Before reporting on the plan of the paper we remind some other basic facts about the Schottky problem.
According to the Novikov's conjecture, a
indecomposable ppav is the Jacobian of a genus $g$ curve if
and only if there exist vectors $U\neq0,V,W\in \CC^g$ such that
$u(x,y,t)=2\partial_x^2\log\theta(Ux+Vy+Wt+z_0,Z)$, satisfies the
Kadomtsev-Petviashvili (KP) equation
$3u_{yy}=(4u_t+6uu_x-u_{xxx})_x$. Relevant progresses on such a
conjecture are due, among the others, to Krichever
\ref\krichever{I.~M.~Krichever, Integration of non-linear equations
by methods of algebraic geometry, {\it Funct. Anal.\ Appl.\ } {\bf
11} (1) (1977), 12-26.}, Dubrovin \ref\dubrovin{B.~Dubrovin, Theta
functions and non-linear equations, {\it Russ.\ Math.\ Surv.\ } {\bf
36} (1981), 11-92.} and Mulase \ref\mulase{M.~Mulase, Cohomological
structure in soliton equations and Jacobian varieties, {\it J.\
Diff.\ Geom.\ } {\bf 19} (1984), 403-430.}. Its proof is due to
Shiota \ref\shiota{T.~Shiota, Characterization of Jacobian varieties
in terms of soliton equations, {\it Invent.\ Math.\ } {\bf 83}
(1986), 333-382.}. A basic step in such a proof concerned the
existence of the $\tau$-function as a global holomorphic function in
the $\{t_i\}$, as clarified by Arbarello and De Concini in
\ref\Arbar{E.~Arbarello and C.~De Concini, Another proof of a
conjecture of S.P. Novikov on periods of abelian integrals on
Riemann surfaces, {\it Duke\ Math.\ Journal\ } {\bf 54} (1987),
163-178.}, where it was shown that only a subset of the KP hierarchy
is needed. Their identification of such a subset is based on basic
results by Gunning \ref\gunning{R. Gunning, Some curves in abelian
varieties, {\it Invent.\ Math.\ } {\bf 66} (1982), 377-389.} and
Welters \ref\welters{G.~E.~Welters, A characterization of
non-hyperelliptic Jacobi varieties, {\it Invent.\ Math.\ } {\bf 74}
(1983), 437-440.} \ref\Welterfg{G.~E.~Welters, A criterion for
Jacobi varieties, {\it Ann.\ Math.\ } {\bf 120} (1984), 497-504.},
characterizing the Jacobians by trisecants (see also
\ref\donagi{R.~Donagi, The Schottky problem, In {\it Theory of
moduli} (Montecatini Terme, 1985), 84-137, Springer Lecture Notes
1337.} and \ref\ArbDeConc{E.~Arbarello and C.~De Concini, On a set
of equations characterizing Riemann matrices, {\it Ann. Math.} {\bf
120} (1984), 119-140.}).

The Schottky problem is still under active investigation, see for
example \ref\Marini{G.~Marini, A geometrical proof of Shiota's
Theorem on a conjecture of S.P. Novikov, {\it Compos.\ Math.\ } {\bf
111} (1998), 305-322.}\ref\polishchuk{A.~Polishchuk, {\it Abelian
varieties, theta functions and the Fourier transform}, Cambridge
Univ. Press, 2003.}\ref\grushevsky{S.~Grushevsky, The degree of the
Jacobian locus and the Schottky problem,
math.AG/0403009.}\ref\Sebastian{S.~Casalaina-Martin, Prym varieties
and the Schottky problem for cubic threefolds, arXiv:
math.AG/0605666.}\ref\KrGr{I.~M.~Krichever and S.~Grushevsky,
Integrable discrete Schrodinger equations and a characterization of Prym varieties by a pair of quadrisecants,
    arXiv:0705.2829v1 [math.AG].}
\ref\CilibderGe{C.~Ciliberto and G.~van der Geer, Andreotti-Mayer loci and the Schottky problem,
        arXiv: math/0701353v1 [math.AG].}
for further developments. In particular,
Arbarello, Krichever and Marini proved that the Jacobians can be
characterized in terms of only the first of the auxiliary linear
equations of the KP equation \ref\ArKrMar{E.~Arbarello,
I.~M.~Krichever and G.~Marini, Characterizing Jacobians via flexes
of the Kummer variety, {\it Math. Res. Lett.} {\bf 13} (1) (2006), 109-123.
}
\ref\KricheverET{
  I.~M.~Krichever,
  Integrable linear equations and the Riemann-Schottky problem,
  {\it Algebraic geometry and number theory}, 497-514, Progr. Math.
  {253}, Birkh\"auser Boston, Boston, 2006.
}. Very recently Krichever
 \ref\Kricheverc{I.~M.~Krichever,
Characterizing Jacobians via trisecants of the Kummer variety,
  arXiv: math.AG/0605625.
}  has proved the conjectures by Welters \Welterfg.

\vskip 6pt

\subsec{Plan of the paper}

\item{--} In section \futuref{notation} we first present some
useful result on tensor products of vector spaces. Then, in subsection
\futuref{detlemmas}, two combinatorial lemmas are proved. Consider a
finite set of functions from an arbitrary set $S$ to a commutative
field, and consider the two-fold products $ff$ of such functions.
Under certain conditions, we will prove that the determinant of such
two-fold products $\det ff_i(x_j)$ evaluated on a suitable number of
elements $x_j\in S$, can be obtained by skew-symmetrization of a
product of determinants of the form $\det f_i(x_j)$. Such a result,
which reveals new structures concerning determinantal properties,
whose interest is not restricted to the theory of algebraic curves,
plays a central role in the present investigation.
\smallskip
\item{--} In section \futuref{Riemanndef}, after introducing some facts on theta functions and the
Jacobian torus we derive relations among higher order theta derivatives and holomorphic
differentials which will be used later.
\smallskip
\item{--} In section \futuref{primecostr}, Petri-like bases of holomorphic $n$-differentials are
defined. Then, we consider the construction of sets of quadratic and
cubic differentials in terms of two- and three-fold products of
holomorphic abelian differential and give necessary and sufficient
conditions for such sets to be bases of $H^0(K_C^2)$ and
$H^0(K_C^3)$. Such considerations will lead to the definitions of
two sections, denoted by $H$ and $K$ (Eqs.\futuref{ilrappo} and
\futuref{iltrap}). The section $K$ will then lead to the
characterization of the special locus of the theta divisor (Theorem
\futuref{zeriK}).
\smallskip
\item{--}{The bases defined in section \futuref{primecostr} are the basic tool for the derivation,
developed at the beginning of section \futuref{deterrel}, of $M-N$
linear independent relations among holomorphic quadratic
differentials obtained as two-fold products of holomorphic abelian
differentials. The coefficients appearing in such relations depend
on the points $p_1,\ldots,p_g$ determining the basis of $H^0(K_C)$
by proposition \futuref{newbasis}. After considering the consistency
conditions on such coefficients, it is shown that such relations
correspond to new non-trivial identities among theta functions and
prime forms evaluated on arbitrary points of the Riemann surface
(Theorems \futuref{main} and \futuref{ththetarel}). It should be
emphasized that such identities follow by the combinatorial Lemmas
\futuref{thcombi} and \futuref{thcombvi}, and by the Fay's trisecant
identity, so that no analogous result is expected for theta
functions defined on arbitrary ppav's (i.e., outside the Jacobian
locus). Furthermore, the coefficients are shown to be related to the
section $K$, so leading to two main consequences. First, the
dependence on a pair of points $p_1,p_2$ can be eliminated
(Corollary \futuref{threla}). Then, the connection among the
relations obtained and Eq.\RiemRel\ is clarified: each point on
$\Theta_s$ corresponds, in some sense, to a zero of the section $K$,
that is, to the failure for the associated set of holomorphic
quadratic differentials, built in proposition \futuref{thlev}, to be
a basis of $H^0(K_C^2)$. The linear relation arising among such
holomorphic quadratic differentials just corresponds to Eq.\RiemRel.

{The case of the curves of genus four is somewhat special, since a unique relation exists. This
subject is considered in a separate subsection, where it is shown that in this case $K$
is deeply related to the Hessian of the theta function evaluated on a singular point.
Furthermore, an expression for the coefficients completely independent of the points
$p_1,\ldots,p_4$ is obtained.
Finally, in the last
subsection, the techniques developed in this section and the bases defined in proposition \futuref{thtrebase} are
applied to the derivation of the relations among holomorphic cubic differentials.
In view of the Enriques-Babbage-Petri Theorem, this is a very interesting subject which deserves
further study along the lines described in this article.}}
\smallskip
\item{--}  In section \futuref{thetkuno}
it is shown that the use of the distinguished bases for the
holomorphic differentials leads to a straightforward derivation of
the Fay's trisecant identity. We then show that the canonical basis
of the abelian holomorphic differentials can be expressed purely in
terms of theta functions. This also implies an expression for the
basic building block of the Mumford form given in terms of theta
functions only.
\smallskip
\item{--} In section \futuref{Siegsec} we first provide the explicit expression of the metric $ds^2_{|\hat\M_g}$
on the moduli space $\hat\M_g$ of genus $g$ canonical curves induced by the Siegel metric, which was known
only in trivial cases $g=2$ and $g=3$. It turns out that $ds^2_{|\hat\M_g}$ can be also expressed
as the Kodaira-Spencer map of the square of the Bergman reproducing kernel (times $4\pi^2$). By Wirtinger Theorem the explicit expression
for the volume form on $\hat\M_g$ is also obtained. Furthermore,
a notable relation satisfied by the determinant of powers of the Bergman reproducing kernel is proved.
Such results are a natural consequence of the present approach,
which also uses, as for the derivation of $ds^2_{|\hat\M_g}$,
the isomorphisms introduced in section \notation.
\smallskip
\item{--} In section \futuref{Mumfsec} we first use the construction of section
\futuref{secdivspin} to derive an expression for the Mumford form
which does not involve any determinant of holomorphic
$1$-differentials. Furthermore, we express the volume form
$d\nu_{|_{\hat\M_g}}$, induced by the Siegel metric, in terms of the
Mumford form. Next, by means of the Mumford isomorphism we
investigate the modular properties of $K(p_3,\ldots,p_g)$ in order
to construct sections of bundles on $\M_g$. For $g=2$ and $g=3$ such
sections reproduce the building blocks for the Mumford form. For
$g=4$, a modular form on the Jacobian locus is obtained, which is
proportional to the Hessian of the theta function evaluated on
$\Theta_s$. This is a remarkable result in view of
\Farkas\grushsalv, where it is shown that the vanishing of such a
Hessian on the Andreotti-Mayer locus $\N_0=\J_4\cup \theta_{null}$,
where $\theta_{null}\subset\ppavmod_4$ is the locus of the ppav's
with a vanishing theta-null, characterizes the intersection
$\J_4\cap \theta_{null}$. This indicates that the sections on $\M_g$
built in terms of $K$ may be considered as generalizations to $g>4$
of such a Hessian, thus providing a tool for the analysis of the
geometry of $\Theta_s$ and of the Andreotti-Mayer locus $\N_{g-4}$.
We explicitly construct such sections for even genus and for the
case $g=5$.

\vskip 6pt

\newsec{Two combinatorial lemmas on determinants of symmetric products}\seclab\notation

Determinants of holomorphic quadratic differentials play a crucial
role in our construction. In particular, in the following sections,
we will construct bases of $H^0(K_C^2)$ in terms of two-fold
products of holomorphic abelian differentials. In this section, we
will consider the purely combinatorial problem concerning the
determinants of a basis of a two-fold symmetric product of a finite
dimensional space of functions. We first introduce a very useful
notation for symmetric tensor products of vector space, which we
will adopt all along the paper; then we derive two lemmas on
determinants which are of interest on their own.

\vskip 6pt

\subsec{Identities induced by the isomorphism
$\CC^{M_n}\leftrightarrow \Sym^n\CC^g$}

\newrem\primedef{Definition}{For each $n\in\ZZ_{>0}$, set $I_n:=\{1,\ldots,n\}$ and let $\perm_n$ denote the group of permutations
of $n$ elements.}

\vskip 6pt

Let $V$ be a $g$-dimensional vector space and let
$$M_n:=\Bigl({g+n-1\atop
n}\Bigr)\ ,
$$ be the dimension of the $n$-fold symmetrized tensor product $\Sym^n V$. We denote by
$$\Sym^nV\ni\eta_1\cdot\eta_2\cdots\eta_n:=\sum_{s\in\perm_n}\eta_{s_1}\otimes\eta_{s_2}\otimes\ldots\otimes
\eta_{s_n}\ ,$$
the symmetrized tensor product of an $n$-tuple $(\eta_1,\ldots,\eta_n)$ of elements of $V$.

\vskip 6pt

Fix a surjection $m:I_g\times I_g\rightarrow I_M$,
$M:=M_2=g(g+1)/2$, such that \eqn\suriezione{m(i,j)=m(j,i)\ ,}
$i,j\in I_g$. Such a surjection corresponds to an isomorphism
$\CC^M\rightarrow\Sym^2\CC^g$ with $\tilde e_{m(i,j)} \mapsto
\sprod{e_i}{e_j}$.

A useful choice for such an isomorphism is considered in the following definition.

\newrem\definizione{Definition}{
Let $A:\CC^M\rightarrow{\rm Sym}^2\CC^g$, $M\equiv M_2$, be the isomorphism
$A(\tilde e_i):=\sprod{e_{\1_i}}{e_{\2_i}}$, with $\{\tilde
e_i\}_{i\in I_M}$ the canonical basis of $\CC^M$ and
$$(\1_i,\2_i):=\left\{\vcenter{\vbox{\halign{\strut\hskip 6pt $ # $ \hfil & \hskip 2cm$ # $ \hfil\cr (i,i)\ , &1\le i\le g\ ,\cr
(1,i-g+1)\ ,&g+1\le i\le 2g-1\ ,\cr(2,i-2g+3)\ , &2g\le i\le 3g-3\
,\cr \hfill\vdots \hfill& \hfill\vdots\hfill\cr (g-1,g)\
,&i=g(g+1)/2\ ,\cr}}}\right.$$ so that $\1_i\2_i$ is the $i$-th
element in the $M$-tuple $(11,22,\ldots,gg,12,\ldots,1g,23,\ldots)$.
Similarly, let $\{\tilde e_i\}_{i\in I_{M_3}}$ be the canonical basis of
$\CC^{M_3}$, and fix an isomorphism $A:\CC^{M_3}\to {\rm Sym}^3\CC^g$, $M_3:=g(g+1)(g+2)/6$,
with $A(\tilde e_i):=(e_{\1_i},e_{\2_i},e_{\3_i})_S$,
whose first $6g-8$ elements are
$$(\1_i,\2_i,\3_i):=\left\{\vcenter{\vbox{\halign{\strut\hskip 6pt $ # $ \hfil & \hskip 2cm$ # $ \hfil\cr
(i,i,i)\ , &1\le i\le g\ ,\cr (1,1,i-g+2)\ ,&g+1\le i\le 2g-2\ ,\cr
(2,2,i-2g+4)\ , &2g-1\le i\le 3g-4\ ,\cr (1,2,i-3g-4)\ , &3g-3\le
i\le 4g-4\ ,\cr (1,i-4g+6,i-4g+6)\ ,& 4g-3\le i\le 5g-6\ ,\cr
(2,i-5g+8,i-5g+8)\ ,& 5g-5\le i\le 6g-8\ .\cr}}}\right.$$}

\vskip 6pt

As we will see, we do not need the explicit expression of $A(\tilde e_i)$ for
$6g-8<i\le M_3$.
In general, one can define an isomorphism $A:\CC^{M_n}\to {\rm Sym}^n\CC^g$, with $A(\tilde
e_i):=(e_{\1_i},\ldots,e_{\n_i})$, by fixing the $n$-tuples $(\1_i,\ldots,\n_i)$, $i\in I_{M_n}$, in
such a way that $\1_i\le \2_i\le \ldots\le \n_i$.

For each vector $u:=\tp(u_1,\ldots,u_g)\in\CC^g$ and matrix $B\in
M_g(\CC)$, set
$$\underbrace{u\cdots u_i}_{n\hbox{ times}}:=\prod_{\m\in\{\1,\ldots,\n\}}u_{\m_i}\ ,\qquad
(\underbrace{B\cdots B}_{n\hbox{
times}})_{ij}:=\sum_{s\in\perm_n}\prod_{\m\in\{\1,\ldots,\n\}}B_{\m_is(\m)_j}\
,$$ $i,j\in I_{M_n}$, where the product is the standard one in
$\CC$. In particular, let us define
$$\chi_i\equiv\chi_i^{(n)}:=\prod_{k=1}^g\biggl(\sum_{\m\in\{\1,\ldots,\n\}}\delta_{k\m_i}\biggr)!=(\delta\cdots\delta)_{ii}\ ,$$ $i\in I_{M_n}$,
(we will not write the superscript $(n)$ when it is clear from the
context) where $\delta$ denotes the identity matrix, so that, for
example,
$$\chi^{(2)}_i=1+\delta_{\1_i\2_i}\
,\qquad\quad
\chi_i^{(3)}=(1+\delta_{\1_i\2_i}+\delta_{\2_i\3_i})(1+\delta_{\1_i\3_i})\
.$$

\vskip 6pt

\noindent Such a single indexing satisfies basic identities,
repeatedly used in the following.

\vskip 6pt

\newth\idddd{Lemma}{Let $V$ be a vector space and $f$ an arbitrary function $f:I^n_g\to
V$, where $I^n_g:=I_g\times\ldots\times I_g$ ($n$ times). Then, the
following identity holds
\eqn\identitt{\sum_{i_1,\ldots,i_n=1}^gf(i_1,\ldots,i_n)=\sum_{i=1}^{M_n}\chi_i^{-1}\sum_{s\in\perm_n}f(s(\1)_i,\ldots,s(\n)_i)\
,} that, for $f$ completely symmetric, reduces to
\eqn\identittddd{\sum_{i_1,\ldots,i_n=1}^gf(i_1,\ldots,i_n)=n!\sum_{i=1}^{M_n}\chi_i^{-1}f(\1_i,\ldots,\n_i)\
.}}

\vskip 6pt

\noindent{\bf Proof.} Use
$$\sum_{i_1,\ldots,i_n=1}^gf(i_1,\ldots,i_n)=\sum_{i_n\ge\ldots
\ge i_1=1}^g\sum_{s\in\perm_n}{f(i_{s_1},\ldots,i_{s_n})\over
\prod_{k=1}^g(\sum_{m=1}^n\delta_{ki_m})! }\ .$$ \hfill$\square$

\vskip 6pt

\noindent Note that $u^{\otimes n}\equiv u\otimes\ldots\otimes u$ is
an element of $\Sym^n\CC^g\cong\CC^{M_n}$, for each $u\in\CC^g$. By
\identitt, the following identities are easily verified
$$u^{\otimes n}\cong \sum_{i=1}^{M_n}\chi_i^{-1}u\cdots u_i \tilde e_i\ ,\quad (Bu)^{\otimes
n}\cong\sum_{i,j=1}^{M_n}\chi_i^{-1}\chi_j^{-1}(B\cdots
B)_{ij}u\cdots u_j\tilde e_i\ ,
$$
where $\CC^{M_n}\ni\tilde e_i\cong e\cdots e_i\in\Sym^n\CC^g$, $i\in
I_{M_n}$. Furthermore,
\eqn\formprod{\sum_{j=1}^{M_n}\chi_j^{-1}(B\cdots B)_{ij}(C\cdots
C)_{jk}=((BC)\cdots (BC))_{ik}\ ,} where $B,C$ are arbitrary
$g\times g$ matrices. This identity yields, for any non-singular $B$
\eqn\forminv{\sum_{j=1}^{M_n}\chi_j^{-1}\chi_{k}^{-1}(B\cdots
B)_{ij}(B^{-1}\cdots
B^{-1})_{jk}=(\delta\cdots\delta)_{ik}\chi_k^{-1}=\delta_{ik}\ ,}
and then \eqn\formdet{\det\nolimits_{ij}\bigl((B\cdots
B)_{ij}\chi_j^{-1}\bigr)\det\nolimits_{ij}\bigl((B^{-1}\cdots
B^{-1})_{ij}\chi_j^{-1}\bigr)=1\ .} Also observe that
\eqn\formsimprod{\prod_{i=1}^{M_n}u\cdots
u_i=\prod_{k=1}^gu_k^{{n\over g}M_n}\ ,}  where the product and the
exponentiation are the standard ones among complex numbers; in
particular,
$$\prod_{i=1}^{M}uu_i=\prod_{k=1}^gu_k^{g+1}\ .$$
In the following we will denote the minors of $(B\cdots B)$ by
$$|B\cdots B|^{i_1\ldots i_m}_{j_1\ldots j_m}:=\det_{{i\in {i_1,\ldots,i_m}\atop j\in
{j_1,\ldots,j_m}}}(B\cdots B)_{ij}\ ,$$
$i_1,\ldots,i_m,j_1,\ldots,j_m\in I_{M_n}$, with $m\in I_{M_n}$.

\vskip 6pt

\newrem\definsiemi{Definition}{Fix $g,n\in \ZZ_{>0}$. Set
$$I_{M_n}\supset I^{\diag}_n:=\{i\in I_{M_n}\mid \1_i=\2_i=\ldots=\n_i\}\ .$$
Fix $l<g$ and $a,a_1,\ldots,a_l\in I_g$ and define the following
subsets of $I_{M_n}$
$$\eqalign{&I^a_n:=\{i\in I_{M_n}\mid \1_i= a\vel
\2_i=a\vel\ldots\vel\n_i=a\}\ ,\cr &I^{a_1\ldots
a_l}_n:=\bigcup_{k\in I_l}I^{a_k}_n\ ,\cr
&I_N^{a_1a_2}:=I_{2}^{\diag}\cup I_2^{a_1a_2}\ ,\cr
&I_{M_n,l}:=I^{1\ldots l}_n\ .}$$}

\vskip 6pt

\subsec{Determinantal combinatorics}

In this section, we will consider a general surjection $m:I_g\times I_g\to I_M$, satisfying
Eq.\suriezione. For example, by using the construction of
section \notation, it is possible to define $m$ so that
$m(\1_i,\2_i)=m(\2_i,\1_i)=i$, $i\in I_M$; the corresponding
isomorphism is $A$, defined in Definition
\definizione.

For each morphism $s:I_M\to I_M$ consider the $g$-tuples $d^k(s)$,
$k\in I_{g+1}$, where \eqn\ild{d^i_j(s)=d^{j+1}_i(s)=s_{m(i,j)}\ ,}
$i\le j\in I_g$. Note that if $s$ is a
monomorphism, then each $g$-tuple consists of distinct integers, and
each $i\in I_M$ belongs to two distinct $g$-tuples.

Consider $\perm_g^{g+1}\equiv\underbrace{\perm_g\times\cdots\times\perm_g}_{g+1\hbox{ times}}$ and
define $\varkappa:\perm_g^{g+1}\times I_M\rightarrow I_M$,
depending on $m$, by
\eqn\ilkappa{\varkappa_{m(i,j)}(r^1,\ldots,r^{g+1})=m(r^i_j,r^{j+1}_i)\
,} $i\le j\in I_g$, where
$(r^1,\ldots,r^{g+1})\in\perm_g^{g+1}$. Note that
$$d^i_j(\varkappa(r^1,\ldots,r^{g+1}))=
d^{j+1}_i(\varkappa(r^1,\ldots,r^{g+1}))=m(r^i_j,r^{j+1}_i)\ ,$$
$i\le j\in I_g$. Consider the subset of $I_M$ determined
by
$$I_{M,n}:=\{m(i,j)|i\in I_n,\,j\in I_g\}\ ,$$ $n\in I_g$,
with the ordering inherited from $I_M$, and denote by
\eqn\Lcard{L:=M-(g-n)(g-n+1)/2 \ ,} its cardinality. The elements
$\varkappa_l(r^1,\ldots,r^{g+1})$, $l\in I_{M,n}$, are independent
of $r^j_i$, with $n+1\le i,j\le g$, and $\varkappa$ can be
generalized to a function
$\varkappa:I_{M,n}\times\tilde\perm^{g,n}\to I_M\ ,$ where
$\tilde\perm^{g,n}:=\perm_g^n\times\perm_n^{g-n+1}$, by
\eqn\ilkappaii{\varkappa_i(\tilde r^1,\ldots,\tilde r^{g+1}):=
\varkappa_i(r^1,\ldots,r^{g+1})\ ,} $i\in I_{M,n}$, $(\tilde
r^1,\ldots,\tilde r^{g+1})\in\tilde\perm^{g,n}$, where
$r^j\in\perm_g$, $j\in I_{g+1}$, are permutations satisfying
$r^j=\tilde r^j$, $j\in I_n$, and $r^j_i=\tilde r^j_i$, $i\in I_n$,
$n+1\le j\le g$. Furthermore, if $\{\varkappa_i(\tilde
r^1,\ldots,\tilde r^{g+1})\}_{i\in I_{M,n}}$ consists of distinct
elements, then it is a permutation of $I_{M,n}$. By a suitable
choice of the surjection \eqn\lam{m(j,i)=m(i,j):=M-(g-j)(g-j-1)/2+i\
,} $j\le i\in I_g$, we obtain $I_{M,n}=I_L$ as an equality among
ordered sets. Note that this choice for $m$, which is convenient to
keep the notation uncluttered, does not correspond to the one
introduced in subsection \notation.

Consider the maps $s:I\rightarrow I$, where $I$ is any ordered
subset of $I_M$; if $s$ is bijective, then it is a permutation of
$I$. We define the function $\sgn(s)$ to be the sign of the
permutation if $s$ is bijective, and zero otherwise.

\vskip 6pt

\subsec{Combinatorial lemmas}\subseclab\detlemmas

Let $F$ be a commutative field and $S$ a non-empty set. Fix a set
$f_i$, $i\in I_g$, of $F$-valued functions on $S$, and $x_i\in S$,
$i\in I_M$. Set
$$ff_{m(i,j)}:=f_if_j\ ,$$
$i,j\in I_g$, and
$$\det f(x_{d^j(s)}):=\det\nolimits_{ik} f_k(x_{d^j_i(s)})\ ,$$
$j\in I_{g+1}$, where $x_i\in S$, $i\in I_M$. Furthermore, for any
ordered set $I\subseteq I_M$, we denote by
$$\det\nolimits_I ff(x_1,\ldots,x_{{\rm Card}(I)})\ ,$$
the determinant of the matrix $(ff_m(x_i))_{^{i\in I_{{\rm
Card}(I)}}_{m\in I}}$.

\newth\thcombi{Lemma}{Choose $n\in I_g$ and $L$
points $x_i$ in $S$, $i\in I_L$, with $L$ given by \Lcard. Fix $g-n$
points $p_i\in S$, $n+1\le i\le g$ and $g$ $F$-valued functions
$f_i$ on $S$, $i\in I_g$. The following $g(g-n)$ conditions
\eqn\combcond{ f_i(p_j)=\delta_{ij}\ ,} $1\le i\le j$, $n+1\le j\le
g$, imply \eqn\combi{\det\nolimits_{I_{M,n}}
ff(x_1,\ldots,x_L)={1\over
c_{g,n}}\sum_{s\in\perm_L}\sgn(s)\prod_{j=1}^n\det
f(x_{d^j(s)})\prod_{k=n+1}^{g+1} \det
f(x_{d^k_1(s)},\ldots,x_{d^k_n(s)},p_{n+1},\ldots,p_g)} where
\eqn\lacost{c_{g,n}:=\sum_{(\tilde r^1,\ldots,\tilde
r^{g+1})\in\tilde\perm^{g,n}} \prod_{k=1}^{g+1}\sgn(\tilde r^k)
\sgn(\varkappa(\tilde r^1,\ldots,\tilde r^{g+1}))\ .} In particular,
for $n=g$ \eqn\combii{c_g\det
ff(x_1,\ldots,x_M)=\sum_{s\in\perm_M}\sgn(s) \prod_{j=1}^{g+1}\det
f(x_{d^j(s)})\ ,} where
$$
c_g:=c_{g,g}=\sum_{r^1,\ldots,
r^{g+1}\in\perm_g}\prod_{k=1}^g\sgn(r^k)\,\sgn(\varkappa(r^1,\ldots,r^g))
\ .
$$}

\vskip 6pt

\noindent {\bf Proof.} It is convenient to fix the surjection $m$ as
in \lam, so that $I_{M,n}=I_L$. Next consider
\eqn\combiii{c_{g,n}\det\nolimits_{I_L}
ff(x_1,\ldots,x_L)=c_{g,n}\sum_{s\in\perm_L}
\sgn(s)ff_1(x_{s_1})\cdots ff_L(x_{s_L})\ .} Restrict the sums in
\lacost\ to the permutations $(\tilde r^1,\ldots,\tilde r^{g+1})\in\perm^{g,n}$, $i\in I_n$, such that
$\sgn(\varkappa(\tilde r^1,\ldots,\tilde r^{g+1}))\neq 0$, and set
$s':=s\circ \varkappa(\tilde r^1,\ldots,\tilde r^{g+1})$, so that
$$ff_1(x_{s_1})\cdots
ff_L(x_{s_L})=ff_{\varkappa_1}(x_{s'_1})\cdots
ff_{\varkappa_L}(x_{s'_L})\ ,$$ where $\varkappa_i$ is to be
understood as $\varkappa_i(\tilde r^1,\ldots,\tilde r^{g+1})$. Note
that, for all $l\in I_M$, there is a unique pair $i,j\in I_g$, $i\le
j$, such that $l=m(i,j)$, and by \ild\ and \ilkappa\ the following
identity
$$ff_{\varkappa_l(r^1,\ldots,r^{g+1})}(x_{s'_l})=ff_{m(r^i_j,r^{j+1}_i)}(x_{s'_{m(i,j)}})=
f_{r^i_j}(x_{d^i_j(s')}) f_{r^{j+1}_i}(x_{d^{j+1}_i(s')})\ ,$$ holds
for all $(r^1,\ldots,r^{g+1})\in\perm_g^{g+1}$. On the other hand,
if $l\in I_L$, then $i\le n$ and by Eq.\ilkappaii\
$$ff_1(x_{s_1})\cdots ff_L(x_{s_L})
=\prod_{i=1}^n f_{\tilde r^i_1}(x_{d^i_1(s')})\cdots f_{\tilde
r^i_g} (x_{d^i_g(s')})\prod_{j=n+1}^{g+1} f_{\tilde
r^j_1}(x_{d^j_1(s')})\cdots f_{\tilde r^j_n} (x_{d^j_n(s')})\ .$$
The condition $f_i(p_j)=\delta_{ij}$, $i\le j$, implies
$$\sum_{\tilde r^j\in\perm_n}\sgn(\tilde r^j)f_{\tilde r^j_1}(x_{d^j_1(s')})
\cdots f_{\tilde r^j_n}(x_{d^j_n(s')})= \det
f(x_{d^j_1(s')},\ldots,x_{d^j_n(s')},p_{n+1},\ldots,p_g)\ ,$$
$n+1\le j\le g+1$. Hence, Eq.\combi\ follows by replacing the sum
over $s$ with the sum over $s'$ in \combiii, and using
$$\sgn(s)=\sgn(s')\,\sgn
(\varkappa(\tilde r^1,\ldots,\tilde r^{g+1}))\ .$$ Eq.\combii\ is an
immediate consequence of \combi.\hfill$\square$

\vskip 6pt

\newrem\rrremark{Remark}{The summation over $\perm_M$ in Eq.\combii\
yields a sum over $(g+1)!$ identical terms, corresponding to
permutations of the $g+1$ determinants in the product. Such an
overcounting can be avoided by summing over the following subset of
$\perm_M$
$$\perm'_M:=\{s\in \perm_M, \hbox{ {\it s.t.} } s_1=1,\ s_2<s_3<\ldots<s_g,\
s_2<s_i,\ g+1\leq i\leq 2g-1\}\ ,$$ and by replacing $c_g$ by
$c_g/(g+1)!$.}

\vskip 6pt

Direct computation gives
\eqn\glic{c_{g,1}=g!\ , \quad c_{g,2}=g!(g-1)!(2g-1)\ , \quad
c_2=6\ , \quad c_3=360\ , \quad c_4=302400\ .} For $g=2$, $c_g/(g+1)!=1$ and
$\perm'_{M=3}=\{(1,2,3)\}$, so that \eqn\peppp{\det
ff(x_1,x_2,x_3)=\det f(x_1,x_2)\det f(x_1,x_3)\det f(x_2,x_3)\ .}

A crucial point in proving Lemma \thcombi\ is that if
$\varkappa_i(\tilde r^1,\ldots,\tilde r^{g+1})$, $i\in I_{M,n}$, are
pairwise distinct elements in $I_M$, then they belong to
$I_{M,n}\subseteq I_M$, with $\varkappa$ a permutation of such an
ordered set. For a generic ordered set $I\subseteq I_M$, one should
consider $\varkappa$ as a function over $g+1$ permutations $\tilde
r^i$, $i\in I_{g+1}$, of suitable ordered subsets of $I_g$. In
particular, $\tilde r^i$ should be a permutation over all the
elements $j\in I_g$ such that $m(i,j)\in I$, for $j\ge i$, or
$m(i-1,j)\in I$, for $j<i$. However, the condition that the elements
$\varkappa_i(\tilde r^1,\ldots,\tilde r^{g+1})$, $i\in I$, are
pairwise distinct does not imply, in general, that they belong to
$I$ and Lemma \thcombi\ cannot be generalized to a determinant of
products $ff_i$, $i\in I$. On the other hand, the subsets
\eqn\laI{I:=I_{M,n}\cup \{m(i,j)\}\ ,} satisfy such a condition for
$n< i,j\le g$ and yield the following generalization of Lemma
\thcombi.

\vskip 6pt

\newth\thcombvi{Lemma}{Assume the hypotheses of Lemma {\rm \thcombi} for $n<g$,
and choose a point $x_{L+1}\in S$, and a pair $i,j$, $n< i,j\le g$.
Then the following relation \eqn\combvi{\eqalign{&\det\nolimits_I
ff(x_1,\ldots,x_{L+1})\cr &={1\over
c'_{g,n}}\sum_{s\in\perm_{L+1}}\sgn(s) \prod_{k=1}^n\det
f(x_{d^k(s)})\det
f(x_{d^{n+1}_1(s)},\ldots,x_{d^{n+1}_{n+1}(s)},p_{n+1},\ldots,\check
p_i,\ldots,p_g)\cr &\cdot\det
f(x_{d^{n+2}_1(s)},\ldots,x_{d^{n+2}_{n+1}(s)},p_{n+1},\ldots,\check
p_j,\ldots,p_g) \prod_{l=n+3}^{g+1}\det
f(x_{d^l_1(s)},\ldots,x_{d^l_n(s)},p_{n+1},\ldots,p_g) \ ,}} where
$$c'_{g,n}:=\sum_{(\tilde r^1,\ldots,\tilde r^{g+1})\in\tilde\perm^I}\prod_{i=1}^{g+1}\sgn(\tilde
r^i)\sgn(\varkappa(\tilde r^1,\ldots,\tilde r^{g+1}))\ ,$$
$\tilde\perm^I:=\perm_g^n\times\perm_{n+1}^2\times\perm_n^{g-n-1}$,
and $I$ is defined in {\rm \laI}, holds.}

\vskip 6pt

\noindent {\bf Proof.} A straightforward generalization of the proof
of Lemma \thcombi. \hfill$\square$

\vskip 6pt

\subsec{Examples of the combinatorial lemmas}

We now show some examples of the combinatorial construction described in the last subsection.
Set $g=4$, so that $M=g(g+1)/2=10$. Fix a surjection
$m:I_4\times I_4\to I_{10}$ with $m(i,j)=m(j,i)$, for example by setting $m(i,j)=[m]_{ij}$, with
$[m]$ the symmetric matrix
$$[m]=
\left(\matrix{
1 & 2 & 3 & 4\cr
2 & 5 & 6 & 7\cr
3 & 6 & 8 & 9\cr
4 & 7 & 9 & 10
}\right)
\ .$$
For each function $s:I_{10}\to I_{10}$, the $4$-tuples $d^i(s)$, $i=1,\ldots,g+1=5$, are determined by
$$d_j^i(s)=d_i^{j+1}(s)=s_{m(i,j)}\ ,$$ $i\le j\in I_g$,
so that, with the above choice of $m$,
$$\eqalignno{
d^1(s)&=(s_1,s_2,s_3,s_4)\ ,\cr d^2(s)&=(s_1,s_5,s_6,s_7)\ ,\cr
d^3(s)&=(s_2,s_5,s_8,s_9)\ ,\cr d^4(s)&=(s_3,s_6,s_8,s_{10})\ ,\cr
d^5(s)&=(s_4,s_7,s_9,s_{10})\ . }$$ Let $\perm_g$ be the group of
permutations of $g$ elements. The function
$\varkappa:\perm_4^5\times I_{10}\to I_{10}$ is defined by
\eqn\ilvarkappa{\varkappa_{m(i,j)}(r^1,\ldots,r^5)=m(r_j^i,r_i^{j+1})\
,} $ i\le j\in I_g$, where $(r^1,\ldots,r^5)\in\perm_4^5$. For
example, fix
$$\eqalignno{
r^1&=(3,4,1,2)\ ,\cr
r^2&=(1,2,4,3)\ ,\cr
r^3&=(2,4,1,3)\ ,\cr
r^4&=(1,2,3,4)\ ,\cr
r^5&=(2,4,1,3)\ .
}$$
To determine $\varkappa_1(r^1,\ldots,r^5)$, note that $1=m(1,1)$, so that, by definition,
$$\varkappa_{m(1,1)}(r^1,\ldots,r^5)=m(r_1^1,r^2_1)=m(3,1)=3\ .$$
As a further example note that $2=m(1,2)=m(2,1)$, so that
$$\varkappa_{m(1,2)}(r^1,\ldots,r^5)=m(r^1_2,r^3_1)=m(4,2)=7\ ,$$
(observe that Eq.\ilvarkappa, which defines $\varkappa$, holds only
for $i\le j$). The $4$-tuples $d^i(\varkappa(r_1,\ldots,r_5))$ are
$$\eqalignno{
d^1(\varkappa)&=(3,7,1,5)\ ,\cr d^2(\varkappa)&=(3,7,7,9)\ ,\cr
d^3(\varkappa)&=(7,7,3,3)\ ,\cr d^4(\varkappa)&=(1,7,3,9)\ ,\cr
d^5(\varkappa)&=(5,9,3,9)\ . }$$ It is readily verified the general
relation
$$d^i_j(\varkappa(r^1,\ldots,r^5))=d^{j+1}_i(\varkappa(r^1,\ldots,r^5))=m(r^i_j,r^{j+1}_i)\ ,$$
$ i\le j\in I_g$. Note that if
$\varkappa(r^1,\ldots,r^5):I_{10}\to I_{10}$, for some fixed
$r^1,\ldots,r^5$, is a monomorphism, then it determines a
permutation of $I_{10}$. Hence, we can define the function
$\sgn(\varkappa(r^1,\ldots,r^5))$ to be the sign of the permutation
$\varkappa(r^1,\ldots,r^5)$ if it is a monomorphism, and zero
otherwise.

\noindent Consider the subset
$$I_{M,n}=\{m(i,j)\mid i\in I_n, j\in I_g\}\ ,$$
for some $n\in I_g$. $\varkappa$ can be
generalized to a function from $\tilde\perm^{g,n}\times I_{M,n}$,
where $\tilde\perm^{g,n}:=\perm_g^n\times\perm_n^{g-n+1}$, into
$I_M$. As an example, consider
$\varkappa:\tilde\perm^{4,2}\times I_{10,2}\to I_{10}$, where
$I_{10,2}=\{1,2,3,4,5,6,7\}$ (the precise form of $I_{10,2}$ depends
on the choice of $m$). Fix $(\tilde r^1,\ldots,\tilde
r^5)\in\tilde\perm^{4,2}=\perm_4^2\times\perm_2^{3}$, $e.g.$ by
$$\eqalignno{
\tilde r^1&=(3,4,1,2)\ ,\cr \tilde r^2&=(1,2,4,3)\ ,\cr \tilde
r^3&=(2,1)\ ,\cr \tilde r^4&=(1,2)\ ,\cr \tilde r^5&=(1,2)\ . }$$ As
a specific case, say $\varkappa_6$, note that $6=m(2,3)=m(3,2)$ and
set $$\varkappa_{m(2,3)}(\tilde r^1,\ldots,\tilde r^5)=m(\tilde
r_3^2,\tilde r^4_2)=m(4,2)=7\ .$$ For general choices of $\tilde
r^1,\ldots,\tilde r^5$, $\varkappa(\tilde r^1,\ldots,\tilde
r^5):I_{10,2}\to I_{10}$ may not be a monomorphism. It can be
verified that if the image $\varkappa(\tilde r^1,\ldots,\tilde
r^5)(I_{10,2})\not\subseteq I_{10,2}$, then $\varkappa(\tilde
r^1,\ldots,\tilde r^5)$ is not a monomorphism. Therefore, if
$\varkappa(\tilde r^1,\ldots,\tilde r^5)$ is a monomorphism, then it
determines a permutation of $I_{10,2}$. Hence, we can define the
function $\sgn(\varkappa(\tilde r^1,\ldots,\tilde r^5))$ to be the
sign of $\varkappa(\tilde r^1,\ldots,\tilde r^5)$ if it is a
monomorphism, and zero otherwise.

Let us apply Lemma \thcombi\ to the previous examples. Consider four
linearly independent functions $f_1,\ldots,f_4:\CC\to\CC$, and set
$$ff_{m(i,j)}(z):=f_i(z)f_j(z)\ .$$
Next, fix $x_1,\ldots,x_{10}\in\CC$ and consider
$$\det\left(\matrix{
ff_1(x_1) & \ldots & ff_{10}(x_1)\cr \vdots & \ddots & \vdots\cr
ff_1(x_{10}) & \ldots & ff_{10}(x_{10}) }\right)=\det\left(\matrix{
f_1(x_1)f_1(x_1) & \ldots & f_4(x_1)f_4(x_1)\cr \vdots & \ddots &
\vdots\cr f_1(x_{10})f_1(x_{10}) &\ldots & f_4(x_{10})f_4(x_{10})
}\right)\ ,$$ so that $m(i,j)$ determines the column where $f_if_j$
appears. It is easily verified that the above determinant is
proportional to
\eqn\formula{\eqalign{
\sum_{s\in\perm_{10}}\sgn(s)\det f_i(x_{d^1_j(s)})
\det f_i(x_{d^2_j(s)}) \det
f_i(x_{d^3_j(s)})\det f_i(x_{d^4_j(s)})
\det f_i(x_{d^5_j(s)})\ .}}

This expression, after expanding each determinant, consists of a summation over products of twenty factors
$f_i(x_j)$, where each $x_k$ appears twice. After skew-symmetrization of the $x_k$'s, this expression is
necessarily proportional to the original determinant.

In Lemma \thcombi\ it is also considered the more general case of determinants made up of functions
$ff_i$, where $i$ varies in a subset $I_{M,n}\subset I_M$ of $L<M$ elements. For example, let us
consider the subset $I_{10,2}=\{1,\ldots,7\}$ and fix the points $x_1,\ldots,x_7\in\CC$. We are interested
in the determinant
\eqn\euna{\det\left(\matrix{
ff_1(x_1) & \ldots & ff_{7}(x_1)\cr
\vdots & \ddots & \vdots\cr
ff_1(x_{7}) & \ldots & ff_{7}(x_{7})
}\right)=\det\left(\matrix{
f_1(x_1)f_1(x_1) & \ldots & f_2(x_1)f_4(x_1)\cr
\vdots & \ddots & \vdots\cr
f_1(x_{7})f_1(x_{7}) &\ldots & f_2(x_{7})f_4(x_{7})
}\right)\ .}
By repeating the above construction, this determinant can be expressed as (a sum over) products of two
determinants of $4\times 4$ matrices
times three determinants of lower-dimensional $2\times 2$ matrices
$$\sum_{s\in\perm_{10}}\sgn(s)\det\nolimits_{I_4} f_i(x_{d^1_j(s)})
\det\nolimits_{I_4} f_i(x_{d^2_j(s)})\det\nolimits_{I_2} f_i(x_{d^3_j(s)})
\det\nolimits_{I_2} f_i(x_{d^4_j(s)})
\det\nolimits_{I_2} f_i(x_{d^5_j(s)})\ ,$$
where $\det_{I_n} f_i(x_j):=\det_{ij\in I_n}f_i(x_j)$.
In order to obtain products of five determinants of $4\times 4$ matrices in the
form similar to Eq.\formula, one has to impose some conditions on the functions $f_i$. In
particular, it is sufficient to require that there exist two points, $p_3,p_4\in\CC$, such that
$$\eqalignno{&f_1(p_i)=f_2(p_i)=0\ ,\qquad i=3,4\ ,\cr
&f_3(p_4)=f_4(p_3)=0\ ,\cr
&f_3(p_3)=f_4(p_4)=1\ .}$$
In this case, the following identity
$$\det \left(\matrix{
f_1(x_1) & f_{2}(x_1)\cr
f_1(x_2) & f_{2}(x_2)\cr
}\right)=\det\left(\matrix{
f_1(x_1) & f_{2}(x_1) & f_3(x_1) & f_{4}(x_1)\cr
f_1(x_2) & f_{2}(x_2) & f_3(x_2) & f_{4}(x_2)\cr
f_1(p_3) & f_{2}(p_3) & f_3(p_3) & f_{4}(p_3)\cr
f_1(p_4) & f_{2}(p_4) & f_3(p_4) & f_{4}(p_4)\cr
}\right)\ ,$$
holds and the determinants in \euna\ are proportional to
\eqn\edue{\eqalign{
\sum_{s\in\perm_{7}}\sgn(s)&\det f_i(x_{d^1_j(s)})
\det f_i(x_{d^2_j(s)})
\det f(x_{d^3_1(s)},x_{d^3_2(s)},p_3,p_4)\cr
\cdot&\det f(x_{d^4_1(s)},x_{d^4_2(s)},p_3,p_4)
\det f(x_{d^5_1(s)},x_{d^5_2(s)},p_3,p_4)\ ,
}}
where $\det f(z_1,\ldots,z_4):=\det_{ij\in I_4}f_i(z_j)$. Lemma \thcombi\ generalizes
such a result to any $g$ and $n$. Proportionality of Eqs.\euna\ and \edue\ can be understood as
follows. Upon expanding the determinants in \edue\ and using the
conditions on $f_i$, this expression corresponds to a summation of products of the form
\eqn\prodotto{
f_1f_2f_3f_4\cdot f_1f_2f_3f_4 \cdot f_1f_2\cdot f_1f_2 \cdot
f_1f_2\ ,}
with the $f_i$'s evaluated at $x_1,\ldots,x_7$ (each $x_i$ appears twice). Such a
product can be re-arranged as
$$ff_{i_1}(x_1)ff_{i_2}(x_2)\ldots ff_{i_7}(x_7)\ ,$$
for some $i_1,\ldots,i_7\in I_{10}$. After skew-symmetrization over the variables $x_i$, only the
products with distinct $i_1,\ldots,i_7$ contribute. But this implies $i_1,\ldots,i_7\in I_{10,2}$,
since the only possibility to construct seven different functions $f_if_j$ out of the fourteen
functions in Eq.\prodotto\ is
\eqn\term{f_1^2(x_1)f_1f_2(x_2)f_1f_3(x_3)f_1f_4(x_4)f_2^2(x_5)f_2f_3(x_6)f_2f_4(x_7)\
,}
up to permutations of the $x_i$'s.
This is strictly related to the observation that if $\varkappa(\tilde r^1,\ldots,\tilde r^5)$ is a
monomorphism, then it corresponds to a permutation of $I_{10,2}$. The skew-symmetrization of
\term\ with respect to $x_1,\ldots,x_7$ is exactly the determinant we were
looking for.

Note that Lemma \thcombi\ may not be generalized to the case of
determinants of matrices with rows $ff_{i_1},\ldots,ff_{i_L}$, when
$I:=\{i_1,\ldots,i_L\}$ is a generic subset of $I_{10}$. One can
always define a generalization of the $\varkappa$ function as
$\varkappa(\tilde r^1,\ldots,\tilde r^5):I\to I_{10}$, with $\tilde
r^1,\ldots,\tilde r^5$ in some suitable subset of $\perm_4^5$.
However, the necessary condition for the generalization of Lemma
\thcombi\ is that if $\varkappa$ is a monomorphism, then
$\varkappa(I)=I$. Such a condition is verified, for example, if
$I=I_{10,n}$, as showed before for $I_{10,2}$. The condition still
holds when $I=I_{10,n}\cup\{j\}$, for all the elements $j\in
I_{10}\setminus I_{10,n}$, which is the content of Lemma \thcombvi.
An example for which the analog of Lemma \thcombi\ does not exist is
for $I=\{1,5,8,10\}$, corresponding to determinants of matrices with
rows $f_1^2,f_2^2,f_3^2,f_4^2$. Actually, defining a formula similar
to \formula\ in order to obtain terms in the form
$f_1^2(x_1)f_2^2(x_2)f_3^2(x_3)f_4^2(x_4)$, some unwanted terms,
such as $f_1f_2(x_1)f_2f_3(x_2)f_3f_4(x_3)f_4f_1(x_4)$, do not
cancel on the RHS.

\vskip 6pt

\newsec{Divisors of higher order theta derivatives on Riemann
surfaces}\seclab\Riemanndef

After reminding some basic facts about theta functions, we
investigate the divisor structures of the theta function and its
derivatives that will be used in the subsequent sections.

Set $\abvar_Z:=\CC^g/L_Z$, $L_Z:={\ZZ}^g +Z{\ZZ}^g$, where $Z$
belongs to the Siegel upper half-space
$$\H_g:=\{Z\in M_g(\CC) |\,{}^tZ=Z,
\,\im Z>0\}\ ,$$ and consider the theta function with
characteristics \eqn\thetaconc{\eqalign{\theta
\left[^a_b\right]\left(z,Z\right):&= \sum_{k\in {\ZZ}^g}e^{\pi i
\tp{(k+a)}Z(k+a)+ 2\pi i \tp{(k+a)} (z+b)}\cr &=e^{\pi i
\tp{a}Za+2\pi i\tp{a(z+b)}}\theta
\left[^0_0\right]\left(z+b+Za,Z\right) \ ,}} where $z\in \abvar_Z$,
$a,b\in{\RR}^g$. It has the quasi-periodicity properties
$$\theta \left[^a_b\right]\left(z+n+Zm,Z\right)=
e^{-\pi i \tp{m}Zm-2\pi i\tp{m}z+2\pi i(\tp{a}n-\tp{b}m)}\theta \left[^a_b\right]\left(z,Z\right) \ ,
$$
$m,n\in\ZZ^g$. Denote by $\Theta\subset \abvar_Z$ the divisor of
$\theta(z,Z):=\theta\left[^0_0\right](z,Z)$ and by
$\Theta_{s}\subset \Theta$ the locus where $\theta$ and its gradient
vanish. Geometrically $\theta\left[^a_b\right](z,Z)$ is the unique
holomorphic section of the bundle ${\cal L}_{\Theta_{ab}}$ on
$\abvar_Z$ defined by the divisor $\Theta_{ab}=\Theta+b+Za$ of
$\theta\left[^a_b\right](z,Z)$. A suitable norm, continuous
throughout $\abvar_Z$, is given by
$$
||\theta||^2\left(z,Z\right)= e^{-2\pi\tp\im z (\im Z)^{-1} \im\bar z}
|\theta|^2\left(z,Z\right) \ .
$$
Computing $c_1({\cal L}_\Theta)$ and using the
Hirzebruch-Riemann-Roch Theorem, it can be proved that $\theta$ is
the unique holomorphic section of ${\cal L}_\Theta$. It follows that
$(\abvar_Z,{\cal L}_\Theta)$ is a principally polarized abelian
variety.

\vskip 6pt

\subsec{Riemann theta functions and the prime form}\subseclab\thetaprime

Let
$\{\alpha,\beta\}\equiv\{\alpha_1,\ldots,\alpha_g,\beta_1,\ldots,\beta_g\}$
be a symplectic basis of $H_1(C,\ZZ)$ and $\{\omega_i\}_{i\in I_g}$
the basis of $H^0(K_C)$ satisfying the standard normalization
condition $\oint_{\alpha_i}\omega_j=\delta_{ij}$, forall $i,j\in
I_g$. Let $\tau\in\H_g$ be the Riemann period matrix of $C$,
$\tau_{ij}:=\oint_{\beta_i}\omega_j$. A different choice of the
symplectic basis of $H_1(C,\ZZ)$ corresponds to a
$\Gamma_g:=Sp(2g,\ZZ)$ transformation
$$\left(\matrix{\alpha\cr
\beta}\right)\mapsto \left(\matrix{\tilde\alpha\cr
\tilde\beta}\right)=\left(\matrix{D & C \cr B & A}\right)
\left(\matrix{\alpha\cr \beta}\right)\ ,\qquad\qquad \left(\matrix{A
& B \cr C & D}\right)\in \Gamma_g \ , $$ \eqn\modull{\tau \mapsto
\tau'=(A\tau+B)(C\tau+D)^{-1}\ .} We denote by
$\ppavmod_g:=\H_g/\Gamma_g$ the moduli space of principally
polarized abelian varieties.

Choose an
arbitrary point $p_0\in C$ and let $I(p):=(I_1(p),\ldots,I_g(p))$
$$I_i(p):=\int_{p_0}^p\omega_i\ ,$$ $p\in C$, be the Abel-Jacobi
map. It embeds $C$ into the Jacobian $J_0(C):={\CC}^g/L_\tau$,
$L_\tau:={\ZZ}^g +\tau {\ZZ}^g$, and generalizes to a map from the
space of divisors of $C$ into $J_0(C)$ as $I(\sum_i n_i
p_i):=\sum_in_iI(p_i)$, $p_i\in C$, $n_i\in\ZZ$. By Jacobi Inversion Theorem
the restriction of $I$ to $C_g$ is a surjective map onto $J_0(C)$.
The Hodge metric of the polarization of $J_0(C)$ has the K\"ahler form
\eqn\khaeler{\nu ={i\over2}\sum_{i,j=1}^g(\imtau^{-1})_{ij} dz_i\wedge d\bar z_j \ ,\qquad
} whose pullback $I^*\nu=g\mu$ defines the Bergman two-form on $C$
\eqn\metricasuc{\mu={i\over2g}\sum_{i,j=1}^g(\imtau^{-1})_{ij}\omega_i\wedge\bar\omega_j \ .\qquad
}

If $\delta',\delta''\in\{0,1/2\}^g$, then $\theta
\left[\delta\right]\left(z,\tau\right):=\theta
\left[^{\delta'}_{\delta''}\right]\left(z,\tau\right)$ has definite parity in $z$
$$\theta \left[\delta\right] \left(-z,\tau\right)=e(\delta)
\theta
\left[\delta\right] \left(z,\tau\right)\ ,$$
where $e(\delta):=e^{4\pi i\!\tp{\delta'} \delta''}$.
There are $2^{2g}$ different characteristics for which $\theta
\left[\delta\right] \left(z,\tau\right)$ has definite parity.
By Abel Theorem each one of such theta characteristics
determines the divisor class of a spin bundle $L_\alpha\simeq
K^{1\over2}_C$, so that we may call them spin structures. There are $2^{g-1}(2^g+1)$ even and $2^{g-1}(2^g-1)$
odd spin structures.

Consider the vector of Riemann constants
\eqn\scriviamolo{\K^p_i:={1\over 2}+{1\over 2}\tau_{ii}-\sum_{j\neq
i}^g\oint_{\alpha_j}\omega_j\int_{p}^x\omega_i\ ,} $i\in I_g$,
for all $p\in C$. For any $p$ we define the Riemann divisor
class by
\eqn\lakkapppaa{
I(\Delta):=(g-1)I(p)-\K^p \ ,
}
which has the property $2\Delta = K_C$.

In the following, we will consider $\theta(D+e):=\theta
\left[^0_0\right](I(D)+e,\tau)$, for all $e\in J_0(C)$, evaluated at
some $0$-degree divisor $D$ of $C$. We will also use the notation
$$\deltadiv(D):=\theta(I(D-n\Delta))\ ,$$
for each divisor $D$ of degree $n(g-1)$, $n\in\ZZ$.

According to the Riemann Vanishing Theorem, for any $p\in C$ and $e\in J_0(C)$

\smallskip

\item{\it i.} if $\theta(e)\ne 0$, then the divisor $D$ of $\theta(x-p-e)$
in $C$ is effective of degree $g$, with index of specialty $i(D)=0$
and $e=I(D-p-\Delta)$;

\smallskip

\item{\it ii.} if $\theta(e)=0$, then for some $\zeta\in C_{g-1}$,
$e=I(\zeta-\Delta)$.

\vskip 6pt

By Riemann's Singularity Theorem it follows that the dimension of $\Theta_{s}$ for $g\geq4$
is $g-3$ in the hyperelliptic case and $g-4$ if the curve is canonical.

Let $\nu$ a non-singular odd characteristic. The holomorphic
1-differential \eqn\hquadro{h^2_\nu(p):=\sum_{1}^g\omega_i(p)
\partial_{z_i}\theta\left[\nu\right](z)_{|_{z=0}}\ ,} $p\in C$,
has $g-1$ double zeros. The prime form \eqn\pojdlk{
E(z,w):={\theta\left[\nu\right](w-z,\tau)\over
h_{\nu}(z)h_{\nu}(w)}\ ,} is a holomorphic section of a line bundle
on $C\times C$, corresponding to a differential form of weight
$(-1/2,-1/2)$ on $\tilde C\times \tilde C$, where $\tilde C$ is the
universal cover of $C$. It has a first order zero along the diagonal
of $C\times C$. In particular, if $t$ is a local coordinate at $z\in
C$ such that $h_\nu=dt$, then
$$
E(z,w)={t(w)-t(z)\over\sqrt{dt(w)}\sqrt{dt(z)}}(1+{\cal O}((t(w)-t(z))^2)) \ .
$$
Note that $I(z+\tp\alpha n+\tp\beta m)=I(z)
+n+\tau m$, $m,n\in\ZZ^g$, and
$$
E(z+\tp\alpha n+\tp\beta m,w)=\chi e^{-\pi i \tp m \tau m -2\pi i \tp mI(z-w)}E(z,w) \ ,
$$
where $\chi:=e^{2\pi i(\tp \nu' n- \tp \nu'' m)}\in\{-1,+1\}$, $m,n\in\ZZ^g$.

We will also consider the multi-valued $g/2$-differential $\sigma(z)$ on $C$ with empty divisor, that is a
holomorphic section of a trivial bundle on $C$, and
satisfies the property
$$
\sigma(z+\tp\alpha n+\tp\beta m)=\chi^{-g}e^{\pi i(g-1)\tp m \tau m+2\pi i \tp m
\K^z}\sigma(z)\ .
$$
Such conditions fix $\sigma(z)$ only up to a factor independent of
$z$; the precise definition, to which we will refer, can be given,
following \ref\FayMAM{J.~Fay, Kernel functions, analytic torsion and
moduli spaces, {\it Mem.\ Am.\ Math.\ Soc.\ } {\bf 96} (1992).}, on
the universal covering of $C$  (see also \ref\jfayy{J.~Fay, {\it
Theta Functions on Riemann surfaces}, Springer Lecture Notes 352
(1973).}). Furthermore, \eqn\equationn{
\sigma(z,w):={\sigma(z)\over\sigma(w)}={\deltadiv(\sum_{1}^gx_i-z)\over
\deltadiv(\sum_{1}^gx_i-w)}\prod_{i=1}^g{E(x_i,w)\over E(x_i,z)} \
,} for all $z,w,x_1,\ldots,x_g\in C$, which follows by observing
that the RHS is a nowhere vanishing section both in $z$ and $w$ with
the same multi-valuedness of $\sigma(z)/\sigma(w)$.

\vskip 6pt

\newrem\nuooov{Definition}{For all $z\in C$ and $\nu$ non-singular theta characteristics,
set
\eqn\nuovasigma{\sigma_{\nu}(z):=h_{\nu}(z)^g\exp\bigg(-\sum_{i=1}^g\oint_{\alpha_i}\omega_i(w)\log
\theta\left[\nu\right](w-z)\bigg)\ .} Such a
$g/2$-differential satisfies the same general properties of
$\sigma$, and is defined directly on $C$.}

\vskip 6pt

\noindent Under the modular transformations $z\to z'=z(CZ+D)^{-1}$,
$Z\to Z'=(AZ+B)(CZ+D)^{-1}$ the theta characteristics
transform as
$$\Bigl(\matrix{a\cr b}\Bigr)\to\Bigl(\matrix{\tilde a\cr \tilde b}\Bigr)=\left(\matrix{D & -C \cr
-B & A}\right) \Bigl(\matrix{a\cr b}\Bigr)\ ,
$$ $G:=\Bigl(\matrix{A & B\cr C & D}\Bigr)\in \Gamma_g$, for all $a,b,z\in \CC^g$, and the theta functions transform
as \ref\Igusa{J.~I.~Igusa, {\it Theta functions}, Springer-Verlag, 1972.}
$$\theta[^a_b](z,Z)\to \theta[^{a'}_{b'}](z',Z')=\epsilon_G(\det (CZ+D))^{1\over2}
e^{2\pi i\phi[^a_b](G)+\pi i
\tp z(CZ+D)^{-1}Cz}\theta[^a_b](z,Z)\ ,$$
where $\epsilon_G$ is an eighth root of 1 depending only on $G$,
$$\Bigl(\matrix{a'\cr b'}\Bigr):=\Bigl(\matrix{\tilde a\cr \tilde b}\Bigr)
+{1\over 2}\Bigl(\matrix{{\rm diag}\,(C\tp D)\cr {\rm
diag}\,(A\tp
B)}\Bigr)\ ,
$$
and
$$2\phi[^a_b](G):=(\tp a \ \tp b)\left(\matrix{ -\tp BD & \tp BC\cr \tp BC & -\tp
AC}\right)\Bigl(\matrix{a\cr b}\Bigr)+{\rm diag}(A\tp B)\cdot (Da-Cb)\ .$$

\vskip 6pt

Let $\omega(z,w)$ be the unique symmetric differential
on $C\times C$, with only a double pole along $z=w$, satisfying $\oint_{\alpha_j}\omega(z,w)=0$ and $\oint_{\beta_j}\omega(z,w)=2\pi i \omega_j$, $j\in I_g$.
The latter conditions imply that under a modular transformation
$$
\hat\omega(z,w)=\omega(z,w)-2\pi i \tp \omega(z)(C\tau+D)^{-1}\omega(w)\ .
$$
Since $E(z,w)$ is the unique antisymmetric solution of $\partial_z\partial_w \log E=\omega(z,w)$ which is consistent
with the expansion
of $\omega(z,w)$ for $z\sim w$, it follows that
$$
\hat E(z,w)=E(z,w)e^{\pi i (C\tau+D)^{-1} C\int_z^w\omega\cdot\int_z^w\omega} \ ,
$$
for all $z,w\in C$.

\vskip 6pt

\newth\riemannconstmod{Lemma}{{\rm (Fay \FayMAM)} If $\{\alpha,\beta\}$ and $\{\tilde\alpha,\tilde\beta\}$ are two markings
of $C$ related by \modull\ and $\K^q$ and $\K^{q\, '}$ denote the respective vectors of Riemann constants for $q\in C$,
then there are $a_0,b_0\in({1\over2}\ZZ)^g$, depending on the markings, such that
$$
a_0-{1\over2}{\rm diag}\,(C\tp D)\in \ZZ^g \ , \qquad b_0-{1\over2}{\rm diag}\,(A\tp B)\in \ZZ^g\ ,
$$
$$
\K^{q\, '}=\tp(C\tau+D)^{-1}\K^q+\tau'a_0+b_0\in \CC^g\ ,
$$
and
$$
\theta(z'+\K^{q\, '},\tau')=
\epsilon'(\det (C\tau+D))^{1\over2}
e^{\pi i
\tp(z+\K^q)(C\tau+D)^{-1}C(z+\K^q)-\pi i\tp a_0\tau' a_0-2\pi i \tp(C\tau+D)^{-1}(z+\K^q)}\theta(z,\tau)\ ,
$$
for all $z=\tp(C\tau+D) z'\in\CC^g$, with $\epsilon'$ an eighth root of 1 depending on the markings.}

\vskip 6pt

\vskip 6pt

Theta functions and, in particular, {\it Thetanullwerte}, i.e. theta constants
$\theta_\nu(0)$, with $\nu$ even characteristics, can be
used to construct modular forms, i.e. meromorphic functions on
$\H_g$ which are invariant under modular transformations. Some
regularity conditions at infinity are also required for $g=1$, which are not necessary
for $g>1$ due to the Koecher principle. More
generally, one considers modular forms of weight $k<0$, i.e.
holomorphic functions $f$ on $\H_g$ which transform as
\eqn\modulla{f(Z')=\det(CZ+D)^{-k}f(Z)\ ,}
under modular transformations or other discrete subgroups of
${\rm Sp}(2g,\RR)/\ZZ_2$, the group of automorphisms of $\H_g$.

\vskip 6pt

\subsec{Determinants in terms of theta functions}

Set
\eqn\essesimm{S(p_1+\ldots+p_g):=
{\theta_\Delta(\sum_1^g p_i-y) \over \sigma(y)\prod_1^gE(y,p_i)} \
,} $y,p_1,\ldots,p_g\in C$.

\newth\lemessesimm{Lemma}{For all $p_1,\ldots,p_g\in C$, $S(p_1+\ldots+p_g)$ is independent of $y$.
For each fixed $\d\in C_{g-1}$, consider the map $\pi_{d}:C\to C_g$,
$\pi_d(p):=p+\d$. The pull-back $\pi_{d}^*S$ vanishes identically if
and only if $\d$ is a special divisor; if $\d$ is not special, then
$\pi_{d}^*S$ is the unique (up to a constant) holomorphic
$1/2$-differential such that $[(\pi_{d}^*S)+\d]$ is the canonical
divisor class.}

\vskip 6pt

\noindent {\bf Proof.} If $p_1+\ldots+p_g$ is a special divisor, the Riemann Vanishing Theorem
implies $S=0$ identically in $y$; if $p_1+\ldots+p_g$ is not special, $S$ is a single-valued
meromorphic section in $y$ with no zero and no pole. It follows that, in any case, $S$ is a constant
in $y$. This also shows that $S(p_1+\ldots+p_g)=0$ if and only if $p_1+\ldots+p_g$ is a special
divisor. Hence, if $\d\in C_{g-1}$ is a special divisor, $S(p+\d)=0$ for all $p\in C$. On the
contrary, if $\d$ is not special, then $h^0(K_C\otimes\O(-\d))=1$, and $S(p+\d)=0$ if and only if $p$ is one of the zeros of the
(unique, up to a constant) holomorphic section of $H^0(K_C\otimes\O(-\d))$, and this concludes the
proof.\hfill$\square$

\vskip 6pt

\newth\thdettheta{Proposition}{Fix $n\in\NN_+$, set
$N_n:=(2n-1)(g-1)+\delta_{n 1}$ and let $\{\phi_i^n\}_{i\in
I_{N_n}}$ be arbitrary bases of $H^0(K_C^n)$. There are constants
$\kenne[\phi^n]$ depending only on the marking of $C$ and on
$\{\phi_i^n\}_{i\in I_{N_n}}$ such that
\eqn\dettheta{\kuno[\phi^1]={
\det\phi_i^1(p_j)\sigma(y)\prod_1^gE(y,p_i)\over
\deltadiv\bigl(\sum_{1}^gp_i-y\bigr)\prod_1^g\sigma(p_i)
\prod_{i<j}^gE(p_i,p_j) }={ \det\phi_i^1(p_j)\over
S\bigl(\sum_{1}^gp_i\bigr)\prod_1^g\sigma(p_i)
\prod_{i<j}^gE(p_i,p_j) }\ ,} and
\eqn\detthetaii{\kenne[\phi^n]={\det\phi_i^{n}(p_j)\over
\theta_\Delta\bigl(\sum_{1}^{N_n}
p_i\bigr)\prod_{1}^{N_n}\sigma(p_i)^{2n-1}\prod_{i<j}^{N_n}
E(p_i,p_j)}\ ,}  for $n\geq2$, for all $y,p_1,\ldots,p_{N_n}\in C$.}

\vskip 6pt

\noindent {\bf Proof.} $\kenne[\phi^n]$ is a meromorphic function
with empty divisor with respect to $y,p_1,\ldots,p_{N_n}$.
\hfill$\square$

\vskip 6pt

\newrem\lasigmanu{Remark}{Replacing the $g/2$-differential $\sigma$ in \dettheta\ and \detthetaii\ by $\sigma_\nu$,
defined in \nuovasigma, defines the new constants
$\kenne_\nu[\phi^n]$.}

\vskip 6pt

For each set $\{\phi_i^n\}_{i\in I_{N_n}}\subset H^0(K_C^n)$,
consider $W[\phi^n](p_1,\ldots,p_{N_n}):=\det\phi^n_i(p_j)$, and the
Wronskian $W[\phi^n](p):=\det\partial_p^{j-1}\phi^n_i(p)$. If
$W[\phi^n](p_1,\ldots,p_{N_n})$ does not vanish identically, then,
for each $\{\phi_i^{n'}\}_{i\in I_{N_n}}\subset H^0(K_C^n)$, we have
the constant ratio \eqn\ratiodef{
{\kenne[\phi^n]\over\kenne[\phi^{n'}]}={W[\phi^{n'}](p_1,\ldots,p_{N_n})\over
W[{\phi}^{n}](p_1,\ldots,p_{N_n})}={W[\phi^{n'}](p)\over
W[{\phi}^{n}](p)}\ .}

\vskip 6pt

\subsec{Relations among higher order theta derivatives and holomorphic
differentials}

By Riemann Vanishing Theorem it follows that
$$\theta(np+\c_{g-n}-y-\Delta)$$
$n\in I_g$, as a function of $y$, has a zero of order $n$ at $p$ for
all the effective divisors $\c_{g-n}$ of degree $g-n$. In
particular, \eqn\vanagain{
\sum_i\theta_i(p+\c_{g-2}-\Delta)\omega_i(p)=0\ . }

\vskip 6pt

\newth\vanteo{Theorem}{
Fix $x_1,\ldots,x_{g-1}\in C$. The following relations hold
$$\eqalignno{
&\sum_i^{\phantom{i_1\ldots
i_{g-1}}}\theta_i(x_1+\ldots+x_{g-1}-\Delta)\omega_i(x_1)=0\ ,\cr
&\sum_{i,j}^{\phantom{i_1\ldots
i_{g-1}}}\theta_{ij}(x_1+\ldots+x_{g-1}-\Delta)\omega_i(x_1)\omega_j(x_2)=0\
,\cr &\qquad\vdots\cr &\sum_{i_1,\ldots, i_{g-1}}\theta_{i_1\ldots
i_{g-1}}(x_1+\ldots+x_{g-1}-\Delta)\omega_{i_1}(x_1)\cdots\omega_{i_{g-1}}(x_{i_{g-1}})=0\
.}$$}

\vskip 6pt

\noindent {\bf Proof.} Without loss of generality, we can assume
distinct $x_1,\ldots,x_{g-1}$; the general case follows by
continuity arguments. The first relation is just Eq.\vanagain. Let
us assume that the equation
$$\sum_{i_1,\ldots, i_n}\theta_{i_1\ldots i_n}(x_1+\ldots+x_{g-1}-\Delta)
\omega_{i_1}(x_1)\ldots\omega_{i_n}(x_n)=0\ ,$$ holds, for all $n\in
I_{N-1}$, with $1<N\le g-1$. Then by taking its derivative with
respect to $x_{n+1}$ one obtains the subsequent relation.
\hfill$\square$

\vskip 6pt

\newth\ilcorol{Corollary}{
Fix $p\in C$ and a set of effective divisors $\c_k$, $k\in I_{g-2}$
of degree $k$. The following relations hold
$$\eqalignno{
&\sum_i^{\phantom{i_1,\ldots,
i_{g-1}}}\theta_i(p+\c_{g-2}-\Delta)\omega_i(p)=0\ ,\cr
&\sum_{i,j}^{\phantom{i_1,\ldots,
i_{g-1}}}\theta_{ij}(2p+\c_{g-3}-\Delta)\omega_i\omega_j(p)=0\ ,\cr
&\qquad\vdots\cr &\sum_{i_1,\ldots, i_{g-1}}\theta_{i_1\ldots
i_{g-1}}((g-1)p-\Delta)\omega_{i_1}\cdots\omega_{i_{g-1}}(p)=0\ .}$$
}

\vskip 6pt

\noindent We denote by $\lambda:=\{\lambda_1,\ldots,\lambda_l\}$ a
partition of length $|\lambda|:=l$ of some integer $d>0$, that is
$$\sum_{i=1}^l\lambda_i=d\ ,\qquad \lambda_1\ge\ldots\ge\lambda_l>0\ .$$
On the set of the partitions of an integer $d$, a total order
relation can be defined by setting
$$\lambda'>\lambda\qquad \Longleftrightarrow\qquad \exists i,\ 0<i\le \min\{|\lambda|,|\lambda'|\},\
\hbox{ s.t. } \left\{\matrix{\lambda'_j=\lambda_j, &1\le j< i\ ,\cr
\lambda'_i>\lambda_i\ .&}\right.$$ With respect to such a relation,
the minimal and the maximal partitions $\lambda^{min}$ and
$\lambda^{max}$ of $d$, are
$$\lambda^{min}_1=\ldots=\lambda^{min}_d=1\ ,\qquad \lambda^{max}_1=d\ .$$ Also observe that
$\lambda^{min}$ and $\lambda^{max}$ have, respectively, the maximal
and minimal lengths $|\lambda^{min}|=d$, $|\lambda^{max}|=1$.

\noindent For a general holomorphic $d$ differential $\eta$, let
$\eta(z)$ be its trivialization around a point $p\in C$, with
respect to some local coordinate $z$ and let us define
$$\eta^{(0)}(p):=\eta(z)\ ,\qquad \eta^{(n)}(p):= {\partial^n\eta\over \partial z^n}(z)\ ,
\qquad n> 0\ .$$

\vskip 6pt
\newth\anonimo{Theorem}{Fix $d\in I_{g-1}$, a point $p\in C$ and a effective divisor
$\c_{g-d}$ of degree $g-d$. Then, for each partition $\lambda$ of
$d$, there exists $c(\lambda)\in\ZZ$ independent of $C,p,\c_{g-d}$,
such that
\eqn\deriv{\eqalign{&\sum_{i_1,\ldots,i_l}^g\theta_{i_1\ldots
i_l}((d-1)p+\c_{g-d}-\Delta)
\omega_{i_1}^{(\lambda_1-1)}\cdots\omega_{i_l}^{(\lambda_l-1)}(p)\cr
&\quad= c(\lambda)\sum_{j_1,\ldots,j_d}^g\theta_{j_1\ldots
j_d}((d-1)p+\c_{g-d}-\Delta)\omega_{j_1}\cdots\omega_{j_d}(p)\ ,}}
where $l:=|\lambda|$. }

\vskip 6pt

\noindent {\bf Proof.} The theorem is just an identity for
$\lambda=\lambda^{min}$, with $c(\lambda^{min})=1$. Let us consider
a partition $\lambda>\lambda^{min}$ of $d$, and set $l:=|\lambda|<d$
($|\lambda|=d$ necessarily implies $\lambda=\lambda^{min}$). Fix
$\c=x_1+\ldots+x_{g-1}$, with $x_1,\ldots, x_{g-1}\in C$, and apply
the derivative operator
$$\D^{(\lambda)}:=
\Bigl({d\over dx_1}\Bigr)^{\lambda_1}\cdots \Bigl({d\over
dx_{l}}\Bigr)^{\lambda_{l}}\ ,$$ to the identity
\eqn\identity{\theta(\c-\Delta)=0\ .} Upon taking the limit
$x_1,\ldots,x_{l}\to p$, we obtain a sum, such that each term can be
associated to a partition $\lambda'$ of $d$ and written as
$$\sum_{i_1,\ldots,i_{l'}}^g\theta_{i_1\ldots i_{l'}}(lp+\c_{g-1-l}-\Delta)
\omega_{i_1}^{(\lambda'_1-1)}\cdots\omega_{i_{l'}}^{(\lambda'_{l'}-1)}(p)\
,
$$
with $l':=|\lambda'|$ and $\c_{g-1-l}=x_{l+1}+\ldots+x_{g-1}$. The
sum is over a set of partitions $\lambda'$ satisfying
$\lambda'\le\lambda$ and $l'\ge l$, so that $\lambda$ is the maximal
partition appearing. Thus, the sum can be rearranged as
\eqn\lambdaespr{\eqalign{ &\sum_{i_1,\ldots, i_{l}}\theta_{i_1\ldots
i_{l}}(lp+\c_{g-1-l}-\Delta)\omega^{(\lambda_1-1)}_{i_1}
\cdots\omega^{(\lambda_{l}-1)}_{i_{l}}(p)\cr
&=\sum_{\lambda'<\lambda}b(\lambda,\lambda')\sum_{i_1,\ldots,
i_{l'}} \theta_{i_1\ldots
i_{l'}}(lp+\c_{g-1-l}-\Delta)\omega^{(\lambda'_1-1)}_{i_1}\cdots\omega^{(\lambda'_{l'}-1)}_{i_{l'}}(p)\
, }} for some coefficients $b(\lambda,\lambda')\in\ZZ$. If the only
non-vanishing contribution to the RHS corresponds to
$\lambda'=\lambda^{min}$, the theorem follows after taking the limit
$x_{l+1},\ldots,x_{d-1}\to p$. Otherwise, for each
$\lambda'>\lambda^{min}$, one can obtain a further identity by
applying the operator $\D^{(\lambda')}$ to the identity \identity\
and taking the limit $x_1,\ldots,x_{l'}\to p$. This procedure leads
to an expression for
$$\sum_{i_1,\ldots, i_{l'}}\theta_{i_1\ldots
i_{l'}}(l'p+\c_{g-1-l'}-\Delta)\omega_{i_1}^{(\lambda'_1-1)}
\cdots\omega_{i_{l'}}^{(\lambda'_{l'}-1)}(p)\ ,$$ analogous to
Eq.\lambdaespr, where the RHS is a sum of terms corresponding to
partitions $\lambda''<\lambda'$. This expression can be used to
replace the term corresponding to $\lambda'$ in Eq.\lambdaespr,
considered in the limit $x_{l+1},\ldots,x_{l'}\to p$, with a sum
over a set of partitions $\lambda''<\lambda'$. After a finite number
of steps, the RHS of Eq.\lambdaespr\ reduces to a term corresponding
to $\lambda^{min}$ times an integer coefficient
$$c(\lambda):=\sum_{\lambda'<\lambda}\sum_{\lambda''<\lambda'}\ldots\sum_{\lambda^{\ldots}}
b(\lambda,\lambda')b(\lambda',\lambda'')\cdots
b(\lambda^{\ldots},\lambda^{min})\ .
$$ The arguments of the $\theta$-functions on both sides are
$$l'p-\c_{g-1-l'}-\Delta\ ,$$ where $l'$
is the length of the minimal partition $\lambda'>\lambda^{min}$
appearing in any intermediate step of the procedure. Therefore,
$l'\le d-1$ and the theorem follows. (Actually, with some more
effort, it can be proved that the bound $d-1$ cannot be improved).
\hfill$\square$

\vskip 6pt

\newth\anoncorol{Corollary}{Fix $d\in I_{g-1}$, a point $p\in C$ and an effective divisor
$\c_{g-d-1}$ of degree $g-d-1$. Then, for each partition $\lambda$
of $d$,
$$\sum_{i_1,\ldots,i_l}^g\theta_{i_1\ldots i_l}(dp+\c_{g-d-1}-\Delta)
\omega_{i_1}^{(\lambda_1-1)}\cdots\omega_{i_l}^{(\lambda_l-1)}(p)=0\
,
$$
where $l:=|\lambda|$. }

\vskip 6pt

\noindent {\bf Proof.} A trivial application of Eq.\deriv, with
$\c_{g-d}:=p+\c_{g-d-1}$, and Corollary \ilcorol. \hfill$\square$

\vskip 6pt

\newsec{Special loci in $C^g$ and the section $K$}\seclab\primecostr

In this section, we first introduce the basis of holomorphic
$1$-differentials (as a particular case of a definition which holds
for bases of holomorphic $n$-differentials, for arbitrary
$n\in\NN_+$), whose $n$-fold product can be used to construct bases of
holomorphic $n$-differentials. Next, we focus on the construction of
sets of two- and three-fold products of holomorphic abelian
differentials and discuss the condition under which they correspond
to bases of $H^0(K_C^2)$ and $H^0(K_C^3)$. We then introduce the
function $K$, which plays a key role in the present investigation,
in particular for what concerns the study of $\Theta_s$. Theorem
\quadgen\ summarizes some of the main results of the present
section, while Theorem \zeriK\ shows that $K$ counts the number of
intersections of special varieties on $J_0(C)$ defined in terms of
$\Theta_s$.

\vskip 6pt

\subsec{Duality between $N_n$-tuples of points and bases of
$H^0(K_C^n)$}

Let $C$ be a canonical curve of genus $g$ and let $C_d$, $d>0$, be
the set of effective divisors of degree $d$. Let $\{\eta\}_{i\in
I_g}$ be a basis of $H^0(K_C)$ and fix the divisor
$\c:=p_1+\ldots+p_{g-1}$ in such a way that the matrix
$[\eta_i(p_j)]_{^{i\in I_g}_{j\in I_{g-1}}}$ be of maximal rank. The
ratio
$$\sigma_{c}(p,q):={\det\eta(p,p_1,\ldots,p_{g-1})\over
 \det\eta(q,p_1,\ldots,p_{g-1})}\ ,$$
is a meromorphic function on $C_{g-1}$ and a meromorphic section of
the bundle $${\cal L}:=\pi^*_1K_C\otimes \pi^*_2 K_C^{-1}\ ,$$ on
$C\times C$, where $\pi_1$, $\pi_2$ are the projections of $C\times
C$ onto its first and second component, respectively. Note that
$\delta^*({\cal L})$, where $\delta:C\to C\times C$ is the diagonal
embedding $\delta(p):=(p,p)$, $p\in C$, is the trivial line bundle
on $C$. Furthermore, $\delta^*\sigma_c$ has neither zeros nor poles,
and $\sigma_c(p,p)/\sigma_{c'}(p,p)=1$, $\c,\c'\in C_{g-1}$. Hence,
there exists an isomorphism $H^0(\delta^*({\cal L}))\to \CC$ such
that $\sigma_{c}(p,p)\to 1$, for all $\c\in C_{g-1}$ for which
$\sigma_c$ is well-defined.

\vskip 6pt

\newth\thnewbasis{Proposition}{Fix $n\in\ZZ_{>0}$ and let $p_1,\ldots,p_{N_n}$ be a set of points
of $C$ such that $$\det\,\phi^n(p_1,\ldots,p_{N_n})\ne 0\ ,$$
with $\{\phi^n_i\}_{i\in I_{N_n}}$ an arbitrary basis of $H^0(K_C^n)$. Then
\eqn\basendiff{\gamma^n_i(z):=
{\det\,\phi^n(p_1,\ldots,p_{i-1},z,p_{i+1},\ldots,p_{N_n})\over
 \det\,\phi^n(p_1,\ldots,p_{N_n})}\ ,} $i\in I_{N_n}$, for all $z\in C$, is a basis of
 $H^0(K_C^n)$ which is independent of the choice of the basis
 $\{\phi^n_i\}_{i\in I_{N_n}}$ and, up to normalization, on the local coordinates on
 $C$.}

\vskip 6 pt

\noindent {\bf Proof.} Since the matrix
$[\phi^n]_{ij}:=\phi^n_i(p_j)$ is non-singular, it follows that
\eqn\remember{\gamma^n_i=\sum_{j=1}^{N_n}[\phi^n]^{-1}_{ij}\phi^n_j\
,} $i\in I_{N_n}$, is a basis of $H^0(K_C^n)$.\hfill$\square$

\vskip 6pt

\noindent Observe that \eqn\detdet{\det\,
\gamma_i^n(p_1,\ldots,p_{j-1},z,p_{j+1},\ldots,p_{N_n})=\gamma_j^n(z)
\ ,} for all $z\in C$, and
\eqn\centralle{\gamma_i^n(p_j)=\delta_{ij}\ ,} $i,j\in I_{N_n}$.
Furthermore, \eqn\detdetover{\det\, \gamma_i^n(z_j)={\det\,
\phi_i^n(z_j)\over\det\, \phi_i^n(p_j)}\ ,} for all
$z_1,\ldots,z_g\in C$. For $n=1$, for each choice of
$p_1,\ldots,p_g\in C$ with $\det\,\eta_i(p_j)\ne 0$, we set
\eqn\newbasis{\sigma_i(z):=\gamma^1_i(z)\ ,\qquad i\in I_g\ .}

\vskip 6pt

The bases $\{\gamma^n_i\}_{i\in I_{N_n}}$, in some sense, can be considered as dual to the $N_n$-tuple of points
$p_1,\ldots,p_{N_n}$ involved in their definition. In facts, each point $p\in C$ corresponds to the
element of
$\PP(H^0(K_C^n)^*)$ given by
\eqn\pdual{p[\phi]:=\phi(p)\ ,\qquad \phi\in H^0(K_C^n)\ .}
By Eq.\pdual, we mean that a generic representative of $p$ in $H^0(K_C^n)^*$ is given by $\phi\to\phi(p)$, where
$\phi(p)$ is determined with respect to the choice of a local trivialization and a local coordinate
around $p$. Hence, the choice of a $N_n$-tuple of points $p_1,\ldots,p_{N_n}\in C$ satisfying the
conditions of Proposition \thnewbasis, together with the choice of local trivializations of $K_C$ and
coordinates around each $p_i$, $i\in I_{N_n}$, determines a basis of $H^0(K_C^n)^*$. On the other
hand, by Proposition \thnewbasis, such data also determine the basis $\{\gamma_i^n\}_{i\in I_N}$ of
$H^0(K_C)$, which satisfies the property
$$p_i[\gamma_j^n]=\delta_{ij}\ ,\qquad \forall i,j\in I_{N_n}\ ,
$$
so that
$$p_i\equiv \gamma^n_i{}^*\ ,$$
where $\{\gamma^n_i{}^*\}_{i\in I_{N_n}}$ if the basis dual to $\{\gamma_i^n\}_{i\in I_N}$.

Such a construction also extends to symmetric products of
$H^0(K^n_C)$; we will only consider the case of $\Sym^2(H^0(K_C))$,
but the generalizations are straightforward. Each basis
$\{\eta_i\}_{i\in I_g}$ of $H^0(K_C)$ naturally defines a basis in
$\Sym^2(H^0(K_C))$, given by
$\spbase{\eta}_i:=\sprod{\eta_{\1_i}}{\eta_{\2_i}}$, $i\in I_M$.
Similarly, each element $p+q\in C_2\cong \Sym^2 C$ corresponds to an
element $p\cdot q\in\PP(\Sym^2(H^0(K_C))^*)$ by \eqn\pqdef{(p\cdot
q)\bigl[\smsum_k
\eta_k\cdot\rho_k\bigr]:=\sum_k(\eta_k(p)\rho_k(q)+\eta_k(q)\rho_k(p))\
, \qquad\smsum_k \eta_k\cdot\rho_k\in\Sym^2H^0(K_C)\ ,} with the
same notation of Eq.\pdual. In this sense, with respect to an
arbitrary choice of local coordinates around $p_1,\ldots,p_g$, we
have \eqn\pqdual{p_i[\sigma_j]=\delta_{ij}\ ,\qquad (p\cdot
p)_k[\sigma\cdot\sigma_l]=\chi_k\delta_{kl}\ ,} $i,j\in I_g$,
$k,l\in I_M$, so that
$$p_i\equiv \sigma_i^*\ ,\qquad \chi_k^{-1}(p\cdot p)_k\equiv (\sigma\cdot\sigma)_k^*\ ,
\qquad \forall i\in I_g,\ k\in I_M\ ,$$
where $\{\sigma_i^*\}_{i\in I_g}$ and $\{(\sigma\cdot\sigma)^*_k\}_{k\in I_M}$
are the dual bases of $\{\sigma_i\}_{i\in I_g}$ and
$\{\sigma\cdot\sigma_k\}_{k\in I_M}$, respectively.

\vskip 6pt

The results of section \Riemanndef\ can be used to derive an explicit expression for the matrix
$[\omega]^{-1}_{ij}$, with $\{\omega_i\}_{i\in I_g}$ the dual basis of the symplectic basis of
$H_1(C,\ZZ)$.

\vskip 6pt

\newrem\defAB{Definition}{For each fixed $g$-tuple
$(p_1,\ldots,p_g)\in C^g$ let us define the following
effective divisors
$$\a:=\sum_{j\in I_g}p_j\ ,\qquad \a_i:=\a-p_i\ ,\qquad \b:=\a-p_1-p_2\ ,
$$ $i\in I_g$. Define the subset of $C^g$
$$\A:=\{(p_1,\ldots,p_g)\in C^g\mid\det\eta_i(p_j)=0\}\ ,
$$ with $\{\eta_i\}_{i\in I_g}$ an arbitrary basis of $H^0(K_C)$. }

\vskip 6pt

Fix $g+1$ arbitrary points $p_1,\ldots,p_g,z\in C$. By taking the
limit $y\to z$ in Eq.\dettheta, we obtain
\eqn\ildet{\det\eta(z,p_1,\ldots,\check p_i,\ldots,p_g)
=\kuno[\eta]\sum_{l=1}^g\theta_{\Delta,l}(\a_i)
\omega_l(z)\prod_{{j,k\neq i \atop j<k}}E(p_j,p_k)\prod_{j\neq
i}\sigma(p_j)\ ,} for all $i\in I_g$, where
$\theta_{\Delta,i}(e):=\partial_{z_i}\theta_\Delta(z)_{|z=e}$,
$e\in\CC^g$. Furthermore,
\eqn\ilwronsk{\theta(\K^w+w-z)={\sigma(z)E(z,w)^g\over\kappa[\omega]\sigma(w)^g}W[\omega](w)\
,} for all $w,z\in C$, with $W[\omega](z)$ the Wronskian of
$\{\omega_i\}_{i\in I_g}$ at $z$. Note that, by \ildet, the
condition $(p_1,\ldots,p_g)\in C^g\setminus\A$ implies
\eqn\connnd{\sum_j\theta_{\Delta,j}(\a_i)\omega_j(p_i)\neq 0\ ,} for
all $i\in I_g$.

\newth\inverr{Proposition}{Fix $(p_1,\ldots,p_g)\in C^g\setminus \A$, with $\A$ defined in \defAB. Set $[\omega]_{ij}:=\omega_i(p_j)$. We have
\eqn\leomega{[\omega]^{-1}_{ij}=\oint_{\alpha_j}\sigma_i={\theta_{\Delta,j}
\left(\a_i\right)\over\sum_k\theta_{\Delta,k}\left(\a_i\right)\omega_k(p_i)}\
,} $i,j\in I_g$, so that
\eqn\siggmas{\sigma_i(z)=\sigma(z,p_i){\deltadiv(\a+z-y-p_i)\over\deltadiv(\a-y)E(z,p_i)}{E(y,p_i)\over
E(y,z)} \prod_{1}^gE(z,p_j)\prod_{j\neq i}{1\over E(p_i,p_j)} \ ,}
and \eqn\kkkk{\kuno[\sigma]={\sigma(y)\prod_1^gE(y,p_i)\over
\deltadiv(\a-y)\prod_{i<j}^gE(p_i,p_j)\prod_1^g\sigma(p_k)} \ ,} for
all $z,y,x_i,y_i\in C$, $i\in I_g$, with $\a,\a_i$ as in Definition
\defAB. Furthermore, fix $p_1,\ldots,p_{N_n}\in C$ such that
$\det\,\phi^n(p_1,\ldots,p_{N_n})\ne 0$, with $\{\phi^n_i\}_{i\in
I_{N_n}}$ an arbitrary basis of $H^0(K_C^n)$. Then,
\eqn\glindiff{\gamma^n_i(z)=\sigma(z,p_i)^{2n-1}{\deltadiv\bigl(\sum_1^{N_n}p_j+z-p_i\bigr)\prod_{^{j=1}_{j\neq
i}}^{N_n}E(z,p_j)\over\deltadiv\bigl(\sum_1^{N_n}p_j\bigr)\prod_{^{j=1}_{j\neq
i}}^{N_n}E(p_i,p_j)}\ , } $i\in I_{N_n}$, and
\eqn\knphi{\kenne[\gamma^n]={1\over
\deltadiv\bigl(\sum_1^{N_n}p_i\bigr)
\prod_{1}^{N_n}\sigma(p_i)^{2n-1}\prod_{^{i,j=1}_{i<j}}^{N_n}
E(p_i,p_j)}\ .}}

\vskip 6pt

\noindent {\bf Proof.}
By \remember\ and \newbasis\ we have
$\sigma_i=\sum_j[\omega]^{-1}_{ij}\omega_j$, and \leomega\ follows by \dettheta\ and \ildet.
Eqs.\siggmas\kkkk\ follow by \dettheta\ and by $\det\sigma_i(p_j)=1$, respectively. Similarly,
\glindiff\ follows by \basendiff\ and \detthetaii. Eq.\knphi\ follows by $\det\gamma^n_i(p_j)=1$. \hfill$\square$

\vskip 6pt

\newth\surprise{Corollary}{Fix $(p_1,\ldots,p_g)\in C^g\setminus \A$. Then
$$\sum_{i\in I_g}{\theta_{\Delta,j}
\left(\a_i\right)\over\sum_l\theta_{\Delta,l}\left(\a_i\right)\omega_l(p_i)}\omega_k(p_i)=\delta_{jk}\
,$$
$j,k\in I_g$.}

\vskip 6pt

\noindent {\bf Proof.} This is just the identity $\sum_{i\in
I_g}[\omega]^{-1}_{ij}[\omega]_{ki}=\delta_{jk}$ with
$[\omega]^{-1}_{ij}$ given by \leomega.\hfill$\square$

\vskip 6pt

\subsec{Special loci in $C^g$ from the condition of linear
independence for holomorphic differentials}

There exist natural homomorphisms from $\Sym^n(H^0(K_C))$ to
$H^0(K_C^n)$, which, for $n=2$, we denote by
$$\eqalign{\psi\colon\Sym^2(H^0(K_C^2))&\to
H^0(K_C^2) \cr \sprod{\eta}{\rho}\quad&\mapsto\eta\rho\ .}$$ By Max
Noether's Theorem, if $C$ is a Riemann surface of genus two or
non-hyperelliptic with $g\ge 3$, then $\psi$ is surjective. Set
\eqn\lev{v_i:=\psi\spbase{\sigma}_i=\sigma_{\1_i}\sigma_{\2_i}\ ,}
$i\in I_M$, so that \eqn\devv{v_i(p_j)=
\left\{\vcenter{\vbox{\halign{\strut\hskip 6pt $ # $ \hfil & \hskip
2cm$ # $ \hfil\cr \delta_{ij}\ , & i\in I_g \ ,\cr 0\ ,& g+1\le i\le
M\ ,\cr}}}\right.} $j\in I_g$. By dimensional reasons, it follows
that for $g=2$ and $g=3$ in the non-hyperelliptic case, the set
$\{v_i\}_{i\in I_N}$ is a basis of $H^0(K_C^2)$ if and only if
$\{\sigma_i\}_{i\in I_g}$ is a basis of $H^0(K_C)$. On the other
hand, for $g\ge 3$ in the hyperelliptic case, there exist
holomorphic quadratic differentials which cannot be expressed as
linear combinations of products of elements of $H^0(K_C)$, so that
$v_1,\ldots,v_N$ are not linearly independent. The other
possibilities are considered in the following proposition.

\vskip 6pt

\newth\thlev{Proposition}{Fix the points $p_1,\ldots,p_g\in C$, with $C$ non-hyperelliptic of genus
$g\ge 4$. If the following conditions are
satisfied
\smallskip
\item{\it i.} $\det\eta_i(p_j)\ne 0$, with $\{\eta_i\}_{i\in I_g}$ an arbitrary basis of $H^0(K_C)$;
\smallskip
\item{\it ii.} $\b:=\sum_{i=3}^gp_i$ is the greatest common divisor of
$(\sigma_1)$ and $(\sigma_2)$, with $\{\sigma_i\}_{i\in I_g}$
defined in \newbasis,
\smallskip
\noindent then
$\{v_i\}_{i\in I_N}$ is a basis of $H^0(K_C^2)$.
\smallskip
\noindent Conversely, if there exists a set
$\{\hat\sigma_i\}_{i\in I_g}$ of holomorphic $1$-differentials, such that
\smallskip
\item{\it a.} $i\neq j\Rightarrow\hat\sigma_i(p_j)=0$, for all $i,j\in I_g$;
\smallskip
\item{\it b.} $\{\hat v_i\}_{i\in I_N}$ is a basis of $H^0(K_C^2)$, with $\hat
v_i:=\hat\sigma\hat\sigma_i$, $i\in I_N$;
\smallskip
\noindent then {\it i}) and {\it ii}) hold.
}

\vskip 6pt

\noindent {\bf Proof.} To prove that {\it i}) and {\it ii}) imply that $\{v_i\}_{i\in I_N}$ is a
basis of $H^0(K_C^2)$, we first prove that $\sigma_i$ is the unique $1$-differential, up to normalization,
vanishing at $\c_i:=(\sigma_i)-\b$, $i=1,2$. Any
$1$-differential $\sigma'_i\in H^0(K_C)$ vanishing at $\c_i$ corresponds to an element $\sigma'_i/\sigma_i$ of
$H^0(\O(\b))$, the space of meromorphic functions $f$ on $C$ such that $(f)+\b$ is an effective divisor. Suppose that
there exists a $\sigma'_i$ such that $\sigma'_i/\sigma_i$ is not a constant, so that $h^0(\O(\b))\ge 2$. By the
Riemann-Roch Theorem
$$h^0(K_C\otimes \O(-\b))=h^0(\O(\b))-\deg \b-1+g\ge 3\ ,$$
there exist at least $3$ linearly independent $1$-differentials
vanishing at the support of $\b$ and, in particular, there exists a linear
combination of such differentials vanishing at $p_1,\ldots,p_g$.
This implies that $\det\eta_i(p_j)= 0$, with
$\{\eta_i\}_{i\in I_g}$ an arbitrary basis of $H^0(K_C)$, contradicting the hypotheses.
Fix $\zeta_i,\zeta_{1i},\zeta_{2i}\in \CC$ in such a way that
$$\sum_{i=3}^g\zeta_i\sigma_i^2+\sum_{i=1}^g\zeta_{2i}
\sigma_1\sigma_i+\sum_{i=2}^g\zeta_{1i}\sigma_2\sigma_i=0\ .$$ Evaluating this
relation at the point $p_j$, $3\le j\le g$ yields
$\zeta_j=0$.
Set
\eqn\leti{t_1:=-\sum_{j=2}^g\zeta_{1j}\sigma_j\ ,\quad
t_2:=\sum_{j=1}^g\zeta_{2j}\sigma_j\ ,}
so that $\sigma_1t_2=\sigma_2 t_1$.
Since the supports of $\c_1$ and $\c_2$ are disjoint,
$t_i$ must be an element of $H^0(K_C\otimes\O(-\c_i))$, $i=1,2$ and then, by the previous remarks,
$t_1/\sigma_1=t_2/\sigma_2=\zeta\in\CC$. By \leti\
$$\zeta\sigma_1+\sum_{j=2}^g\zeta_{1j}\sigma_j=0\ ,\quad
\zeta\sigma_2-\sum_{k=1}^g\zeta_{2k}\sigma_k=0\
,$$ and, by linear independence of $\sigma_1,\ldots,\sigma_g$, it follows that
$\zeta=\zeta_{1j}=\zeta_{2k}=0$,
$2\le j\le g$, $k\in I_g$.

\vskip 3pt

\noindent Let us now assume that {\it a}) and {\it b}) hold for some set $\{\hat\sigma_i\}_{i\in
I_g}$. Then
$\{\hat\sigma_i\}_{i\in I_g}$ is a basis of $H^0(K_C)$ if and only if $\det\eta_i(p_j)\neq 0$. If
$\{\hat\sigma_i\}_{i\in I_g}$ is not a basis of $H^0(K_C)$, the corresponding $\hat
v_i$, $i\in I_N$, cannot span a $N$-dimensional vector space. Then {\it i}) is satisfied and the
basis $\{\hat\sigma_i\}_{i\in I_g}$ corresponds, up to a non-singular diagonal transformation, to
the basis $\{\sigma_i\}_{i\in I_g}$, defined in \newbasis.

Without loss of generality, to prove {\it
ii}) we can assume
that $\hat\sigma_i\equiv\sigma_i$, $i\in I_g$ and then $\hat v_i\equiv v_i$, $i\in I_N$. Suppose
there exists $p\in C$ such that $p+\b\leq (\sigma_i)$,
for all $i\in I_2$. If $p\equiv p_1$ or $p\equiv p_2$, then $\sigma_i(p)=0$, for all $i\in I_g$,
and therefore $\{\sigma_i\}_{i\in I_g}$ would not be a basis, which
contradicts {\it b}).

Suppose there exists $i$, $3\le i\le g$, with
$p\equiv p_i$. In this case, each $v_j$, $j\in I_N\setminus \{i\}$, has a double zero in $p_i$,
whereas $v_i(p_i)\neq 0$; therefore, an element of $H^0(K_C^2)$ with a single zero in $p_i$ (such
as, for example,
$\sigma_i\sigma_j$, with $3\le j\le g$, $j\neq i$) cannot be expressed as a linear combination of
$v_1,\ldots,v_N$, in contradiction with the assumptions.

Finally, suppose that $p\neq p_i$, for all $i\in I_g$. In this case,
there exists at least one $\sigma_i$, $3\le i\le g$, with
$\sigma_i(p)\neq 0$, since, on the contrary, $\{\sigma_i\}_{i\in
I_g}$ would not be a basis of $H^0(K_C)$. Suppose that
$\sigma_i(p)\neq 0$ and $\sigma_j(p)\neq 0$ for some $3\le i,j\le
g$, $i\neq j$. Then $\sigma_i\sigma_j$ cannot be expressed as a
linear combination of $v_k$, $k\in I_N$. In fact,
$\sigma_i\sigma_j(p_k)=0$, for all $k\in I_g$, would imply that
$\sigma_i\sigma_j=\sigma_1\rho_1+\sigma_2\rho_2$, for some
$\rho_1,\rho_2\in H^0(K_C)$; but this is impossible, since
$\sigma_1(p)=0=\sigma_2(p)$, whereas $\sigma_i\sigma_j(p)\neq 0$.
Therefore, there should exist exactly one $i\in I_g$ with
$\sigma_i(p)\neq 0$. It follows that $\sigma_j(p)=0=\sigma_j(p_i)$,
for all $j\in I_g\setminus \{i\}$; then
$h^0(K_C\otimes\O(-p-p_j))\ge g-1$ and, by Riemann-Roch Theorem,
there exists a non-constant meromorphic function on $C$, with only
single poles in $p$ and $p_j$. But this would imply that $C$ is
hyperelliptic, in contradiction with the hypotheses. \hfill$\square$

\vskip 6pt

The proof that {\it i}) and {\it ii}) imply that $\{v_i\}_{i\in
I_N}$ is a basis is due to Petri \petriuno\ (see also
\ref\ottimo{E.~Arbarello, M.~Cornalba, P.A.~Griffiths and J.~Harris,
{\it Geometry of algebraic curves, I,} Grundlehren Math. Wiss., vol.
267, Springer-Verlag, 1985.}). It can be proved that on a
non-hyperelliptic curve there always exists a set of points
$\{p_1,\ldots,p_g\}$ satisfying the hypotheses of Proposition
\thlev. This is related to the classical result $\dim\Theta_s=g-4$
for non-hyperelliptic surfaces of genus $g\ge 4$, as will be shown
in Corollary \futuref{thdimthsing}.

In view of Theorem \thlev, it is useful to introduce the following
subset of $C^g\equiv\underbrace{C\times \ldots \times C}_{g\hbox{
times}}$.

\newrem\defB{Definition}{Let $\B$ be the subset of $C^g$
$$\B:=\{(p_1,\ldots,p_g)\in C^g\mid\det\eta_i(p_j)=0\vel
\gcd((\sigma_1),(\sigma_2))\neq \b \}\ ,$$ for an arbitrary basis
$\{\eta_i\}_{i\in I_g}$ of $H^0(K_C)$. }

\vskip 6pt

\newth\rango{Corollary}{Fix $(p_1,\ldots,p_g)\in C^g\setminus\A$ such that the greatest common divisor
of $(\sigma_1)$ and $(\sigma_2)$ be $\b+q_1+\ldots+q_n$, for some $q_1,\ldots,q_n\in C$, $n\ge 1$. Then the dimension
$r$ of the vector space generated by $\{v_i\}_{i\in I_N}$ is $r=N-n$.}

\vskip 6pt

\noindent {\bf Proof.} Let us prove that $n$ is the number ($N-r$)
of independent linear relations among $v_1,\ldots,v_N$. Set
$\d:=q_1+\ldots +q_n$. By $\det\eta_i(p_j)\neq 0$, the quadratic
differentials $\sigma_i^2$, $i\in I_g$, are linearly independent and
independent of $\sigma_1\sigma_2,\sigma_1\sigma_i,\sigma_2\sigma_i$,
$i\in I_g\setminus \{1,2\}$. Therefore, all the independent linear
relations have the form \eqn\lineare{\sigma_1 t_2=\sigma_2 t_1\ ,}
for some $t_1,t_2\in H^0(K_C)$, with the condition $t_1(p_1)=0$ in
order to exclude the trivial relation $t_i=\sigma_i$, $i=1,2$.
Consider the effective divisors $\hat\c_1,\hat\c_2$ of degree $g-n$
with no common points, defined by $\hat\c_i:=(\sigma_i)-\d-\b$,
$i=1,2$. By $\det\eta_i(p_j)\ne0$, it follows that
$h^0(K_C\otimes\O(-\b))=2$, so that $h^0(K_C\otimes\O(-\b-\d))=2$
too. This implies that $\sigma_1/\sigma_2$ and $\sigma_2/\sigma_1$
are the unique elements of $H^0(\O(\hat\c_1))$ and
$H^0(\O(\hat\c_2))$, respectively. Then, by Riemann-Roch Theorem, we
have $h^0(K_C\otimes\O(-\hat\c_i))=n+1$, $i=1,2$. By Eq.\lineare,
the divisors of $t_1$, $t_2$ satisfy
$$\hat\c_1+(t_2)=\hat\c_2+(t_1)\ ,$$
so that $t_i\in H^0(K_C\otimes \O(-\hat c_i))$. In particular, a basis $\sigma_1,\alpha_1,\ldots,\alpha_n$
of $H^0(K_C\otimes \O(-\hat c_1))$ can be chosen in such a way that $\alpha_i(p_1)=0$, for all $i\in I_n$. Hence, $t_1$ is a linear combination
of $\alpha_1,\ldots,\alpha_n$ and there are at most $n$ linearly independent relations of the form
\lineare. This implies $N-r\le n$.

Let us now prove that such $n$ linearly independent relations exist. By the Riemann-Roch Theorem, since
$h^0(K_C\otimes\O(-\b-\d))=2$, we obtain $h^0(\O(\b+\d))=n+1$; a basis for $H^0(\O(\b+\d))$ is given by
$\alpha_1/\sigma_1,\ldots,\alpha_n/\sigma_1$ and the constant function. On the other hand, if
$\sigma_2,\beta_1,\ldots,\beta_n$ is a basis for $H^0(K_C\otimes\O(-\hat\c_2))$, then
$\beta_1/\sigma_2,\ldots,\beta_n/\sigma_2$ are $n$ linearly independent elements of
$H^0(\O(\b+\d))$. Hence, there exist $n$ linearly independent relations
$${\beta_i\over \sigma_2}=\sum_{j=1}^nc_{ij}{\alpha_j\over\sigma_1}+c_{i0}\ ,
$$
$i\in I_n$, for some $c_{ij}\in\CC$, $0\le j\le n$. By multiplying both sides by $\sigma_1\sigma_2$, we obtain
$$\sigma_1\beta_i=\sum_{j=1}^nc_{ij}\sigma_2\alpha_j+c_{i0}\sigma_1\sigma_2\ .$$
Therefore, $N-r\ge n$ and the corollary follows.
\hfill$\square$

\vskip 6pt

\noindent Consider the holomorphic $3$-differentials (with the notation defined in section \notation)
\eqn\trebasis{\varphi_i=\sigma\sigma\sigma_i:=\sigma_{\1_i}\sigma_{\2_i}
\sigma_{\3_i}\ ,} $i\in I_{M_3}$, with
$\{\sigma_i\}_{i\in I_g}$ a basis of $H^0(K_C)$. By
the Max Noether's Theorem and dimensional reasons, it follows that the first
$N_3:=5g-5$ of such differentials are a basis of $H^0(K_C^3)$ for $g=3$ in the non-hyperelliptic
case,
whereas they are not linearly independent for $g\ge 2$ in the hyperelliptic case. The other possibilities
are considered in the following proposition.

\vskip 6pt

\newth\thtrebase{Proposition}{Fix the points $p_1,\ldots,p_g\in C$, with $C$ non-hyperelliptic of genus $g\ge 4$.
If the following conditions are satisfied for a fixed $i\in
I_g\setminus\{1,2\}$:
\smallskip
\item{\it i. } $\det\eta_j(p_k)\neq 0$, with $\{\eta_j\}_{j\in I_g}$ an arbitrary basis of $H^0(K_C)$;
\smallskip
\item{\it ii. } $\b:=\sum_{j=3}^gp_j$ is the greatest common divisor of $(\sigma_1)$ and $(\sigma_2)$,
with $\{\sigma_j\}_{j\in I_g}$ defined in \newbasis;
\smallskip
\item{\it iii. } $p_k$ is a single zero for $\sigma_1$, for all $k\neq i$, $3\le k\le g$;
\smallskip
\noindent then the set $\{\varphi_j\}_{j\in I_{N_3-1}}\cup
\{\varphi_{i+5g-8}\}$ is a basis of $H^0(K_C^3)$. \noindent In
particular, if {\it i}), {\it ii}) and
\smallskip
\item{\it iii'.} $p_3,\ldots,p_g$ are single zeros for $\sigma_1$,
\smallskip
\noindent are satisfied, then, for each $i$, $3\le i\le g$, the set $\{\varphi_j\}_{j\in
I_{N_3-1}}\cup \{\varphi_{i+5g-8}\}$ is a basis of $H^0(K_C^3)$.
\smallskip
\noindent Conversely, if for some fixed $i\in I_g\setminus\{1,2\}$ there exists a set $\{\hat\sigma_j\}_{j\in
I_g}$ of holomorphic $1$-differentials, such that
\smallskip
\item{\it a.} $j\neq k\Rightarrow\hat\sigma_j(p_k)=0$, for all $j,k\in I_g$;
\smallskip
\item{\it b.} $\{\hat\varphi_j\}_{j\in
I_{N_3-1}}\cup \{\hat\varphi_{i+5g-8}\}$ is a basis of $H^0(K_C^3)$, with $\hat
\varphi_j:=\hat\sigma\hat\sigma\hat\sigma_j$, $j\in I_{M_3}$;
\smallskip
\noindent then {\it i}), {\it ii}) and {\it iii}) hold.}

\vskip 6pt

\noindent {\bf Proof.} We first prove that if {\it i}), {\it ii})
and {\it iii}) hold for a fixed $i$, $3\le i\le g$, then
$\{\varphi_j\}_{j\in I_{N_3-1}}\cup \{\varphi_{i+5g-8}\}$ is a basis
of $H^0(K_C^3)$. To this end it is sufficient to prove that the
equation
$$\sum_{j=3}^g(\zeta_j\sigma_j^3+\zeta_{1j}\sigma_1\sigma_j^2+\zeta_{12j}\sigma_1\sigma_2\sigma_j)
+\sigma_1^2\mu+\sigma_2^2\nu +\zeta_{2i}\sigma_2\sigma_i^2=0\ ,$$ is
satisfied if and only if
$\zeta_j,\zeta_{1j},\zeta_{2i},\zeta_{12j}\in\CC$, $3\le j\le g$,
and $\mu,\nu \in H^0(K_C)$ all vanish identically (no non-trivial
solution). Evaluating such an equation at $p_j\in C$, $3\le j\le g$,
gives $\zeta_j=0$. Furthermore, note that, by condition {\it iii}),
for each $j\neq i$, $3\le j\le g$, $\sigma_1\sigma_j^2$ is the
unique $3$-differential with a single zero in $p_j$, so that
$\zeta_{1j}=0$. We are left with
\eqn\conditio{\zeta_{1i}\sigma_1\sigma_i^2+\zeta_{2i}\sigma_2\sigma_i^2+\sigma_1^2\mu+
\sigma_2^2\nu +\sum_{j=3}^g \zeta_{12j}\sigma_1\sigma_2\sigma_j=0\
.} By Riemann-Roch Theorem, for each $k$, $3\le k\le g$,
$h^0(K_C\otimes\O(-\b-p_k))\ge 1$; the condition {\it ii}) implies
that $h^0(K_C\otimes\O(-\b-p_k))\le 1$, so that, in particular,
there exists a unique (up to a constant) non-vanishing $\beta$ in
$H^0(K_C\otimes\O(-\b-p_i))$. Furthermore,
$$H^0(K_C\otimes\O(-\b))\neq \bigcup_{k=3}^gH^0(K_C\otimes\O(-\b-p_k))\ ,$$
because the LHS is a $2$-dimensional space and the RHS is a finite
union of $1$-dimensional subspaces; then, there exists $\alpha\in
H^0(K_C\otimes\O(-\b))$ such that $p_3,\ldots,p_g$ are single zeros
for $\alpha$. Note that $\alpha$ and $\beta$ span
$H^0(K_C\otimes\O(-\b))$ and $\alpha^2$, $\beta^2$ and $\alpha\beta$
span $H^0(K^2_C\otimes\O(-2\b))$. Hence, the existence of
non-trivial $\zeta_{1i},\zeta_{2i},\zeta_{12j},\nu,\mu$ satisfying
Eq.\conditio\ is equivalent to the existence of non-trivial
$\nu',\mu'\in H^0(K_C)$ and
$\zeta_{\alpha},\zeta_{\beta},\zeta_{\alpha\beta j}\in \CC$
satisfying
$$\zeta_{\alpha}\alpha\sigma_i^2+\zeta_{\beta}\beta\sigma_i^2+\alpha^2\mu'+\beta^2\nu' +\sum_{j=3}^g
\zeta_{\alpha\beta j}\alpha\beta\sigma_j=0\ .$$ Note that
$\alpha\sigma_i^2$ is the unique $3$-differential with a single zero
in $p_i$, so that $\zeta_{\alpha}=0$. Condition {\it ii}) implies
that $\b$ is the greatest common divisor of $(\alpha)$ and
$(\beta)$. Then $\alpha\neq 0$ on the support of $\c_\beta$, where
$\c_\beta:=(\beta)-\b-p_i$. Hence, $\mu'\in
H^0(K_C\otimes\O(-\c_{\beta}))$, which, by Riemann-Roch Theorem, is
a $1$-dimensional space, so that $\mu'=\zeta'_\mu\beta$, for some
$\zeta'_\mu\in\CC$. Since, by construction, $\beta\neq 0$, we have
$$\zeta_{\beta}\sigma_i^2+\zeta'_\mu\alpha^2+\beta\nu' +\sum_{j=3}^g \zeta_{\alpha\beta j}\alpha\sigma_j=0\ .$$
By evaluating such an equation at $p_i$ gives $\zeta_\beta=0$.
Furthermore, since $\beta\neq 0$ on the support of $\c_\alpha$,
where $\c_\alpha:=(\alpha)-\b$, it follows that
$\nu'=\zeta'_\nu\alpha$, for some $\zeta'_\nu\in\CC$. Since
$\alpha\neq 0$
$$\zeta'_\mu\alpha+\zeta'_\nu\beta +\sum_{j=3}^g \zeta_{\alpha\beta j}\sigma_j=0\ ,$$
which implies that $\zeta'_\mu=\zeta'_\nu=\zeta_{\alpha\beta j}=0$,
for all $3\le j\le g$.

\vskip 3pt

\noindent Conversely, suppose that {\it a}) and {\it b}) hold for
some fixed $i$, with $3\le i\le g$, and for some set
$\{\hat\sigma_j\}_{j\in I_g}$. If $\det\eta_j(p_k)=0$, then
$\{\hat\sigma_j\}_{j\in I_g}$ is not a basis of $H^0(K_C)$ and
$\{\hat\varphi_j\}_{j\in I_{N_3-1}}$ cannot span a
$(N_3-1)$-dimensional vector space. Then {\it i}) is satisfied and
the basis $\{\hat\sigma_j\}_{j\in I_g}$ corresponds, up to a
non-singular diagonal transformation, to the basis
$\{\sigma_j\}_{j\in I_g}$, defined in \newbasis.

\noindent Without loss of generality, we can prove {\it ii}) and
{\it iii}) for $\hat\sigma_j\equiv\sigma_j$, $j\in I_g$ and then
$\hat \phi_j\equiv \phi_j$, $j\in I_{M_3}$. Since the
$3$-differentials $\sigma_1 v_j$, $j\in I_N$, are distinct elements
of a basis of $H^0(K_C^3)$, then $v_j$, $j\in I_N$, are linearly
independent elements of $H^0(K_C^2)$ and, by Proposition \thlev,
also condition {\it ii}) is satisfied.

\noindent Finally, assume that there exists $k\neq i$, $3\le k\le
g$, such that $\sigma_1$ has a double zero in $p_k$. Then, apart
from $\varphi_k\equiv\sigma_k^3$, which satisfies
$\varphi_k(p_k)\neq 0$, all the other $3$-differentials of the basis
have a double zero in $p_k$. Therefore, an element of $H^0(K^3_C)$
with a single zero in $p_k$ cannot be a linear combination of the
elements of such a basis, which is absurd. (An example of a
holomorphic $3$-differential with a single zero in $p_k$ is
$\sigma_2\sigma_k^2$, since, by condition {\it ii}), $\sigma_2$
cannot have a double zero in $p_k$). \hfill$\square$

\vskip 6pt

\subsec{Combinatorial lemmas and determinants of holomorphic
differentials}

Applying Lemmas \thcombi\ and \thcombvi\ to determinants of
symmetric products of holomorphic $1$-differentials on an algebraic
curve $C$ of genus $g$ leads to combinatorial relations. By
Eq.\dettheta\ and \detthetaii, such combinatorial relations yield
non trivial identities among products of theta functions.

\vskip 6pt

\newth\thdetvan{Proposition}{The following identities \eqn\detdue{\det \eta\eta(x_1,x_2,x_3)=\det
\eta(x_1,x_2)\det \eta(x_1,x_3)\det \eta(x_2,x_3)\ ,\qquad\quad g=2\
,} \eqn\dettre{\det \eta\eta(x_1,\ldots,x_6)={1\over
15}\sum_{s\in\perm'_6}\sgn(s) \prod_{i=1}^4\det
\eta(x_{d^i_1(s)},x_{d^i_2(s)},x_{d^i_3(s)})\ , \qquad\quad g=3\ ,}
\eqn\detvan{\sum_{s\in\perm_M}\sgn(s) \prod_{i=1}^{g+1}\det
\eta(x_{d^i(s)})=0\ , \qquad\quad g\geq4\ ,} where $\{\eta_i\}_{i\in
I_g}$ is an arbitrary basis of $H^0(K_C)$ and $x_i$, $i\in I_M$,
are arbitrary points of $C$, hold. Furthermore, they are equivalent
to \eqn\due{\det
\eta\eta(x_1,x_2,x_3)=-\kuno[\eta]^3{\prod_{i=1}^3\deltadiv(\sum_{j=1}^3x_j-2x_i)\prod_{1}^3\sigma(x_j)\over
\prod_{i<j}E(x_i,x_j)}\ ,} \vskip 3pt \noindent for $g=2$ \vskip 3pt
\eqn\tre{\det
\eta\eta(x_1,\ldots,x_6)={\kuno[\eta]\over
15}^4\prod_{i=1}^6\sigma(x_i)^2\sum_{s\in\perm'_6}\sgn(s)\prod_{k=1}^4
{\deltadiv\bigl(\sum_{i=1}^3x_{d^k_i(s)}-y_{k,s}\bigr)\prod_{i<j}^3E(x_{d_i^k(s)},x_{d_j^k(s)})
\over \prod_{i=1}^3E(y_{k,s},x_{d^k_i(s)})\sigma(y_{k,s})}\ , }
\vskip 3pt \noindent for $g=3$
\vskip 3pt
\eqn\thetavan{{\kuno[\eta]\over
c_g}^{g+1}\prod_{l=1}^M\sigma(x_l)^2\sum_{s\in\perm_M}\sgn(s)\prod_{k=1}^{g+1}
{\deltadiv\bigl(\sum_{i=1}^gx_{d^k_i(s)}-y_{k,s})\prod_{i<j}^gE(x_{d_i^k(s)},x_{d_j^k(s)})
\over \prod_{i=1}^gE(y_{k,s},x_{d^k_i(s)})\sigma(y_{k,s})}=0\ ,}
\vskip 3pt
\noindent for
$g\geq4$, where $y_{k,s}$, $k\in I_{g+1}$, $s\in \perm_M$, are
arbitrary points of $C$.}

\vskip 6pt

\noindent {\bf Proof.} Eqs.\detdue-\detvan\ follow by applying Lemma
\thcombi\ to $\det\eta\eta(x_1,\ldots,x_M)$ and noting that it
vanishes for $g\ge 4$. Eqs.\due-\thetavan\ then follow by
Eq.\dettheta. \hfill$\square$

\vskip 6pt

In \ref\DHokerQP{
  E.~D'Hoker and D.~H.~Phong,
  Two-loop superstrings. IV: The cosmological constant and modular forms,
  {\it Nucl.\ Phys.\ } B {\bf 639} (2002), 129-181.
  }
D'Hoker and Phong made the interesting observation that for $g=2$
\eqn\dhph{\det \omega\omega(x_1,x_2,x_3)=\det \omega(x_1,x_2)\det
\omega(x_1,x_3)\det \omega(x_2,x_3)\ ,} that proved by first
expressing the holomorphic differentials in the explicit form and
then using the product form of the Vandermonde determinant.
Eq.\dhph\ corresponds to \detdue\ when the generic basis
$\eta_1,\eta_2$ of $H^0(K_C)$ is the canonical one. On the other
hand, the way \detdue\ has been derived shows that \dhph\ is an
algebraic identity since it does not need the explicit hyperelliptic
expression of $\omega_1$ and $\omega_2$. Eq.\dhph\ is the first case
of the general formulas, derived in Lemmas \thcombi\ and \thcombvi,
expressing the determinant of the matrix $ff_i(x_j)$ in terms of a
sum of permutations of products of determinants of the matrix
$f_i(x_j)$. In particular, by \dettre, for $g=3$ we have
$$\det
\omega\omega(x_1,\ldots,x_6)={1\over15}\sum_{s\in\perm_6'}\sgn(s)
\prod_{i=1}^4\det\omega(x_{d^i(s)})\ .$$

For $n<g$, a necessary condition for Eq.\combi\ to hold is the
existence of the points $p_i$, $3\le i\le g$, satisfying
Eq.\combcond; in particular, Lemmas \thcombi\ and \thcombvi\ can be
applied to the basis $\{\sigma_i\}_{i\in I_g}$, of
$H^0(K_C)$, defined in Eq.\newbasis.

\vskip 6pt

\newth\thdetvp{Theorem}{Fix the points $p_1,\ldots,p_g\in C$,
and $\hat\sigma_i\in H^0(K_C)$, $i\in I_g$, in such a way that
$\hat\sigma_i(p_j)=0$, for all $i\neq j\in I_g$. Define
$\hat v_i\in H^0(K_C^2)$, $i\in I_N$, by
$$\hat v_i:=\psi\spbase{\hat\sigma}_i=
\hat\sigma_{\1_i}\hat\sigma_{\2_i}\ ,$$
and let $\{\eta_i\}_{i\in I_g}$ be an
arbitrary basis of $H^0(K_C)$. Then, the following identity \eqn\dettv{\eqalign{& \det \hat
v(p_3,\ldots,p_g,x_1,\ldots,x_{2g-1})\Bigl({\det\eta_i(p_j)\over\hat\sigma_1(p_1)\hat\sigma_2(p_2)}\Bigr)^{g+1}
\prod_{i=3}^g\hat\sigma_i(p_i)^{-4}\cr\cr
&={(-)\over c_{g,2}}^{g+1}\sum_{s\in\perm_{2g-1}}\sgn(s)\det
\eta(x_{d^1(s)})\det\eta(x_{d^2(s)})\prod_{i=3}^{g+1} \det
\eta(x_{d^i_1(s)},x_{d^i_2(s)},p_3,\ldots,p_g)\ ,}} holds for all
$x_1,\ldots,x_{2g-1}\in C$, where, according to \glic, $c_{g,2}=g!(g-1)!(2g-1)$.
}

\vskip 6pt

\noindent {\bf Proof.} Assume that $p_1,\ldots,p_g$ satisfy
the hypotheses of Proposition \thnewbasis, so that
$\{\hat\sigma_i\}_{i\in I_g}$ is a basis of $H^0(K_C)$ and $\hat\sigma_i(p_i)\neq 0$, for all $i\in
I_g$. Since the points $p_1,\ldots,p_g$ satisfying such a condition are a dense set in $C^g$, it
suffices to prove Eq.\dettv\ in this case and then conclude by continuity arguments.
A relation analogous to \devv\ holds
$$v_i(p_j)=
\left\{\vcenter{\vbox{\halign{\strut\hskip 6pt $ # $ \hfil & \hskip
2cm$ # $ \hfil\cr \hat\sigma_i(p_i)^2\delta_{ij}\ , & i\in I_g \ ,\cr 0\ ,& g+1\le
i\le M\ ,\cr}}}\right.$$
$j\in I_g$, so that
$$\det \hat
v(p_3,\ldots,p_g,x_1,\ldots,x_{2g-1})=(-)^{g+1}\prod_{i=3}^g\hat\sigma_i(p_i)^2\det_{I_{M,2}}\hat\sigma\hat\sigma(x_1,\ldots,x_{2g-1})\ .$$
By
Lemma \thcombi\ for $n=2$, $\det_{I_{M,2}}\hat\sigma\hat\sigma(x_1,\ldots,x_{2g-1})$ is equal to the RHS of \combi\ divided by
$\prod_{i=3}^g\hat\sigma_i(p_i)^{g-1}$. Eq.\dettv\ then follows by the identity
$$\det\hat\sigma_i(z_j)={\det\eta_i(z_j)\over \det\eta_i(p_j)}\det\hat\sigma_i(p_j)={\det\eta_i(z_j)\over \det\eta_i(p_j)}\prod_{i=1}^g\hat\sigma_i(p_i)\ .$$
\hfill$\square$

\vskip 6pt

\newrem\thdetvxxdx{Remark} If $\det\eta_i(p_j)\neq 0$, then Theorem
\thdetvp\ holds for $\hat\sigma_i\equiv\sigma_i$, so that $\hat\sigma_i(p_i)=1$, $i\in I_g$, and $\hat
v_i\equiv v_i$, $i\in I_N$.

\vskip 6pt

\newth\thdetvx{Corollary}{Let $\b:=\sum_{i=3}^gp_i$ be a fixed divisor of $C$ and define $\hat v_i$, $i\in I_N$, as in Theorem
\thdetvp. Then for all $x_1,\ldots,x_N\in C$
\eqn\detvi{\eqalign{&\det \hat v(x_1,\ldots,x_N)=-{F\over
c_{g,2}}{\deltadiv\bigl(\sum_{1}^Nx_i\bigr)\prod_{i=2g}^N(\sigma(x_i)^3
\prod_{j=1}^{i-1}E(x_j,x_i))\over
\deltadiv\bigl(\sum_{1}^{2g-1}x_i+\b\bigr)\prod_{i=3}^{g}
\prod_{j=1}^{2g-1}E(p_i,x_j)}\prod_{i=1}^{2g-1}\sigma(x_i)^2\cr\cr
\cdot&\sum_{s\in\perm_{2g-1}}\sgn(s)S\bigl(\smsum_{i=1}^gx_{s_i}\bigr)
S\bigl(\smsum_{i=g}^{2g-1}x_{s_i}\bigr)
\prod_{^{i,j=1}_{i<j}}^gE(x_{s_i},x_{s_j})\prod_{^{i,j=g}_{i<j}}^{2g-1}
E(x_{s_i},x_{s_j}) \cr\cr \cdot
&\prod_{k=1}^{g-1}\Bigl(S(x_{s_k}+x_{s_{k+g}}+\b)
E(x_{s_k},x_{s_{k+g}})\prod_{i=3}^gE(x_{s_k},p_i)E(x_{s_{k+g}},p_i)\Bigr)\
,}} where $F\equiv F(p_1,\ldots,p_g)$ is
$$F:=\Bigl(
{\hat\sigma_1(p_1)\hat\sigma_2(p_2)\over
S(\a)\sigma(p_1)\sigma(p_2)E(p_1,p_2)}\Bigr)^{g+1}
\prod_{i=3}^g{\hat\sigma_i(p_i)^4\over \sigma(p_i)^5
(E(p_1,p_i)E(p_2,p_i))^{g+1}\prod_{j>i}^gE(p_i,p_j)^3} \ .$$}

\vskip 6pt

\noindent {\bf Proof.} Apply Eq.\dettv\ to
$$\det \hat
v(x_1,\ldots,x_N)={\det \rho(x_1,\ldots,x_N)\over\det
\rho(p_3,\ldots,p_g,x_1,\ldots,x_{2g-1})} \det \hat
v(p_3,\ldots,p_g,x_1,\ldots,x_{2g-1})\ ,$$
with $\{\rho_i\}_{i\in I_N}$ an arbitrary
basis of $H^0(K_C^2)$. Eq.\detvi\ then follows by Eqs.\dettheta\detthetaii. \hfill$\square$

\vskip 6pt

\subsec{Characterization of the $\B$ locus and the divisor of $K$}

Proposition \thnewbasis\ shows that
$\det\eta_i(p_j)\neq 0$, for an arbitrary basis $\{\eta_i\}_{i\in
I_g}$ of $H^0_C(K)$, is a necessary and sufficient condition on the points
$p_1,\ldots,p_g$ for the existence of a basis of holomorphic $1$-differentials
$\{\hat\sigma_i\}_{i\in I_g}$, such that $i\neq j\Rightarrow\sigma_i(p_j)=0$, $i,j\in I_g$. By
Eq.\dettheta\ and \essesimm\ it follows that the subset $\A\subset C^g$, for which such a condition is not satisfied,
corresponds to the set of solutions of the equation
$$S(\a)
\prod_{i<j}^gE(p_i,p_j)=0\ .$$

It is more difficult to characterize the locus $\B\subset C^g$,
whose elements are the $g$-tuples of points $p_1,\ldots,p_g$ which
do not satisfy the conditions of Proposition \thlev. The following
theorems show that such a locus can be characterized as the set of
solutions of the equation $H=0$ for a suitable function
$H(p_1,\ldots,p_g)$.

\vskip 6pt

\newth\thlevsing{Theorem}{Fix $g-2$ distinct points $p_3,\ldots,p_g\in
C$ such that \eqn\ipotesi{\{I(p+\b-\Delta)|p\in C\}\cap
\Theta_{s}=\emptyset\ ,} $\b:=\sum_3^gp_i$. Then, for each $p_2\in
C\setminus \{p_3,\ldots,p_g\}$, there exists a finite set of points
$S$, depending on $\b$ and $p_2$, with $\{p_2,\ldots,p_g\}\subset
S\subset C$, such that, for all $p_1\in C\setminus S$, the
holomorphic $1$-differentials $\{\sigma_i\}_{i\in I_g}$, associated
to the points $p_1,\ldots,p_g$ by Proposition \thnewbasis, is a
basis of $H^0(K_C)$ and the corresponding quadratic differentials
$\{v_i\}_{i\in I_N}$ is a basis of $H^0(K_C^2)$. Conversely, if for
some fixed $g-2$ arbitrary points $p_3,\ldots,p_g\in C$, there exist
$p_1,p_2\in C$ such that the associated $\{\sigma_i\}_{i\in I_g}$
and $\{v_i\}_{i\in I_N}$ are bases of $H^0(K_C)$ and $H^0(K_C^2)$,
then \ipotesi\ holds.}

\vskip 6pt

\noindent{\bf Proof.}
Eq.\ipotesi\ implies that $h^0(K_C\otimes\O(-\b-p))=1$, for all $p\in
C$.  Hence,
$h^0(K_C\otimes\O(-\b))= 2$ and, for each pair of linearly independent
elements $\sigma_1,\sigma_2$ of $H^0(K_C\otimes\O(-\b))$, the supports of $(\sigma_1)-\b$ and $(\sigma_2)-\b$ are disjoint.
Fix $p_2\in C\setminus\{p_3,\ldots,p_g\}$
and let $\sigma_1$ be a non-vanishing element of $H^0(K_C\otimes\O(-\b-p_2))$. Define the finite set $S$ as the support
of $(\sigma_1)$ or, equivalently, as the union of $\{p_2,\ldots,p_g\}$ and the set of zeros of
$S(x+p_2+\b)$. Then, for all $p_1\in C\setminus S$, fix $\sigma_2\in H^0(K_C\otimes\O(-\b-p_1))$ so
that $\sigma_1$ and $\sigma_2$ are linearly independent.
Then $p_1,\ldots,p_g$ satisfy the conditions {\it i}) and {\it ii}) of Proposition
\thlev, and $\{v_i\}_{i\in I_N}$, as defined in \lev, is a basis of
$H^0(K_C^2)$. Conversely, if $I(p+\b-\Delta)\in\Theta_{s}$ for some $p\in C$, then, for each pair $\sigma_1,\sigma_2\in
H^0(K_C\otimes\O(-\b))$, their greatest common divisor satisfies $\gcd(\sigma_1,\sigma_2)\ge p+\b$ and the condition {\it ii}) of Proposition
\thlev\ does not hold.\hfill$\square$

\vskip 6pt

The classical result that the dimension of $\Theta_s$ is $g-4$ for a non-hyperelliptic Riemann
surface of genus $g\ge 4$, immediately gives the following corollary by simple dimensional considerations.

\vskip 6pt

\newth\thdimthsing{Corollary}{In a non-hyperelliptic Riemann surface $C$ of genus $g\ge 4$, there
always exist $g$ points $p_1,\ldots,p_g\in C$ such that the corresponding $\{v_i\}_{i\in I_N}$ is a
basis of $H^0(K_C^2)$.}\

\vskip 6pt

\noindent{\bf Proof.} By Theorem \thlevsing, it is sufficient to
prove that there exists $\b\in C_{g-2}$ satisfying the condition
\ipotesi. Suppose, by absurd, that this is not true. Then, a
translation of $W_{g-2}=I(C_{g-2})$ is a subset of $\Theta_s\ominus
W_1:=\{e-I(p)\mid e\in\Theta_s, p\in C\}$. The corollary then
follows by observing that $W_{g-2}$ has dimension $g-2$, whereas the
dimension of each component of $\Theta_s\ominus W_1$ is less than
$\dim\Theta_s+\dim W_1=g-3$.\hfill$\square$

\vskip 6pt

\newth\thlaH{Theorem}{Fix $p_1\ldots,p_g\in C$. The function $H\equiv H(p_1,\ldots,p_g)$
\eqn\ilrappo{\eqalign{H&:= {S(\a)^{5g-7}E(p_1,p_2)^{g+1}\over
\deltadiv\bigl(\b+\sum_{1}^{2g-1}x_i\bigr)\prod_{i=1}^{2g-1}
\sigma(x_i)}
\prod_{i=3}^g{E(p_1,p_i)^4E(p_2,p_i)^4\prod_{j>i}^gE(p_i,p_j)^5\over
\sigma(p_i)} \cr\cr &\cdot \sum_{s\in\perm_{2g-1}}
{S\bigl(\sum_{i=1}^gx_{s_i}\bigr)
S\bigl(\sum_{i=g}^{2g-1}x_{s_i}\bigr)\over
\prod_{i=3}^gE(x_{s_g},p_i)}
\prod_{i=1}^{g-1}{S(x_{s_i}+x_{s_{i+g}}+\b)\over \prod_{_{j\neq
i}^{j=1}}^{\ss g-1}E(x_{s_i}, x_{s_{j+g}})} \ ,}} is independent of
the points $x_1,\ldots,x_{2g-1}\in C$. Furthermore, the set
$\{v_i\}_{i\in I_N}$, defined as in {\rm \lev}, is a basis of
$H^0(K_C^2)$ if and only if $H\neq 0$.}

\vskip 6pt

\noindent {\bf Proof.} Consider the holomorphic
$1$-differentials
$$\hat\sigma_i(z):=A_i^{-1}\sigma(z)S(\a_i+z)
\prod_{^{j=1}_{j\neq
i}}^gE(z,p_j)=A_i^{-1}\sum_{j=1}^g\theta_{\Delta,j}(\a_i)\omega_j(z)\
,$$ $i\in I_g$, with $\a_i$ as in Definition \defAB\ and
$A_1,\ldots,A_g$ non-vanishing constants. If the points
$p_1,\ldots,p_g$ satisfy the hypotheses of Proposition \thnewbasis,
then $\{\hat\sigma_i\}_{i\in I_g}$ corresponds, up to a non-singular
diagonal transformation, to the basis defined in
\newbasis. Let $\{\rho_i\}_{i\in I_N}$ be an arbitrary
basis of $H^0(K_C^2)$. By \detthetaii\ the following identity
$$\eqalignno{&\det\rho(p_3,\ldots,p_g,x_1,\ldots,x_{2g-1})\cr\cr & =
\kdue[\rho]
\sgn(s)\prod_{_{i<j}^{i,j=1}}^{2g-1}E(x_{s_i},x_{s_j})
\deltadiv\Bigl(\smsum_{\ss 1}^{\ss
2g-1}x_i+\b\Bigr)\prod_{i=1}^{2g-1}\sigma(x_i)^3\prod_{i=3}^g\sigma(p_i)^3
\prod_{_{i<j}^{i,j=3}}^gE(p_i,p_j)\prod_{i=3}^g\prod_{j=1}^{2g-1}E(p_i,x_j)
\ ,}
$$
holds for all $s\in\perm_{2g-1}$. Together with Eq.\detvi\ and the
above expression for $\hat\sigma_i$, it implies that
\eqn\ratioo{H=\kappa[\rho]c_{g,2} (A_1A_2)^{g+1}\prod_{i=3}^gA_i^4
{\det \hat v(p_3,\ldots,p_g,x_1,\ldots,x_{2g-1})\over
\det\rho(p_3,\ldots,p_g,x_1,\ldots,x_{2g-1})}\ .} Hence, $H$ is
independent of $x_1,\ldots,x_{2g-1}$, and $H\neq0$ if and only if
$\{\hat v_i\}_{i\in I_N}$ is a basis of $H^0(K_C^2)$. On the other
hand the vector $(\hat v_1,\ldots,\hat v_N)$ corresponds, up to a
non-singular diagonal transformation, to $(v_1,\ldots,v_N)$, with
$v_i$, $i\in I_N$, defined in \lev. \hfill$\square$

\vskip 6pt

\newrem\gfgfg{Remark}{By \ratioo\
$$\kappa[\hat v]={H(p_1,\ldots,p_g)\over c_{g,2}
(A_1A_2)^{g+1}\prod_{i=3}^gA_i^4}\ .$$
Furthermore, if $(p_1,\ldots,p_g)\notin \A$, then one can choose
$$A_i=\sigma(p_i)S(\a)\prod_{^{j=1}_{j\neq i}}^gE(p_i,p_j)=\sum_{j=1}^g\theta_{\Delta,j}(\a_{i})\omega_j(p_i)\ ,$$
to obtain $\hat\sigma_i\equiv\sigma_i$, $i\in I_g$, and
\eqn\kkl{\eqalign{\kappa[v]&={H(p_1,\ldots,p_g)\over c_{g,2}
\prod_{i=1}^2\bigl(\sum_{j=1}^g\theta_{\Delta,j}(\a_i)\omega_j(p_i)
\bigr)^{g+1}
\prod_{i=3}^g\bigl(\sum_{j=1}^g\theta_{\Delta,j}(\a_i)\omega_j(p_i)\bigr)^4}\cr\cr
&={H(p_1,\ldots,p_g)\over
c_{g,2}S(\a)^{6g-6}\prod_{i=1}^2\bigl(\sigma(p_i)\prod_{^{j=1}_{j\neq
i}}^gE(p_i,p_j)\bigr)^{g+1}
\prod_{i=3}^g\bigl(\sigma(p_i)\prod_{^{j=1}_{j\neq
i}}^gE(p_i,p_j)\bigr)^4}\ .}} Observe that $\A\subset \B$. Theorem
\thlevsing\ shows that if $(p_1,\ldots,p_g)\notin\A$, a necessary
and sufficient condition for $(p_1,\ldots,p_g)$ to be in $\B$ is
that there exists $p\in C$ such that $I(\b+p-\Delta)\in\Theta_{s}$.
Hence, $\B$ is the union of $\A$ together with the pull-back of a
divisor in $C^{g-2}$ by the projection $C^{g}\to C^{g-2}$ which
``forgets'' the first pair of points:
$(p_1,\ldots,p_g)\to(p_3,\ldots,p_g)$. Such a divisor is
characterized by the equation $K=0$, where $K$ is defined in the
following corollary.}

\vskip 6pt

\newth\corollo{Corollary}{Define
\eqn\iltrap{\eqalign{K(p_3,\ldots,p_g)&:={1\over
\deltadiv\bigl(\b+\sum_{1}^{2g-1}x_i\bigr)\prod_{1}^{2g-1}\sigma(x_i)\prod_{i=3}^g\sigma(p_i)
}\cr\cr &\cdot \sum_{s\in\perm_{2g-1}}
{S\bigl(\sum_{i=1}^gx_{s_i}\bigr)
S\bigl(\sum_{i=g}^{2g-1}x_{s_i}\bigr) \over
\prod_{i=3}^gE(x_{s_g},p_i)}\prod_{i=1}^{g-1}{
S(x_{s_i}+x_{s_{i+g}}+\b) \over
\prod_{_{j\neq
i}^{j=1}}^{\ss g-1}E(x_{s_i}, x_{s_{j+g}})} \ .}}
\smallskip
\item{\it a.} $K\equiv K(p_3,\ldots,p_g)$ is independent of
$x_1,\ldots,x_{2g-1}\in C$.
\smallskip
\item{\it b.} For any
$p_1,\ldots,p_g\in C$ such that $\det\eta_i(p_j)\neq 0$, the set
$\{v_i\}_{i\in I_N}$, defined in {\rm \lev}, is a basis of
$H^0(K_C^2)$ if and only if $K\neq 0$.
\smallskip
\item{\it c.}
\eqn\baaabba{S(p_1+p_2+\b)=0\ ,\; \forall p_1,p_2\in C \qquad \Longrightarrow \qquad K=0 \ .}
\smallskip
\item{\it d.} If $p_3,\ldots,p_g$ are pairwise distinct and $K\neq
0$, then there exist $p_1,p_2\in C$ such that $H\neq 0$.
\smallskip}

\vskip 6pt

\noindent {\bf Proof.}
\item{--} {\it a.} The ratio
\eqn\baaa{{H\over K}=S(\a)^{5g-7}E(p_1,p_2)^{g+1}\prod_{i=3}^g(E(p_1,p_i)E(p_2,p_i))^4
\prod_{_{i<j}^{i,j=3}}^gE(p_i,p_j)^5\ ,}
is independent of
$x_1,\ldots,x_{2g-1}$, so that $a$) follows by Theorem \thlaH\ or, equivalently, noticing that
by Eqs.\kkkk\ilrappo\kkl\ and \baaa\
\eqn\Kkappa{K(p_3,\ldots,p_g):=(-)^{g+1}c_{g,2}{\kdue[v]\over\kuno[\sigma]^{g+1}}\prod_{_{i<j}^{i,j=3}}^gE(p_i,p_j)^{2-g}
\prod_{i=3}^g\sigma(p_i)^{3-g}\ .}
\smallskip
\item{--} {\it b.} By \dettheta\ and \baaa\ the condition
$\det\eta_i(p_j)\neq0$ implies $H/K\neq 0$. In this case $K\neq 0$
if and only if $H\neq 0$, and $b$) follows by Theorem \thlaH.
\smallskip
\item{--} {\it c.} If $S(p_1+p_2+\b)=0$, for all $p_1,p_2\in C$, then the numerators in each term of the sum in
\iltrap\ vanish for all $x_1,\ldots,x_{2g-1}\in C$. Since $K$ is
independent of $x_1,\ldots,x_{2g-1}$, it follows that the proof of
point {\it c}) is equivalent to prove that there exist
$x_1,\ldots,x_{2g-1}\in C$ such that the denominators in \iltrap\ do
not vanish. On the other hand, the possible zeros of such
denominators are the ones corresponding of the zeros of the primes
forms, which are avoided by simply choosing
$p_3,\ldots,p_g,x_1,\ldots,x_{2g-1}$ pairwise distinct, and the ones
of $\theta_\Delta(\b+\sum_1^{2g-1}x_i)$. Fix an arbitrary $y\in C$
and set $w:=I(\b+\sum_{g+1}^{2g-1}x_i+y-2\Delta)$. Then
$$\theta\bigl(\b+\smsum_1^{2g-1}
x_i-3\Delta\bigr)=\theta\bigl(w+\smsum_{1}^{g}x_i-y-\Delta\bigr)\
,$$ and, by the Jacobi Inversion Theorem, by varying the points
$x_1,\ldots,x_g\in C$ one can span the whole Jacobian variety. Then,
one can always choose $x_1,\ldots,x_{2g-1}$ pairwise distinct and
distinct from $p_3,\ldots,p_g$ in such a way that
$\theta\bigl(w+\sum_{1}^{g}x_i-y-\Delta\bigr)\neq 0$, so that the
denominator does not vanish and {\it c}) follows.
\smallskip
\item{--} {\it d.} Since $K\neq 0$, by {\it c}) there exist $p_1,p_2\in C$ such that $S(p_1+p_2+\b)\neq
0$. By continuity arguments, it follows that there exist some
neighbourhoods $U_i\subset C$ of $p_i$, $i=1,2$, such that
$S(x_1+x_2+\b)\neq 0$ for all $(x_1,x_2)\in U_1\times U_2$. Hence,
we can choose $p_1,p_2$ so that $S(p_1+p_2+\b)\neq 0$ and
$p_1,\ldots,p_g$ are pairwise distinct. Then, by Eq.\baaa, $H/K\neq
0$ and, since $K\neq 0$, we conclude that $H\neq 0$. \hfill$\square$

\vskip 6pt

\noindent In view of Eq.\Kkappa, it is useful to define
\eqn\kpiccola{k(p_3,\ldots,p_g):=K(p_3,\ldots,p_g)\prod_{_{i<j}^{i,j=3}}^gE(p_i,p_j)^{g-2}
\prod_{i=3}^g\sigma(p_i)^{g-3}=(-)^{g+1}c_{g,2}{\kdue[v]\over\kuno[\sigma]^{g+1}}\
,} which is a holomorphic $(g-3)$-differential in each of its $g-2$
arguments.

\newth\quadgen{Theorem}{Fix $p_1,\ldots,p_g\in C$, with $C$ non-hyperelliptic of genus $g\ge 4$ and
let $\{\hat\sigma_i\}_{i\in I_g}$ be a set of non-vanishing holomorphic $1$-differentials such that
$i\neq j\,\Rightarrow\, \hat\sigma_i(p_j)= 0$, for all $i,j\in I_g$. The
following statements are equivalent
\smallskip
\item{\it i.} The conditions
\itemitem{\it $i'$.} $(p_1,\ldots,p_g)\notin \A$;
\itemitem{\it $i''$.} $\b:=\sum_{i=3}^gp_i$ is the greatest common divisor of
$(\sigma_1)$ and $(\sigma_2)$;
\smallskip
are satisfied;
\smallskip
\item{\it ii.} $H(p_1,\ldots,p_g)\neq 0$, where $H$ is defined in Eq.\futuref{ilrappo};
\smallskip
\item{\it iii.} $\{\hat v_i\}_{i\in I_N}$ is a basis of $H^0(K_C^2)$, with $\hat
v_i:=\hat\sigma\hat\sigma_i$, $i\in I_M$.
\smallskip
\noindent More generally, fix $p_3,\ldots,p_g\in C$. The
following statements are equivalent:
\smallskip
\item{\it iv.} $p_3,\ldots,p_g$ are pairwise distinct and $\{I(p+\b-\Delta)|p\in C\}\cap\Theta_{s}=\emptyset$;
\smallskip
\item{\it v.} $p_3,\ldots,p_g$ are pairwise distinct and $K(p_3,\ldots,p_g)\neq 0$, where $K$ is defined in Eq.\futuref{iltrap};
\smallskip
\item{\it vi.} There exist $p_1,p_2\in C$ such that $p_1,\ldots,
p_g$ satisfy {\it i}), {\it ii}) and {\it
iii});
\smallskip
\item{\it vii.} For all $p\in C$, $S(x+p+\b)$ does not vanish
identically as a function of $x$; furthermore, for each $p_2\in C\setminus\{p_3,\ldots,p_g\}$, the
points $p_1,\ldots,
p_g$ satisfy {\it i}), {\it ii}) and {\it
iii}) if and only if $p_1$ is distinct from $p_2,\ldots,p_g$ and from the $g-1$ zeros of
$S(x+p_2+\b)$.
\smallskip
}

\vskip 6pt

\noindent{\bf Proof.}

\item{--} {\it i}) $\Leftrightarrow$ {\it iii}) is proved in Proposition \futuref{thlev} (in the direction {\it i})
$\Rightarrow$ {\it iii}), only the case of normalized
$1$-differentials $\sigma_i(p_i)=1$, for all $i\in I_g$ is
considered; however, by the hypothesis {\it $i'$}), the general case
can be reduced to this choice by a non-singular diagonal
transformation on $\{\hat\sigma_i\}_{i\in I_g}$);

\item{--} {\it ii}) $\Leftrightarrow$ {\it iii}) is proved in Theorem \futuref{thlaH};

\item{--} {\it vii}) $\Rightarrow$ {\it vi}) is obvious;

\item{--} {\it iv}) $\Leftrightarrow$
{\it vii}) follows by first noting that $S(x+p+\b)$ identically
vanishes as a function of $x$ if and only if
$I(p+\b-\Delta)\in\Theta_{s}$, and then by Theorem \thlevsing; in
particular, in such a theorem it is proved that for each fixed
$p_2\in C\setminus\{p_3,\ldots,p_g\}$, the points $p_1,\ldots,p_g$
satisfy {\it i}) if and only if the conditions
$p_1\notin\{p_2,\ldots,p_g\}$, $S(p_1+p_2+\b)\neq 0$ and {\it iv})
hold;

\item{--} {\it vi})
$\Rightarrow$ {\it iv}) also follows by Theorem \thlevsing, where it is proved that if {\it iv}) does not
hold, then {\it $i''$}) cannot be satisfied;

\item{--} {\it v}) $\Leftrightarrow$ {\it vi}), finally, follows by Corollary
\futuref{corollo}, where it is proved that {\it $i'$}) and {\it v}) are equivalent to {\it
ii}) and that if {\it v}) holds, then there
exist $p_1,p_2\in C$ such that $p_1,\ldots,p_g$ satisfy {\it ii}). \hfill$\square$

\vskip 6pt

The function $K(p_3,\ldots,p_g)$ defined in Eq.\futuref{iltrap},
whose zero divisor is characterized in the theorem above, is the
fundamental tool in the proof of the following theorem. Such a
result heavily relies on the properties of $\Theta_s$ in the case
the sublying ppav is the Jacobian torus of a canonical curve. By the
Riemann Singularity Theorem,
$$\Theta_s=W_{g-1}^1+\K^{p_0}\equiv I(C_{g-1}^1)-I(\Delta)\ ,$$
where $p_0\in C$ is the base point for $I$, $W_{g-1}^1:=
I(C_{g-1}^1)$ and $C^1_{g-1}\subset C_{g-1}$ is the subvariety of
codimension $2$ in $C_{g-1}$, whose elements are the special
effective divisors of degree $g-1$. Note that each effective divisor
$\d\in C_{g-3}$ of degree $g-3$ canonically determines an embedding
$\pi_d:C_2\hookrightarrow C_{g-1}$, $C_2\ni\c\mapsto \c+\d\in
C_{g-1}$ of $C_2$ as a subvariety of dimension $2$ in $C_{g-1}$.
Hence, by a simple dimensional counting, we expect the intersection
$C_{g-1}^1\cap \pi_d(C_2)$ to have (in general) dimension $0$. The
following theorem shows that, in the general case in which such an
intersection does not contain any component of dimension greater
than $0$, $C_{g-1}^1\cap \pi_d(C_2)$ corresponds (set-theoretically)
to a set of $g(g-3)/2$ points; furthermore, a remarkable relation of
such a set of points with the canonical divisor is given. Since the
restriction of the Abel-Jacobi map to $C_2$ is an injection (because
$C$ is non-hyperelliptic), such points are in one to one
correspondence with the points in the intersection $\Theta_s\cap
(W_2+I(\d-\Delta))$, where $W_2=I(C_2)$. \vskip 6pt

\newth\zeriK{Theorem}{Let $C$ be non-hyperelliptic of genus $g\ge 4$ and fix $p_4,\ldots,p_g\in C$.
Then, either:
\smallskip
\item{a.} For each point $p\in C$, there exists a point $q\in C$ such that
$$I(p+q+p_4+\ldots+p_g-\Delta)\in\Theta_{s}\ ;$$
\vskip 0pt
\noindent or:
\item{b.} There exist $k:=g(g-3)/2$ effective
divisors
$\c_1,\ldots,\c_k\in C_2$ of degree $2$, such that
\eqn\eeeii{e_i:=I(\c_i+p_4+\ldots+p_g-\Delta)\in\Theta_{s}\ ,\qquad \forall i\in I_k\ .}
Moreover, $\sum_{i=1}^k\c_i+(g-2)\sum_{i=4}^gp_i$ is the divisor of a holomorphic
$(g-3)$-differential on $C$.
\smallskip}

\vskip 6pt

\noindent {\bf Proof.} Consider $K(z,p_4,\ldots,p_g)$ as a function
of $z$. It vanishes at $z\equiv p$ if and only if there exists a
point $q\in C$ such that
$I(p+q+p_4+\ldots+p_g-\Delta)\in\Theta_{s}$. Then, $K= 0$ for all
$z\in C$ if and only if statement $a)$ holds.

\noindent Now, assume that $K(z,p_4,\ldots,p_g)$ is not identically
vanishing and consider
\eqn\kdiff{\phi(z):=K(z,p_4,\ldots,p_g)\prod_{i=4}^g
E(z,p_i)^{g-2}\sigma(z)^{g-3}\ .} By \futuref{Kkappa}, $\phi$ is a
holomorphic $(g-3)$-differential on $C$. Therefore, the divisor $\d$
of $K(z,p_4,\ldots,p_g)$ is effective ($K$ has no poles) of degree
$g(g-3)$ and $\d+(g-2)\sum_{i=4}^gp_i$ is the divisor of a
$(g-3)$-differential. It only remains to prove that $\d$ is the sum
of all the effective divisors of degree $2$ satisfying Eq.\eeeii. By
the equivalence of {\it iv}) and {\it v}) in Theorem
\futuref{quadgen}, if $\c:=q_1+q_2$ satisfies Eq.\eeeii, then $q_1$
and $q_2$ are both zeros of $K$. By construction,
$K(z,p_4,\ldots,p_g)$ can be written as
$$K(z,p_4,\ldots,p_g)=F(z,p_4,\ldots,p_g,x_1,\ldots,x_{2g-1})\det\varphi_i(x_j)\ ,$$ where
$\varphi_1,\ldots,\varphi_{2g-1}$ is a set of generators (depending
on $z,p_4,\ldots,p_g$) of $H^0(K_C^2\otimes\O(-z-p_4-\ldots-p_g))$
and $x_1,\ldots,x_{2g-1}$ are arbitrary points in $C$; $F$ is such
that, by Corollary \corollo, $K$ do not depend on
$x_1,\ldots,x_{2g-1}$. It is easy to verify that $K$ vanishes only
if $\det\varphi_i(x_j)=0$ for all $x_1,\ldots,x_{2g-1}\in C$; the
multiplicity of such a zero is $2g-1-r$, where
$r:=h^0(K_C^2\otimes\O(-z-p_4-\ldots-p_g))$. The space
$H^0(K_C^2\otimes\O(-z-p_4-\ldots-p_g))$ is generated by elements
$\sigma_1\eta$, $\sigma_2\rho$, as $\eta,\rho$ vary in $H^0(K_C)$;
here, $\sigma_1,\sigma_2$ is a basis for the $2$-dimensional space
$H^0(K_C\otimes\O(-z-p_4-\ldots-p_g))$ (note that if there exists
$q\in C$ such that $q+p_4+\ldots+p_g$ is special, then
$K(z,p_4,\ldots,p_g)$ identically vanishes). Proposition \thlev\
shows that $K(z,p_4,\ldots,p_g)\neq 0$, that is $r=2g-1$, if and
only if $h^0(K_C\otimes\O(-q-z-p_4-\ldots-p_g))=1$ (or,
equivalently, $q+z+p_4+\ldots+p_g$ is not special) for all $q\in C$.
Let $q_1$ be a zero of $K$ and denote by $n$ the maximal integer for
which there exist $n-1$ points $q_2,\ldots,q_n\in C$ such that
$h^0(K_C\otimes\O(-q_1-\ldots-q_n-p_4-\ldots-p_g)=2$. By the
considerations above, since $q_1$ is a zero, $n\ge 2$; furthermore,
$q_2,\ldots,q_n$ are zeros of $K$ too. Corollary \rango\ shows that
$$r\equiv h^0(K_C^2\otimes\O(-q_1-p_4-\ldots-p_g))=h^0(K_C^2\otimes\O(-q_1-\ldots-q_n-p_4-\ldots-p_g))=2g-n\ ,$$
so that the multiplicity of each $q_i$, $i\in I_n$, is $2g-1-r=n-1$.
Now, consider a zero $q'_1$ of $K(z,p_4,\ldots,p_g)$, distinct from
$q_1,\ldots,q_n$; by the same construction, if $q_1'$ has
multiplicity $n'-1$, with $n'\ge 2$, then it is an element of a set
of $n'$ (possibly coincident) zeroes $\{q_1',\ldots,q_n'\}$ with the
same multiplicity. By repeating this procedure, we obtain a finite
number $l$ of disjoint sets of zeroes; for each $i\in I_l$, the
$i$-th set contains $n_i\ge 2$ zeroes, we denote by
$q_1^i,\ldots,q_{n_i}^i$, each one with multiplicity $n_i-1$.
 Therefore, we have
$$\d=\sum_{i=1}^l\sum_{j=1}^{n_i}(n_i-1)q^i_j=\sum_{i=1}^l\sum_{j<k}^{n_i}(q^i_j+q^i_k)\ ,$$
and, since $h^0(K_C\otimes\O(-q^i_j-q^i_k-p_4-\ldots-p_g))=2$, each
$\c:=q^i_j+q^i_k$ satisfies Eq.\eeeii; conversely, it follows
immediately that if an element of $C_2$ satisfies Eq.\eeeii, then it
is the sum of a pair of zeroes of $K(z,p_4,\ldots,p_g)$ in the same
set.\hfill$\square$

\vskip 6pt

\newsec{Determinantal relations and combinatorial
$\theta$-identities}\seclab\deterrel

Denote by $\tilde\phi^n:H^0(K_C^n)\rightarrow \CC^{N_n}$ the
isomorphism $\tilde\phi^n(\phi^n_i)=e_i$, with $\{e_i\}_{i\in
I_{N_n}}$ the canonical basis of $\CC^{N_n}$. The isomorphism
$\tilde\eta$ induces an isomorphism
$\spisom{\eta}:\Sym^2(H^0(K_C^2))\rightarrow {\rm Sym}^2\CC^g$. The
natural map $\psi:\Sym^2(H^0(K_C^2))\rightarrow H^0(K_C^2)$ is
surjective if $C$ is canonical.

The choice of a basis $\{\eta_i\}_{i\in I_g}$ of $H^0(K_C)$
determines an embedding of the curve $C$ in $\PP_{g-1}$ by $p\mapsto
(\eta_1(p),\ldots,\eta_g(p))$, so that the elements of $\{\eta_i\}_{i\in I_g}$ correspond to a set
of homogeneous coordinates $X_1,\ldots,X_g$ on $\PP_{g-1}$. Each holomorphic $n$-differential
corresponds to a homogeneous $n$-degree polynomial in $\PP_{g-1}$ by
$$
\phi^n:=\sum_{i_1,\ldots,i_n}B_{i_1,\ldots,i_n}\eta_{i_1}\cdots\eta_{i_n}\mapsto
\sum_{i_1,\ldots,i_n}B_{i_1,\ldots,i_n}X_{i_1}\cdots X_{i_n}\ ,$$
where $X_1,\ldots,X_g$ are homogeneous coordinates on $\PP_{g-1}$. A
basis of $H^0(K_C^n)$ corresponds to a basis of the homogeneous
polynomials of degree $n$ in $\PP_{g-1}$ that are not zero when
restricted to $C$. The curve $C$ is identified with the ideal of all
the polynomials in $\PP_{g-1}$ vanishing at $C$. Enriques-Babbage
and Petri's Theorems state that, with few exceptions, such an ideal
is generated by quadrics
$$\sum_{j=1}^MC^\eta_{ij}XX_j=0\ ,$$
$N+1\le i\le M$,
where $XX_j:=X_{\1_j}X_{\2_j}$. Here, $\{C^\eta_i\}_{N<i\le M}$, with
$C^\eta_i:=(C^\eta_{i1},\ldots,C^\eta_{iM})$, is a set of linearly independent
elements of $\PP(\Sym^2\CC^g)\cong\PP_M$, each one defining a quadric. The isomorphism $\spisom{\eta}$
induces the identification
$\PP(\Sym^2(H^0(K_C)))\cong\PP_M$, under which each quadric corresponds to an element of $\ker\psi$
or, equivalently, to a relation among
holomorphic quadratic differentials
$$\sum_{j=1}^MC^\eta_{ij}\eta\eta_j=0\ .$$

Canonical curves that are not
cut out by such quadrics are trigonal or isomorphic to smooth plane
quintic. In these cases, Petri's Theorem assures that the ideal is
generated by the quadrics above together with a suitable set of
cubics.

This section is devoted to the study of such relations among quadratic and cubic differentials.

\vskip 6pt

\subsec{Relations among holomorphic quadratic differentials}

In the following we derive the matrix form of the map $\tilde
v\circ\psi\circ(\spisom{\sigma})^{-1}$, with respect to the basis
$\{\sigma_i\}_{i\in I_g}$ constructed in the previous subsection.
This will lead to the explicit expression of $\ker\psi$. Set
\eqn\dij{ \tilde\psi_{ij}:=
{\kdue[v_1,\ldots,v_{i-1},v_j,v_{i+1},\ldots,v_N]\over\kdue[v]} \ .}
$i\in I_N$, $j\in I_M$.

\vskip 6pt

\newth\thlemma{Lemma}{$v_1,\ldots,v_M$ satisfy the following $(g-2)(g-3)/2$ linearly independent relations
\eqn\lemma{v_i=\sum_{j=1}^N\tilde\psi_{ji}v_j=\sum_{j=g+1}^N\tilde\psi_{ji}v_j\ ,}}
$i=N+1,\ldots,M$.

\vskip 6pt

\noindent {\bf Proof.} The first equality trivially follows by the
Cramer rule. The identities \devv\
imply $\tilde\psi_{ji}=0$ for $j\in I_g$ and $i=N+1,\ldots,M$, and
the lemma follows. \hfill$\square$

\vskip 6pt

Eq.\lemma\ implies that the diagram
$$\matrix{\hfill\Sym^2(H^0(K_C^2))\, & \sopra{\longrightarrow}{\psi} & H^0(K_C^2) \cr\cr
 \spisom{\sigma}\downarrow\;\hfill & & \downarrow \tilde v
\cr\cr \hfill\CC^M\! \qquad\quad & \sopra{\longrightarrow}{\tilde\psi} &
\CC^N}$$ where $\tilde\psi:\CC^M\rightarrow \CC^N$ is the
homomorphism with matrix elements $\tilde\psi_{ij}$ and ${\rm
Sym}^2\CC^g$ is isomorphic to $\CC^M$ through $A$, introduced in
Definition \definizione, commutes.

Let $\iota:\CC^N\rightarrow \CC^M$ be the injection
$\iota(e_i)=\tilde{e}_i$, $i\in I_N$. The matrix elements of the map
$\iota\circ\tilde\psi:\CC^M\rightarrow \CC^M$ are
$$(\iota\circ\tilde\psi)_{ij}=
\left\{\vcenter{\vbox{\halign{\strut\hskip 6pt $ # $ \hfil & \hskip
2cm$ # $ \hfil\cr \tilde\psi_{ij}\ , &1\le i\le N\ ,\cr 0\ ,&N+1\le
i\le M\ ,\cr}}}\right.$$ $j\in I_M$. Noting that
$(\iota\circ\tilde\psi)_{ij}=\delta_{ij}$, for all $i,j\in I_N$, we
obtain
$$\sum_{i=1}^M(\iota\circ\tilde\psi)_{ji}(\iota\circ\tilde\psi)_{ik}=
\sum_{i=1}^N(\iota\circ\tilde\psi)_{ji}\tilde\psi_{ik}=(\iota\circ\tilde\psi)_{jk}\
,$$ $j,k\in I_M$. Hence, $\iota\circ\tilde\psi$ is a
projection of rank $N$ and, since $\iota$ is an injection,
\eqn\ilker{\ker\tilde\psi=\ker \iota\circ\tilde\psi= ({\rm
id}-\iota\circ\tilde\psi)(\CC^M)\ .}

\vskip 6pt

\newth\thutilde{Lemma}{The set $\{\tilde u_{N+1},\ldots,\tilde u_M\}$, $\tilde u_i:=\tilde e_i-\sum_{j=1}^N\tilde
e_j\tilde\psi_{ji}$, $N+1\leq i\le M$, is a basis of $\ker\tilde\psi$.}

\vskip 6pt

\noindent{\bf Proof.} Since $({\rm
id}-\iota\circ\tilde\psi)(\tilde{e_i})=0$, $i\in I_N$, by \ilker,
the $M-N$ vectors $\tilde u_i=({\rm
id}-\iota\circ\tilde\psi)(\tilde{e_i})$, $N<i\le M$, are a set of
generators for $\ker\psi$ and, since ${\rm dim}\ker\psi =M-N$, the
lemma follows. \hfill$\square$

\vskip 6pt

\noindent Set $\eta\eta_i:=\psi\spbase{\eta}_i$, $i\in I_g$, and let
$X^\eta$ be the automorphism on $\CC^M$ in the commutative diagram
$$\matrix{\hfill\Sym^2(H^0(K_C^2))\, & \sopra{\longrightarrow}{\rm id} & \,\Sym^2(H^0(K_C^2))\hfill \cr\cr
 \spisom{\sigma}\downarrow\;\hfill & & \hfill\;\downarrow \spisom{\eta}
\cr\cr \hfill\CC^M\qquad\quad & \sopra{\longrightarrow}{X^\eta} &\qquad\qquad
\CC^M\hfill }$$ whose matrix elements are
\eqn\XXXX{X^\eta_{ji}=\chi_j^{-1}([\eta]^{-1}[\eta]^{-1})_{ij}={[\eta]^{-1}_{\1_i\1_j}[\eta]^{-1}_{\2_i\2_j}+
[\eta]^{-1}_{\1_i\2_j}[\eta]^{-1}_{\2_i\1_j}\over
1+\delta_{\1_j\2_j}}\ ,}
$i,j\in I_M$, so that
\eqn\vXeta{v_i=\sum_{j=1}^MX^\eta_{ji}\,\eta\eta_j\ ,} $i\in I_M$. Since $\eta\eta_i$, $i\in I_M$, are linearly dependent,
the matrix $X^\eta_{ij}$ is not univocally determined by Eq.\vXeta.
More precisely, an endomorphism $X^{\eta'}\in{\rm End}(\CC^M)$
satisfies Eq.\vXeta\ if and only if the diagram
$$\matrix{\hfill\CC^M\, & \sopra{\longrightarrow}{\tilde\psi} & \,\CC^N\hfill \cr\cr
 X^{\eta'}\downarrow\;\hfill & & \hfill\;\downarrow {\rm id}
\cr\cr \hfill\CC^M\! & \sopra{\longrightarrow}{B^\eta} &\!
\CC^N\hfill }$$ where $B^\eta:=\tilde\psi\circ (X^\eta)^{-1}$,
commutes or, equivalently, if and only if
\eqn\lacond{(X^{\eta'}-X^\eta)(\CC^M)\subseteq
X^\eta(\ker\tilde\psi)\ .} Next theorem provides an explicit
expression for such a homomorphisms. Consider the following
determinants of the $d$-dimensional submatrices of $X^\eta$
$$
\Xmin{\eta}{\ss j_1 \ldots j_d\cr\ss i_1 \ldots i_d}:=
\det\left(\matrix{X^\eta_{i_1j_1} & \ldots & X^\eta_{i_1j_d}\cr
\vdots & \ddots & \vdots \cr X^\eta_{i_dj_1} & \ldots &
X^\eta_{i_dj_d}}\right)\ , $$
$i_1,\ldots,i_d,j_1,\ldots,j_d\in
I_M$, $d\in I_M$.

\vskip 6pt

\newth\thcorol{Theorem}
{\eqn\corol{\sum_{j=1}^MC^\eta_{ij}\eta\eta_j=0\ ,} $N+1\leq i\le
M$, where \eqn\leC{ C^\eta_{ij}:=\sum_{k_1,\ldots,k_N=1}^M
\Xmin{\eta}{\ss 1 \hfill \ldots \hfill N i \cr \ss k_1\hfill \ldots
\hfill k_N j\hfill}\;
{\kdue[\eta\eta_{k_1},\ldots,\eta\eta_{k_N}]\over\kdue[v]} \ ,} are
$M-N$ independent linear relations among holomorphic quadratic
differentials. Furthermore, for all $p\in C$
\eqn\detv{W[v](p)=\sum_{i_1,\ldots,i_N=1}^M\Xmin{\eta}{\ss 1 \hfill
\ldots \hfill N \cr \ss i_1\hfill \ldots \hfill
i_N}\;W[\eta\eta_{k_1},\ldots,\eta\eta_{k_N}](p)\ .}}

\vskip 6pt

\noindent {\bf Proof.} By \lemma\ and \vXeta\
$$\sum_{j=1}^M(X^\eta_{ji}-\sum_{k=1}^N\tilde\psi_{ki}X^\eta_{jk})\eta\eta_j=0\ ,$$
for all $N+1\le i\le M$, and by \dij\
$$\sum_{j=1}^M\biggl[\sum_{k=1}^N(-)^k
{\kdue[v_i,v_1,\ldots,\check{v}_k,\ldots,
v_N]\over\kdue[v]}X^\eta_{jk}+X^\eta_{ji}\biggr]\eta\eta_j=0\ .$$ By
\vXeta\
$${\kdue[v_{i_1},\ldots, v_{i_N}]\over\kdue[v]}=\sum_{k_1,\ldots,k_N=1}^M
\Xmin{\eta}{\ss i_1\ldots i_N\cr \ss k_1\ldots
k_N}{\kdue[\eta\eta_{k_1},\ldots,\eta\eta_{k_N}]\over\kdue[v]}\ ,$$
$i_1,\ldots,i_N\in I_M$, and we get \corol\ with
$$C^\eta_{ij}=\sum_{k_1,\ldots,k_N=1}^M
\biggl[\sum_{l=1}^N(-)^lX^\eta_{jl}\;\Xmin{\eta}{\ss
i\,1\hfill\ldots \check{l} \ldots \hfill N\cr \ss k_1\hfill\ldots
\ldots\hfill k_N}+ X^\eta_{ji}\;\Xmin{\eta}{\ss 1\hfill\ldots\hfill
N\cr \ss k_1\hfill\ldots\hfill k_N}\biggr]
{\kdue[\eta\eta_{k_1},\ldots,\eta\eta_{k_N}]\over\kdue[v]}\ ,$$
which is equivalent to \leC\ by the identity
$$\sum_{l=1}^N(-)^lX^\eta_{jl}\;\Xmin{\eta}{\ss i\,1\hfill\ldots \check{l}
\ldots \hfill N\cr \ss k_1\hfill\ldots \ldots\hfill k_N}+
X^\eta_{ji}\; \Xmin{\eta}{\ss 1\hfill\ldots\hfill N\cr \ss
k_1\hfill\ldots\hfill k_N} \,=\,\Xmin{\eta}{\ss i 1 \hfill \ldots
\hfill N \cr \ss j k_1\hfill \ldots \hfill k_N}\ .$$ Eq.\detv\
follows by \vXeta. \hfill$\square$

\vskip 6pt

The homomorphisms $(X^{\eta'}-X^\eta)\in{\rm End}(\CC^M)$, satisfying
\lacond, are the elements of a $M(M-N)$ dimensional vector space,
spanned by
$$(X^{\eta'}-X^\eta)_{ij}=\sum_{k=N+1}^M\Lambda_{jk}C^\eta_{ki}\
,$$ $i,j\in I_M$, with $\Lambda_{jk}$ an arbitrary $M\times (M-N)$
matrix. An obvious generalization of \lemma\ yields
\eqn\teor{\eta\eta_i =\sum_{j=1}^N v_jB^\eta_{ji}\ ,} $i\in I_M$,
implying that $B^\eta_{ij}=
\kdue[v_1,\ldots,v_{j-1},\eta\eta_i,v_{j+1},\ldots,v_N]/\kdue[v]$,
are the matrix elements of the homomorphism $B^\eta=\tilde\psi\circ
(X^\eta)^{-1}$. Such coefficients can be expanded as
\eqn\leB{B^\eta_{ij}=\sum_{k_1,\ldots,k_{N-1}=1}^M(-)^{j+1}
\Xmin{\eta}{\ss 1\ldots \check{\jmath} \ldots N\hfill\cr \ss
k_1\hfill\ldots\hfill k_{N-1}}
{\kdue[\eta\eta_i,\eta\eta_{k_1},\ldots,\eta\eta_{k_{N-1}}]\over\kdue[v]}\
.}

\vskip 6pt

Define $C_{kl}^{\eta,ij}$, $3\le i<j\le g$, $k,l\in I_g$, by
$$C^{\eta,\1_m\2_m}_{\1_n\2_n}:=C^\eta_{mn}\ ,$$
$m,n\in I_M$, $m>N$.

The following result is a direct consequence of the Petri-like
approach. The bound $r\le 6$ for the rank of quadrics is not sharp,
however: M. Green proved that the ideal of quadrics of a canonical
curve is generated by elements of rank $4$ \Green.

\vskip 6pt

\newth\ranko{Theorem}{All the relations among holomorphic quadratic differentials have rank $r\le 6$.}

\vskip 6pt

\noindent{\bf Proof.} The statement is trivial for $g\le 6$, so let us assume $g\ge 7$. Each relation
can be written as
$$\eqalign{0=&\sigma_i\sigma_j+C^{\sigma,ij}_{12}\sigma_1\sigma_2+C^{\sigma,ij}_{1i}\sigma_1\sigma_i+
C^{\sigma,ij}_{1j}\sigma_1\sigma_j+C^{\sigma,ij}_{2i}\sigma_2\sigma_i+
C^{\sigma,ij}_{2j}\sigma_2\sigma_j\cr
&+\sum_{k\neq 1,2,i,j}C^{\sigma,ij}_{1k}\sigma_1\sigma_k+
\sum_{k\neq 1,2,i,j}C^{\sigma,ij}_{2k}\sigma_2\sigma_k\ ,}$$
where $3\le i<j\le g$ and $C^{\eta,\1_i\2_i}_{\1_j\2_j}:=C^\eta_{ij}$. Set $\eta_1\equiv\sigma_1$, $\eta_2\equiv\sigma_2$, $\eta_3\equiv\sigma_i$,
$\eta_4\equiv\sigma_j$, $\eta_5\equiv\sum_{k\neq 1,2,i,j}C^{\sigma,ij}_{1k}\sigma_k$, $\eta_6\equiv
\sum_{k\neq 1,2,i,j}C^{\sigma,ij}_{2k}\sigma_k$. Then the relations can be written as
$$\sum_{k<l}^6C^{\eta,ij}_{kl}\eta_k\eta_l=0\ ,$$
for suitable $C^{\eta,ij}_{kl}$, and the theorem follows. \hfill$\square$

\vskip 6pt

\subsec{Consistency conditions on the quadrics coefficients}

In the construction in section \primecostr, the points $p_1$ and
$p_2$ play a special role with respect to $p_3,\ldots,p_g$.
Relations among holomorphic quadratic differentials can be obtained
by replacing $p_1$ and $p_2$ with $p_a$ and $p_b$, $a,b\in I_g$,
$a<b$, $(a,b)\neq (1,2)$. In the following of this section, we will
consider the relationships between the coefficients $C^\sigma$
obtained in section \primecostr\ and the analogous coefficients
obtained upon replacing $(1,2)$ by $(a,b)$.

\newth\tutti{Proposition}{There exist $g$ distinct points $p_1,\ldots,p_g\in C$ such that
$$K(p_1,\ldots,\check
p_i,\ldots,\check p_j,\ldots,p_g)\neq 0\ ,$$
for all $i,j\in I_g$, $i\neq j$.}

\vskip 6pt

\noindent {\bf Proof.} Consider the function in $C^g$
$$F(p_1,\ldots,p_g):=\prod_{i< j}K(p_1,\ldots,\check
p_i,\ldots,\check p_j,\ldots,p_g)\ ,$$ and set
$Z:=\{(p_1,\ldots,p_g)\in C^g\mid F(p_1,\ldots,p_g)=0\}$. Note that
$Z=\bigcup_{i<j}\{(p_1,\ldots,p_g)\in\C^g\mid K(p_1,\ldots,\check
p_i,\ldots,\check p_j,\ldots,p_g)=0\}$, so that it is a finite union
of varieties of codimension $1$ in $C^g$ and, in particular, $Z\neq
C^g$. Suppose that $C^g\setminus (\bigcup_{i<j}\Pi_{ij})\subseteq
Z$, where $\Pi_{ij}:=\{(p_1,\ldots,p_g)\in C^g\mid p_i=p_j\}$, $1\le
i<j\le g$. Since $C^g\setminus (\bigcup_{i<j}\Pi_{ij})$ is dense in
$C^g$, it would follow that $Z\equiv C^g$, which is absurd. Hence,
there exist pairwise distinct $p_1,\ldots,p_g\in C$ such that
$F(p_1,\ldots,p_g)\neq 0$. \hfill$\square$

\vskip 6pt

By Proposition \tutti\ and Proposition \futuref{thlev}, one can choose the points $p_1,\ldots,p_g$
in such a way that
$$\{v^{(ab)}_i\}_{i\in I_N}:=\{\sigma_i^2\}_{i\in
I_g}\cup\{\sigma_a\sigma_b\}\cup\{\sigma_a\sigma_i,\sigma_b\sigma_i\}_{i\in I_g\setminus\{a,b\}}\ ,$$
is a basis of $H^0(K_C^2)$. Furthermore, one can obtain $M-N$ independent linear relations
\eqn\lerrela{\sum_{1\le k\le l\le g}(ab)^{ij}_{kl}\sigma_k\sigma_l=0\ ,}
where $i,j\in I_g\setminus\{a,b\}$, $i\neq j$. The coefficients $(ab)^{ij}_{kl}$ are defined by
setting $(ab)_{ij}^{ij}:=1$,
\eqn\leCCi{(ab)^{ij}_{kl}:={\kdue[v^{(ab)}_1,\ldots,\check\sigma_k\check\sigma_l,\sigma_i\sigma_j,\ldots,v^{(ab)}_N]\over
\kdue[v^{(ab)}_1,\ldots,v^{(ab)}_N]}\ ,}
if $k\neq l$ and $\sigma_k\sigma_l\in\{v^{(ab)}_i\}_{i\in I_N}$, and $(ab)^{ij}_{kl}:=0$ for all the
other $(k,l)\in I_g\times I_g$.
In this notation, the coefficients $C^\sigma_{ij}$ defined in \futuref{leC}, with $N< i\le M$, $j\in I_M$, correspond to
$(12)^{\1_i\2_i}_{\1_j\2_j}$. Eqs.\lerrela\ and \leCCi\ can be derived by a trivial generalization
of the same construction
considered in section 2 in the particular case $a=1$, $b=2$.

\newth\syzy{Proposition}{The coefficients $(ab)^{ij}_{kl}$ satisfy the following consistency
conditions
\eqn\lesizy{\eqalign{(ij)^{ab}_{kl}=\sum_{m\le n}(ij)^{ab}_{mn}(ab)^{mn}_{kl}&=\sum_{m\le
n}(ij)^{ab}_{mn}(ai)^{mn}_{kl}=\sum_{m\le n}(ij)^{ab}_{mn}(aj)^{mn}_{kl}\cr
&=\sum_{m\le
n}(ij)^{ab}_{mn}(bi)^{mn}_{kl}=\sum_{m\le n}(ij)^{ab}_{mn}(bj)^{mn}_{kl}\ ,}}
for all $i,j,a,b\in I_g$ pairwise distinct, and for all $k,l\in I_g$.
}

\vskip 6pt

\noindent {\bf Proof.}
Choose $i,j,a,b\in I_g$, with $a<b<i<j$, and consider the relations $\sum_{k\le l}
(ij)^{ab}_{kl}\sigma_k\sigma_l=0$ and $\sum_{k\le
l}(ab)^{ij}_{kl}\sigma_k\sigma_l=0$, that is
\eqn\sizyidea{\eqalign{0=(ij)^{ab}_{ij}\sigma_i\sigma_j+\sigma_a\sigma_b&+(ij)^{ab}_{ai}\sigma_a\sigma_i+
(ij)^{ab}_{aj}\sigma_a\sigma_j+(ij)^{ab}_{bi}\sigma_b\sigma_i\cr & +(ij)^{ab}_{bj}\sigma_b\sigma_j
+\sum_{k\neq a,b,i,j}(ij)^{ab}_{ik}\sigma_i\sigma_k +\sum_{k\neq
a,b,i,j}(ij)^{ab}_{jk}\sigma_j\sigma_k\ ,}}
\eqn\sizyideb{\eqalign{0=\sigma_i\sigma_j+
(ab)^{ij}_{ab}\sigma_a\sigma_b&+(ab)^{ij}_{ai}\sigma_a\sigma_i+
(ab)^{ij}_{aj}\sigma_a\sigma_j+(ab)^{ij}_{bi}\sigma_b\sigma_i\cr &+ (ab)^{ij}_{bj}\sigma_b\sigma_j
+\sum_{k\neq a,b,i,j}(ab)^{ij}_{ak}\sigma_a\sigma_k +\sum_{k\neq
a,b,i,j}(ab)^{ij}_{bk}\sigma_b\sigma_k\ .
}}
Replace the differentials $\sigma_i\sigma_k$ and $\sigma_j\sigma_k$, $k\neq
i,j,a,b$, in Eq.\sizyidea\ by
$$\sigma_i\sigma_k=-\sum_{^{m\le n}_{(m,n)\neq (i,k)}}(ab)^{ij}_{mn}\sigma_m\sigma_n\ ,\quad k\neq
i,j,a,b\ ,$$
and the analogous expression for $\sigma_j\sigma_k$. Then multiply Eq.\sizyideb\ by
$(ij)^{ab}_{ij}$ and consider the difference between \sizyidea\ and \sizyideb. We obtain
\eqn\quasisizy{\eqalign{0=&\Bigl((ij)^{ab}_{ab}-\sum_{m\le n}(ij)^{ab}_{mn}(ab)^{mn}_{ab}\Bigr)\sigma_a\sigma_b+
\sum_{k\neq a,b}\Bigl((ij)^{ab}_{ak}-\sum_{m\le
n}(ij)^{ab}_{mn}(ab)^{mn}_{ak}\Bigr)\sigma_a\sigma_k\cr
&+\sum_{k\neq a,b}\Bigl((ij)^{ab}_{bk}-\sum_{m\le
n}(ij)^{ab}_{mn}(ab)^{mn}_{bk}\Bigr)\sigma_b\sigma_k \ .
}}
Since the holomorphic quadratic differentials appearing in Eq.\quasisizy\ are linearly
independent, it follows that each coefficient vanishes, yielding the first identity in \lesizy, in
the cases in which at least one between $k$ and $l$ is equal to $a$ or $b$. On the other hand, in the case
$k,l\neq a,b$, the only non-vanishing term in the sum $\sum_{m\le n}(ij)^{ab}_{mn}(ab)^{mn}_{kl}$ is
$(ij)^{ab}_{kl}(ab)^{kl}_{kl}=(ij)^{ab}_{kl}$, and the first identity in \lesizy\ follows. The
other identities can be proved by applying the analogous procedure to the relation $\sum_{k\le
l}(ij)^{ab}_{kl}\sigma_k\sigma_l=0$ and one of the relations $\sum_{k\le
l}(ai)^{bj}_{kl}\sigma_k\sigma_l=0$, $\sum_{k\le
l}(bi)^{aj}_{kl}\sigma_k\sigma_l=0$, and so on.
\hfill$\square$

\vskip 6pt

\subsec{A correspondence between quadrics and $\theta$-identities}

\newth\thleX{Theorem}{Fix $p_1,\ldots,p_g$ satisfying the equivalent conditions {\it i}), {\it ii}), {\it iii}) of Theorem \quadgen.
Then, the associated holomorphic quadratic differentials $v_i$,
$i\in I_M$, satisfy
\eqn\vXww{v_i=\sum_{j=1}^MX^\omega_{ji}\,\omega\omega_j\ ,} $i\in
I_N$, where \eqn\leX{X^\omega_{ij}={\theta_{\Delta,\1_j}(\a_{\1_i})
\theta_{\Delta,\2_j}(\a_{\2_i})+\theta_{\Delta,\1_j}(\a_{\2_i})
\theta_{\Delta,\2_j}(\a_{\1_i})\over
(1+\delta_{\1_j\2_j})\sum_{l,m}\theta_{\Delta,l}(\a_{\1_i})
\theta_{\Delta,m}(\a_{\2_i})\omega_l(p_{\1_i})\omega_m(p_{\2_i})}\
,} $i,j\in I_M$, with $\a_i$ as in Definition \defAB, correspond to
the coefficients defined in \XXXX\ for $\eta_i\equiv\omega_i$, $i\in
I_g$. Furthermore, the $M-N$ independent linear relations
\eqn\coroll{\sum_{j=1}^MC^\omega_{ij}\omega\omega_j=0\ ,} $N+1\le
i\le M$, hold, where
\eqn\leCi{C^\omega_{ij}=\sum_{k_1,\ldots,k_N=1}^M \Xmin{\omega}{\ss
1 \hfill \ldots \hfill N i \cr \ss k_1\hfill \ldots \hfill k_N
j\hfill}\;
{\kdue[\omega\omega_{k_1},\ldots,\omega\omega_{k_N}]\over\kdue[v]}\
. } correspond to the coefficients defined in \leC. }

\vskip 6pt

\noindent {\bf Proof.} Eq.\leomega\ implies that Eq.\leX\ is equivalent to \XXXX, and the
theorem follows by Theorem \thcorol. \hfill$\square$

\vskip 6pt

\newrem\rimarcone{Remark}{Choose $p_1,\ldots,p_g$ as in Corollary \rango, with $n=1$ and set
$q:=q_1$. Then, there exists a non-trivial relation
$$a\sigma_1t_2+b\sigma_2 t_1+c\sigma_1\sigma_2=0\ ,
$$
where $a,b,c\in \CC$. Without loss of generality, we can assume that
$t_1(p_1)=0$ and $t_2(p_2)=0$. Set
$$(\sigma_1)=p_2+p_3+\ldots+p_g+q+\sum_{i=1}^{g-2}r_i\ ,$$
and
$$(\sigma_2)=p_1+p_3+\ldots+p_g+q+\sum_{i=1}^{g-2}s_i\ ,$$
for some $r_i,s_i\in C$, $i\in I_{g-2}$.
Then, $(t_1)>p_1+\sum_{i=1}^{g-2}r_i$ and $(t_2)>p_2+\sum_{i=1}^{g-2}s_i$, so that
$$\eqalign{t_1&\sim-{\deltadiv(p_1+\sum_ir_i+z-y)\over
\sigma(y)E(y,z)E(y,p_1)\prod_iE(y,r_i)}\sigma(z)E(z,p_1)\prod_iE(z,r_i)\cr
&\sim-{\deltadiv(p_2+\b+q+y-p_1-z)\over
\sigma(y)E(y,z)E(y,p_1)\prod_iE(y,r_i)}\sigma(z)E(z,p_1)\prod_{i=3}^gE(z,r_i)\cr
&=\sum_{i=1}^g\theta_{\Delta,i}(p_2+\b+q-p_1)\omega_i(z)\ ,}$$
where, in the second line, we used $I(\sum_i
r_i)=I(2\Delta-p_2-\b-q)$ in $J_0(C)$. An analogous calculation
yields
$$t_2\sim\sum_{i=1}^g\theta_{\Delta,i}(p_1+\b+q-p_2)\omega_i(z)\ .
$$ (By the symbol $\sim$, we denote the equality up to a factor independent of $z$; such a factor is not meaningful, since it can be compensated
by a redefinition of
the constants $a,b$.)}

\vskip 6pt

\newth\main{Theorem}{Let $C$ be a canonical curve of genus $g\geq4$ and $\{\omega_i\}_{i\in I_g}$ the
canonically normalized basis of $H^0(K_C)$, and let the points
$p_3,\ldots,p_g\in C$ satisfy
one of the equivalent conditions {\it i}), {\it ii}) and {\it iii})
in Theorem \quadgen. Then, the following $(g-2)(g-3)/2$ independent relations
\eqn\bene{\eqalign{\sum_{s\in\perm_{2g}}\sgn(s)&\det\omega(x_{s_1},\ldots,x_{s_g})
\det\omega(x_{s_g},\ldots,x_{s_{2g-1}}) \det
\omega(x_{s_1},x_{s_{g+1}},x_{s_{2g}},p_{3},\ldots,\check
p_i,\ldots,p_g)\cr &\cdot\det
\omega(x_{s_2},x_{s_{g+2}},x_{s_{2g}},p_{3},\ldots,\check
p_j,\ldots,p_g) \prod_{k=3}^{g-1}\det
\omega(x_{s_k},x_{s_{k+g}},p_{3},\ldots,p_g)=0\ , }} $3\leq i<j\leq
g$, hold for all $x_k\in C$, $1\le k\le 2g$. }

\vskip 6pt

\noindent {\bf Proof.} Fix $i,j$, $3\le i<j\le g$, and choose
$p_1,p_2$ in such a way that $\{\sigma_i\}_{i\in I_g}$ is a basis of
$H^0(K_C)$. Observe that, due to Eq.\lemma, $\det_I
\sigma\sigma(x_1,\ldots,x_{2g})=0$, for all $x_1,\ldots,x_{2g}\in
C$, where $I:=I_{M,2}\cup\{m(i,j)\}$. Applying Lemma \thcombvi, with
$n=2$, such an identity corresponds to Eq.\bene\ with the canonical
basis $\{\omega_i\}_{i\in I_g}$ of $H^0(K_C)$ replaced by
$\{\sigma_i\}_{i\in I_g}$. Eq.\bene\ is then obtained by simply
changing the base. \hfill$\square$

\vskip 6pt

The relations of Theorem \main\ can be directly expressed in terms
of theta functions.

\vskip 6pt

\newth\ththetarel{Theorem}{Fix $p_3,\ldots,p_g\in C$ in such
a way that the equivalent conditions {\it iv}), {\it v}), {\it vi}),
and {\it vii}) of Theorem \quadgen\ are satisfied. The following
$(g-2)(g-3)/2$ independent relations
\eqn\thetarel{\eqalign{V_{i_1i_2}(p_3,\ldots,p_g,x_1,\ldots,x_{2g}):=&\cr\cr
\sum_{s\in\perm_{2g}}\sgn(s) \Biggl\{&\prod_{k=1}^2{S(\hat
x_{k}+\hat x_{{g+k}}+\hat x_{{2g}}+\b_{i_k})E(\hat x_{k},\hat
x_{{2g}}) E(\hat x_{{k+g}},\hat x_{{2g}})\over E(\hat
x_{k},p_{i_k})E(\hat x_{{k+g}},p_{i_k})E(\hat
x_{{2g}},p_{i_k})}\cr\cr &\cdot\prod_{k=1}^{g-1}\bigl(E(\hat
x_{k},\hat x_{{k+g}}) \prod_{j=3}^gE(\hat x_{k},p_j) E(\hat
x_{{k+g}},p_j)\bigr) \cr\cr &\cdot S\bigl(\smsum_{k=1}^g\hat
x_{k}\bigr) \prod_{_{k<j}^{k,j=1}}^gE(\hat x_{k},\hat x_{j})
S\bigl(\smsum_{k=g}^{2g-1}\hat
x_{k}\bigr)\prod_{_{k<j}^{k,j=g}}^{2g-1}E(\hat x_{k},\hat x_{j})
\cr\cr &\cdot\prod_{k=3}^{g-1}S(\hat x_{k}+\hat
x_{{k+g}}+\b)\prod_{j=3}^gE(\hat x_{{2g}},p_j)^2\Biggr\}=0 \
,}}$3\le i_1<i_2\le g$, where $\hat x_i:=x_{s_i}$, $i\in I_{2g}$,
$\b_i:=\b-p_i$, $3\le i\le g$, hold for all $x_i\in C$, $i\in
I_{2g}$.}

\vskip 6pt

\noindent {\bf Proof.} By \dettheta\ $V_{ij}(p_3,\ldots,p_g,x_1,\ldots,x_{2g})$
is equivalent to \bene. \hfill$\square$

\vskip 6pt

\newrem\fhgtgfgfg{Remark}{Note that $V_{ii}\neq 0$ for $i=3,\ldots,g$, since for $i=j$ the
LHS of \bene\ is proportional to a determinant of $2g$ linearly
independent holomorphic quadratic differentials on $C$, evaluated at
general points $x_i\in C$, $i\in I_{2g}$.}

\vskip 6pt

By a limiting procedure we derive the original Petri's relations,
now written in terms of the canonical basis $\{\omega_i\}_{i\in
I_g}$ of $H^0(K_C)$ and with the coefficients expressed in terms of
theta functions.

\vskip 6pt

\newth\threla{Corollary}{Fix $p_1,\ldots,p_g\in C$ in such a way that $p_3,\ldots,p_g$
satisfy the equivalent conditions {\it iv}), {\it v}), {\it vi}),
and {\it vii}) of Theorem \quadgen. The following $(g-2)(g-3)/2$
linearly independent relations \eqn\relai{\sum_{j=1}^M
C^\omega_{ij}\omega\omega_j(z):={\kuno[\sigma]\over\kdue[v]}^{g+1}F(p,x){V_{\1_i\2_i}(p_3,\ldots,p_g,x_1,\ldots,x_{2g-1},z)\over
\deltadiv\bigl(\sum_{1}^{2g-1}x_j+\b\bigr)}=0 \ ,} $N+1\le i\le M$,
where
$$
F(p,x):=c'_{g,2}{\prod_{_{j<k}^{j,k=3}}^gE(p_j,p_k)^{g-4}\prod_{_{j\neq
\1_i}^{j=3}}^gE(p_{\1_i},p_j) \prod_{_{j\neq
\2_i}^{j=3}}^gE(p_{\2_i},p_j)\over
\prod_{j=1}^{2g-1}(\sigma(x_j)\prod_{k=3}^gE(x_j,p_k)
\prod_{k=j+1}^{2g-1}E(x_j,x_k))}\ ,
$$
hold for all $z\in C$. Furthermore, $C^\omega_{ij}$ are independent of
$p_1,p_2,x_1,\ldots,x_{2g-1}\in C$ and correspond to the coefficients defined in \futuref{leC}
(with $\eta_i\equiv\omega_i$, $i\in I_g$) or,
equivalently, in \futuref{leCi}.}

\vskip 6pt

\noindent {\bf Proof.} Consider the identity
\eqn\ratttio{{\det_I \sigma\sigma(x_1,\ldots,x_{2g-1},z) \over \det
v(p_3,\ldots,p_g,x_1,\ldots,x_{2g-1})}=0\ ,} $I:=I_{M,2}\cup \{i\}$,
$N+1\le i\le M$. Upon applying Lemma \thcombvi, with $n=2$, and
Eq.\dettheta\ to the numerator and Eq.\detthetaii\ to the
denominator of \ratttio, Eq.\relai\ follows by a trivial
computation. On the other hand, for arbitrary points $z,y_1,\ldots,y_{g-1}\in C$,
$$S(y_1+\ldots+y_{g-1}+z)={
\sum_{i=1}^g\theta_{\Delta,i}(y_1+\ldots+y_{g-1})\omega_i(z)\over
\sigma(z)\prod_{1}^{g-1}E(z,y_i)}\ .$$ Upon replacing each term of
the form $S(\d_{g-1}+z)$ in
$V_{\1_i\2_i}(p_3,\ldots,p_g,x_1,\ldots,x_{g-1},z)$ by its
expression above, for any effective divisor $\d_{g-1}$ of degree
$g-1$, the dependence on $z$ only enters through
$\omega_i\omega_j(z)$ and the relations \relai\ can be written in
the form of Eq.\corol.

\vskip 6pt

To prove that $C^\omega_{ij}$ are the coefficients in \leC, with $\eta_i\equiv\omega_i$, $i\in I_g$,
first consider the identity
$${\kdue[\omega\omega_{k_1},\ldots,\omega\omega_{k_N}]\over\kdue[v]}=
{\det_{i\in\{k_1,\ldots,k_N\}}\omega\omega_i(p_3,\ldots,p_g,x_1,\ldots,x_{2g-1})\over \det
v(p_3,\ldots,p_g,x_1,\ldots,x_{2g-1})}\ ,$$
then recall that
$$v_i:=\sigma\sigma_i=\sum_{j=1}^MX^\omega_{ji}\omega\omega_j\ ,$$
$i\in I_M$, so that one obtains
$$\sum_{k_1,\ldots,k_N,j=1}^M
\Xmin{\omega}{\ss 1 \hfill \ldots \hfill N i \cr \ss k_1\hfill
\ldots \hfill k_N j\hfill}\;
{\kdue[\omega\omega_{k_1},\ldots,\omega\omega_{k_N}]\over\kdue[v]}\omega\omega_j(z)={\det_I
\sigma\sigma(x_1,\ldots,x_{2g-1},z) \over \det
v(p_3,\ldots,p_g,x_1,\ldots,x_{2g-1})}\ ,$$ as an algebraic identity
(in the sense that it holds as an identity in ${\rm
Sym}^2(H^0_C(K))$ after replacing
$\sigma_i\sigma_j\to(\sigma_i\otimes\sigma_j)_S$ and
$\omega_i\omega_j\to(\omega_i\otimes\omega_j)_S$, $i,j\in I_g$).
Hence, the coefficients of $\omega\omega_j(z)$ on the LHS, given by
\leC\ or, equivalently, by \leCi\ and the ones on the RHS, given by
\relai, are the same.

Eq.\leC\ explicitly shows that the coefficients $C^\omega_{ij}$ are
independent of $x_1,\ldots,x_{2g-1}$. By \relai\ it follows that
they may depend on $p_1$ and $p_2$ only through the term
$\kuno[\sigma]^{g+1}/\kdue[v]$. The dependence of $\kuno[\sigma]$
and $\kdue[v]$ on $p_1$ and $p_2$ is due to the dependence of the
basis $\{\sigma_i\}_{i\in I_g}$ and $\{v_i\}_{i\in I_N}$ on the
choice of $p_1,\ldots,p_g\in C$. On the other hand, Eq.\Kkappa\
implies that $\kuno[\sigma]^{g+1}/\kdue[v]$ is independent of
$p_1,p_2$ and the proof of the corollary is complete.
\hfill$\square$

\vskip 6pt

\subsec{$K=0$ as a quadric from a double point on $\Theta_s$}

Choose $p_3,\ldots,p_g\in C$ pairwise distinct and such that
$K(p_3,\ldots,p_g)\neq 0$. Let $C_2\ni \c:=u+v$, $u,v\in C$, be an
effective divisor of degree $2$, such that $u$ is distinct from
$p_3,\ldots,p_g$ and $\sum_{i=3}^gp_i+\c$ is special. Then there
exists $x\in C$ such that $(x,u,p_3,\ldots,p_g)\in C^g\setminus\A$
(or, otherwise, $K(p_3,\ldots,p_g)$ would vanish); let
$\{\sigma_i\}_{i\in I_g}$ be the basis of $H^0(K_C)$ associated to
$x,u,p_3,\ldots,p_g$ by Proposition \futuref{thnewbasis}.

Let $A(\c)\subset I_g\setminus\{1,2\}$ be the set
$$A(\c):=\{i\in I_g\setminus\{1,2\}\mid \sigma_i(v)\neq 0\}\ ,$$
and $\bar
A(\c):=\{3,\ldots,g\}\setminus A(\c)$ its complement.

\vskip 6pt

\newth\basepoints{Lemma}{The set $A(\c)$ is independent of $x$, provided that
$(x,u,p_3,\ldots,p_g)\in C^g\setminus \A$. Furthermore, for each
subset $A'\subseteq I_g\setminus\{1,2\}$, the divisor $\sum_{i\in
A'}p_i+\c$ is special if and only if $A(\c)\subseteq A'$, and
$A(\c)$ is the unique set satisfying such a property.}

\vskip 6pt \noindent{\bf Proof.} An effective divisor $D$, with $\deg D\le g$, is special if and only
if $h^0(K_C\otimes\O(-D))>g-\deg D$. Consider the divisor $\d:=\sum_{i\in A(c)}p_i+\c$ of degree
$\deg\d=a+2$, where $a$ is the cardinality of $A(\c)$. Since $H^0(K_C\otimes\O(-\d))$ is generated
by $\sigma_1$ and by the elements of $\{\sigma_i\}_{i\in \bar A(c)}$,
$$h^0(K_C\otimes\O(-\d))=g-1-a>g-2-a=g-\deg\d\ ,$$
and $\d$ is special. It follows that if $A(\c)\subseteq A'\subseteq \{3,\ldots,g\}$, then
$\sum_{i\in A'}p_i+\c\ge \d$ is special.

Conversely, set $\d:=\sum_{i\in A'}p_i$ and suppose that $\d+\c$ is special. Note that,
since $\d+u$ is not special,
$$h^0(K_C\otimes\O(-\d-u))=g-\deg\d-1\le h^0(K_C\otimes\O(-\d-\c))\ ,$$ and by
$H^0(K_C\otimes\O(-\d-\c))\subseteq H^0(K_C\otimes\O(-\d-u))$, it follows that
$H^0(K_C\otimes\O(-\d-\c))=H^0(K_C\otimes\O(-\d-u))$; in other words, each element of
$H^0(K_C\otimes\O(-\d-u))$ also vanishes at $v$. Now,
$H^0(K_C\otimes\O(-\d-u))$ is generated by $\sigma_1$ and by the elements of $\{\sigma_i\}_{i\in
\bar A'}$, where $\bar A':=\{3,\ldots,g\}\setminus A'$. Then, $\sigma_i(v)=0$ for all $i\in \bar
A'$, so that $\bar A'\subseteq \bar A(\c)$ and then $A(\c)\subseteq A'$.

Uniqueness follows by noting that if $\tilde A$ satisfies the same property, then
$\tilde A\subseteq A(\c)$ (because $\sum_{i\in A(c)}p_i+\c$ is special) and $A(\c)\subseteq\tilde A$
(because $\tilde A\subseteq\tilde A$ implies that $\sum_{i\in \tilde A}p_i+\c$ is special).

Finally, by defining $A(\c)$ as the unique set satisfying such a property, it follows that $A(\c)$
is independent of $x$.\hfill$\square$

\vskip 6pt

\newth\orderzero{Lemma}{Suppose that $\bar A(\c)\neq\emptyset$ and fix $i\in \bar A(\c)$ and
$j\neq i$, $3\le j\le g$. Let $k+1$, with $k\ge 0$, be the order of the zero of $\sigma_1$ in $p_j$. Then, the holomorphic $1$-differential
$$\lambda^{(c)}_i(z):=\sum_{a,b\in I_g}\theta_{ab}(\c+\sum_{l\neq
i}p_l-\Delta)\omega_a(p_i)\omega_b(z)\ ,$$ has a zero of order $n\ge
k$ in $z=p_j$, and $n>k$ if and only if $j\in \bar A(\c)$.}

\vskip 6pt

\noindent{\bf Proof.} Define the points $\tilde x_1,\ldots,\tilde
x_{g-2-k}$ by
$$(\sigma_1)=\sum_{l=3}^gp_l+u+v+kp_j+\sum_{l=1}^{g-2-k}\tilde x_l\ ,$$
so that $I(\sum_{l=3}^gp_l+u+v+kp_j+\sum_{l=1}^{g-2-k}\tilde
x_l-2\Delta)=b+\tau a$, for some $a,b\in \ZZ^g$. Consider the
identities
$$\eqalign{&
\sum_{l\in I_g}\theta_{l,\Delta}(u+v+\sum_{m=3}^gp_m-w)
\omega_l(z)\cr\cr &=-\eta_ae^{-2\pi i\tp
aI(\Delta-w-kp_j-\sum_m\tilde x_m)}{\deltadiv(\sum_{m}\tilde
x_m+kp_j +w+z-y)E(z,p_j)^kE(z,w)\prod_l E(z,\tilde
x_l)\sigma(z)\over E(y,z)E(y,w)E(y,p_j)^k\prod_l E(y,\tilde
x_l)\sigma(y) }\cr\cr &=e^{-2\pi i\tp aI(z-y)}
{\deltadiv(y+u+v+\sum_{m}p_m-w-z)E(z,p_j)^kE(z,w)\prod_l E(z,\tilde
x_l)\sigma(z)\over E(y,z)E(y,w)E(y,p_j)^k\prod_l E(y,\tilde
x_l)\sigma(y) }\ ,}$$ where $\eta_a:=e^{\tp a\tau a}$, which hold
for arbitrary $w,y\in C$. Dividing by $E(p_i,w)$ and taking the
limit $w\to p_i$ one obtains
$$\lambda_i^{(c)}(z)={e^{-2\pi i\tp aI(z-y)}E(z,p_j)^kE(p_i,z)\prod_{m}
E(z,\tilde x_m)\sigma(z)\over E(y,z)E(y,p_i)E(y,p_j)^k\prod_m
E(y,\tilde x_m)\sigma(y)}\sum_{l\in
I_g}\theta_{l,\Delta}(y+u+v+\sum_{m\neq i}p_m-z)\omega_l(p_i)\ .
$$
Since the right hand side does not depend on $y$, the factor
$E(z,p_j)^k$ cannot be compensated by any factor in the denominator
and the $1$-differential has a zero of order at least $k$ in
$z=p_j$. Furthermore, such a zero if of order strictly greater than
$k$ if and only if
$$\sum_{l\in
I_g}\theta_l(y+u+v+\sum_{m\in A'}p_m-\Delta)\omega_l(p_i)=0\ ,$$ for
all $y\in C$, with $A':=\{3,\ldots,g\}\setminus \{i,j\}$. In
particular, for $y\equiv x$, this implies that the holomorphic
1-differential
$$\sum_{l\in
I_g}\theta_l(x+u+v+\sum_{m\in A'}p_m-\Delta)\omega_l(z)\ ,$$
vanishes at $p_i$. Therefore, such a differential vanishes at
$x,u,v$ and $p_l$, for all $l\neq j$, $3\le l\le g$; hence, it is
proportional to $\sigma_j$, which is the generator of
$H^0(K_C\otimes\O(-u-x-\sum_{l\neq j}p_l))$, and it must be
$\sigma_j(v)=0$, so that $j\in\bar A(\c)$.
Conversely, if $j\in \bar A(\c)$, then $A(\c)\subseteq A'$ and,
by Lemma \basepoints, $y+u+v+\sum_{l\in A'}p_l$
is a special divisor for all $y\in C$. Then, for each $y\in C$, there exist
$q_1,\ldots,q_{g-2}\in C$ such that $I(y+u+v+\sum_{l\in A'}
p_l)=I(p_i+\sum_lq_l)$, so that
$$\sum_{l\in
I_g}\theta_l(y+u+v+\sum_{m\neq
i,j}p_m-\Delta)\omega_l(p_i)=\sum_{l\in
I_g}\theta_l(p_i+\sum_mq_m-\Delta)\omega_l(p_i)=0\ ,
$$
for all $y\in C$, and the lemma follows. \hfill $\square$

\vskip 6pt

Set
\eqn\iLambdaij{\Lambda^{(i)}_{jk}(\c):=\sum_{a,b\in
I_g}\theta_{ab}(\c+\smsum_{l\neq
i}p_l-\Delta)\omega_a(p_j)\omega_b(p_k)\ ,}
$i,j,k\in I_g\setminus \{1,2\}$. Note that, if $i\in\bar A(\c)$, then $\Lambda^{(i)}_{jk}(\c)=0$ for $j=k$ and for $j,k\neq i$, and
$$\Lambda^{(i)}_{ij}(\c)=\lambda^{(c)}_i(p_j)\ ,
$$
$j\neq i$.

\vskip 6pt

\newth\coefflambda{Theorem}{Choose $p_3,\ldots,p_g\in C$, $C_2\ni \c:=u+v$ and $x\in C$ as above.
Suppose $\bar A(\c)\neq\emptyset$ and fix $i\in\bar A(\c)$. If $u$
is a single zero for $K(\ \cdot\ ,p_3,\ldots,\check p_i,\ldots,
p_g)$, then the holomorphic quadratic differentials
$\sigma\sigma_k$, $k\in I^{1i}_N$ (see Definition
\definsiemi\ for notation), satisfy a unique linear relation
$$\sum_{k\in I^{1i}_N}\tilde C^{\sigma(i)}_{k}(\c)\sigma\sigma_k=0\ ,$$
where
$$\tilde
C^{\sigma(i)}_{k}(\c)=\sum_{^{j\in
I^i_2}_{j>N}}\Lambda^{(i)}_{j}(\c)C^{\sigma}_{jk} \ ,$$ $k\in
I^{1i}_N$, with
$\Lambda_{j}^{(i)}(\c):=\Lambda^{(i)}_{\1_j\2_j}(\c)$, $j\in I_M$,
defined in Eq.\iLambdaij.}

\vskip 6pt

\noindent{\bf Proof.} By Theorem \zeriK\ and Corollary \rango, since $u$ is a single zero of $K(\
\cdot\ ,p_3,\ldots,\check p_i,\ldots, p_g)$, then $\sigma\sigma_k$, $k\in I^{1i}_N$, span a
$(N-1)$-dimensional vector space in $H^0(K^2_C)$, and then satisfy a relation
$$\sum_{k\in I^{1i}_N}\tilde C^{\sigma(i)}_{k}(\c)\sigma\sigma_k=0\ .$$
Such a relation determines, up to normalization, an element
$$\ker\psi\ni\phi:=\sum_{k\in I^{1i}_N}\tilde C^{\sigma(i)}_{k}\sigma\cdot\sigma_k\ ,
$$ where $\psi:\Sym^2H^0(K_C)\to H^0(K_C^2)$;
by Theorem \thcorol, $\ker\psi$ is spanned by
$\{\sum_{k=1}^MC^\sigma_{ik}\sigma\cdot\sigma_k\}_{N<i\le M}$, so
that \eqn\relinsym{\sum_{k\in I_N^{1i}}\tilde
C^{\sigma(i)}_{k}(\c)\sigma\cdot\sigma_k=
\sum_{j=N+1}^ML^{(i)}_{j}(\c)\sum_{l\in
I_M}C^{\sigma}_{jl}\sigma\cdot\sigma_l\ ,} for some complex
coefficients $L^{(i)}_{j}(\c)$, $N<j\le M$. Note that, for all
$j,k$, with $N<j,k\le M$, $C^{\sigma}_{jk}=\delta_{jk}$. Then, by
applying $(p\cdot p)_j$ (see Eq.\pqdef), $j=N+1,\ldots,M$, to both
sides of \relinsym, and by using Eq.\pqdual, we obtain
$$L^{(i)}_{j}=\left\{\vcenter{\vbox{\halign{\strut\hskip 6pt $ # $ \hfil & \hskip
5cm$ # $ \hfil\cr \tilde C^{\sigma(i)}_{j}(\c)\ , & \hbox{ for }
j\in I^{1i}_N\ ,\cr 0\ ,& \hbox{ for }j\not\in I^{1i}_N\
,\cr}}}\right.\qquad\qquad N<j\le M \ .$$ Observe that if $j\in
I^{1i}_N$ and $j>N$, then $j\in I^i_2$ (see Def.
\definsiemi), that is, at least one
between $\1_j$ and $\2_j$ is equal to $i$; furthermore, the
condition $j>N$ implies $\1_j\neq \2_j$ and $\1_j,\2_j\neq 1,2$.
Therefore, it remains to prove that $L^{(i)}_{j}(\c)\equiv
C^{\sigma(i)}_j =\Lambda^{(i)}_{\1_j\2_j}(\c)$ for all $j\in I^i_2$,
$j>N$, with respect to a suitable normalization of $\phi$.

The vector $\phi$ can be expressed as \eqn\relperx{\phi\equiv
\sum_{k\in I_N^{1i}}\tilde C^{\sigma(i)}_{k}(\c)\sigma\cdot\sigma_k=
\sigma_1\cdot\eta+\sigma_i\cdot\rho+c\sigma_1\cdot\sigma_i\ , } for
some $\eta,\rho\in H^0(K_C)$, $c\in \CC$, so that the relation
$\psi(\phi)=0$ corresponds to
\eqn\relaperx{\sigma_1\eta+\sigma_i\rho+c\sigma_1\sigma_i=0\ .} Note
that, by the redefinition $\eta\to\eta+\alpha\sigma_i$, $c\to
c-\alpha$, for a suitable $\alpha\in\CC$, we can assume
$\eta(p_i)=0$. Applying $p_i\cdot p_j$, $3\le j\le g$, $j\neq i$, to
both sides of \relperx,  it follows that
$$L^{(i)}_{ij}(\c)=p_i\cdot p_j[\phi]=\rho(p_j)\ ,$$
where
$L^{(i)}_{\1_k\2_k}(\c)=L^{(i)}_{\2_k\1_k}(\c):=L^{(i)}_k(\c)$,
$N<k\le M$. Define $\d\in C_{g-2}$ in such a way that
$$(\sigma_1)=\b+\c+\d\ ,$$ and observe that, by \relaperx,
$\rho\in H^0(K_C\otimes\O(-\d))$ (since $u$ is a single zero for
$K(\cdot,p_3,\ldots,\check p_i,\ldots,p_g)$, it follows that the
$\gcd$ of $(\sigma_1)$ and $(\sigma_i)$ is $\c+\sum_{k\neq i}p_k$).
Furthermore, $\rho$ cannot be a multiple of $\sigma_1$, since, in
this case, the only possibility for Eq.\relaperx\ to hold would be
$\phi=0$. Finally, $L^{(i)}_{ij}(\c)$ is invariant under the
redefinition $\rho\to \rho+a\sigma_1$, since $\sigma_1(p_j)=0$ for
all $j=3,\ldots,g$. Then, we can fix an arbitrary $y\in
C\setminus\supp(\sigma_1)$ and assume that $\rho$ is an element of
the $1$-dimensional space $H^0(K_C\otimes \O(-\d-y))$. By using the
relation $I(\b+\c-y-\Delta)=-I(\d+y-\Delta)$, such an element can be
expressed as follows \eqn\birho{\rho(z)= {a(y)\over A}{\sum_{k\in
I_g}\theta_k(\b+\c-y-\Delta)\omega_k(z)\over E(y,p_i)}\ ,} where the
normalizing constant $A$ can be arbitrarily fixed, and $a$ is a
function such that \eqn\lirho{L^{(i)}_{ij}={a(y)\over A}{\sum_{k\in
I_g}\theta_k(\b+\c-y-\Delta)\omega_k(p_j)\over E(y,p_i)}\ ,} $3\le
j\le g$, $j\neq i$, is independent of $y$. In other words, we assume
that, under the change
$$\eqalign{y&\to \tilde y\ ,\qquad y,\tilde y\in C\setminus\supp(\sigma_1)\cr
\rho &\to\tilde\rho\ ,}$$
$\rho(p_i)$ is equal to $\tilde\rho(p_i)$; this property, together with the fact that
$\tilde\rho\in H^0(K^C\otimes\O(-\d))$, which is generated by $\sigma_1$ and $\rho$, implies that
\eqn\ydep{\tilde\rho=\rho+f(y,\tilde y)\sigma_1\ ,
} for some function $f$. Though Eq.\birho\ only holds for $y\in C\setminus\supp(\sigma_1)$,
the RHS of Eq.\lirho\ is a constant and can be continued to all $y\in C$ and, in particular, in the
limit $y\to p_i$.

It is now sufficient to prove that
$a(p_i):=\lim_{y\to p_i}a(y)$ is finite and non-vanishing (by Eq.\lirho\ such a limit necessarily
exists); in fact, in this case, after fixing the
normalization $A\equiv a(p_i)$, we obtain
$$L^{(i)}_{ij}=\lim_{y\to p_i}{\sum_{k\in
I_g}\theta_k(\b+\c-y-\Delta)\omega_k(p_j)\over
E(y,p_i)}=\Lambda^{(i)}_{ij}(\c)\ .$$ Then, to conclude, it remains
to prove that $\lim_{y\to p_i}a(y)\neq 0,\infty$. Since
$L^{(i)}_{ij}$ and $\Lambda_{ij}^{(i)}$ are finite, $\lim_{y\to
p_i}a(y)= 0$ would imply that $L^{(i)}_{ij}=0$ for all $j$ and then
that Eq.\relperx\ is trivial, which is absurd.

In order to prove that $\lim_{y\to
p_i}a(y)\neq \infty$, let us choose $j\neq i$, $3\le j\le g$, in such a way that, at the point $p_j$,
$\sigma_1$ has a zero of order $k+1$ and $\lambda_i^{(c)}(z)$ has a zero
of order $k$, for some $k\ge 0$. Suppose, by absurd, that such a $j$ does not
exist. Then, by Lemma \orderzero, $\sigma_l(v)=0$, for all $l\in
I_g\setminus\{2\}$. On the other hand, such differentials also
vanish at $u$, so that $h^0(K_C\otimes\O(-u-v))=g-1$. By the
Riemann-Roch Theorem, this would imply that $h^0(\O(u+v))=1$ and
then $C$ would be hyperelliptic, counter the hypotheses.

As discussed above, the hypotheses of the theorem imply that the
greater common divisor of $(\sigma_1)$ and $(\sigma_i)$ is
$\c+\sum_{m\neq i}p_m$; in particular, if $k>0$, then $p_j$ is a
single zero for $\sigma_i$. Hence, by Eq.\relaperx, $\rho(z)$ has a
zero of order at least $k$ in $p_j$. By expanding $\rho(z)$ in the
limit $z\to p_j$, we obtain
$$\rho(z)\sim\beta\zeta^kd\zeta+o(\zeta^k)\ ,$$
with respect to some coordinates $\zeta(z)$ centered in $p_j$. Here, $\beta$ does not depend on $y$,
since, by Eq.\ydep, $\rho(z)$ depends on $y$ only through a term proportional to $\sigma_1(z)$, which is of order
$\zeta^{k+1}$.

By using Eq.\birho, in the limit $z\to p_j$ we have
$${\sum_{a\in I_g}\theta_a(u+v+\sum_{m=3}^gp_m-y-\Delta)
\omega_a(z)\over E(p_i,y)}\sim {A\beta\over a(y)}\zeta^kd\zeta+o(\zeta^k)\ .
$$
In the limit $y\to p_i$, the LHS gives $\lambda_i^{(c)}(z)$, which, by Lemma
\orderzero, has a zero of order exactly $k$ in $z=p_j$. Therefore,
$$\lim_{y\to p_i}{A\beta\over a(y)}\neq 0\ ,$$
that concludes the proof. \hfill$\square$

\vskip 6pt

A classical result known by Riemann is the relation
$$\sum_{a,b\in I_g}\theta_{ab}(e)\omega_a\omega_b=0\ ,$$
which holds for an arbitrary $e\in\Theta_s$. The connection of such a relation to the ones
considered in this paper is given by the following lemma.

\vskip 6pt

\newth\FaKraRel{Lemma}{Choose $p_1,\ldots,p_g$ satisfying conditions i), ii) or iii) of Theorem
\quadgen. Then, for all $e\in\Theta_s$, the relation
$$\sum_{a,b\in I_g}\theta_{ab}(e)\omega_a\omega_b=0\ ,$$
is equivalent to
$$\sum_{i=N+1}^MA_i(e)\sum_{j\in I_M}C^\sigma_{ij}\sigma\sigma_j=0\ ,
$$
where
$$A_i(e):=\sum_{a,b\in I_g}\theta_{ab}(e)\omega_a(p_{\1_i})\omega_b(p_{\2_i})\ ,$$
$i\in I_M$.
}

\vskip 6pt

\noindent{\bf Proof.} Two relations are equivalent if they correspond to the same vector in
$\ker\psi$, up to normalization. Since $\ker\psi$ is spanned by
$\{\sum_{k=1}^MC^\sigma_{ik}\sigma\cdot\sigma_k\}_{N<i\le M}$, then
$$\sum_{a,b\in I_g}\theta_{ab}(e)\omega_a\cdot\omega_b=
\sum_{i=N+1}^MA_i(e)\sum_{j\in I_M}C^\sigma_{ij}\sigma\cdot\sigma_j\ ,
$$ for some complex coefficients $A_i(e)$, $i\in I_M$. By applying $p\cdot p_i$, $i=N+1,\ldots,M$,
to both sides of this equation, and using $C^\sigma_{ij}=\delta_{ij}$, for $N<i,j\le M$, we conclude.\hfill$\square$

\vskip 6pt

\newth\voila{Theorem}{Choose $p_3,\ldots,p_g\in C$, $C_2\ni \c:=u+v$ and $x\in C$ as above.
Suppose $\bar A(\c)\neq\emptyset$ and fix $i\in\bar A(\c)$. If $u$
is a single zero for $K(\ \cdot\ ,p_3,\ldots,\check p_i,\ldots,
p_g)$, then the linear relation
$$\sum_{k\in I^{1i}_N}\tilde C^{\sigma(i)}_{k}(\c)\sigma\sigma_k=0\ ,$$
is equivalent to
$$\sum_{a,b\in I_g}\theta_{ab}(\c+\sum_{j\neq i}p_j-\Delta)\omega_a\omega_b=0\ .
$$
}

\vskip 6pt

\noindent{\bf Proof.} By construction, $I(\c+\sum_{j\neq i}p_j-\Delta)\in \Theta_s$. Then, use
Theorem \coefflambda\ and Lemma \FaKraRel, and note that
$$A_k(I(\c+\sum_{j\neq i}p_j-\Delta))=\Lambda_k^{(i)}(\c)\ ,
$$
$k=N+1,\ldots,M$.\hfill$\square$

\vskip 6pt

\newth\trigon{Theorem}{If $C$ is a trigonal curve, then there exist $2g-4$ pairwise distinct
 points $p_3,\ldots,p_g,u_3,\ldots,u_g\in C$
such that $K(p_3,\ldots,p_g)\neq0$ and $K(u_j,p_3,\ldots,\check
p_i,\ldots,p_g)=0$ if and only if $j\neq i$, for all $i,j\in
I_g\setminus\{1,2\}$. Furthermore, if, for each $i\in
I_g\setminus\{1,2\}$, the points $u_j$, $j\in
I_g\setminus\{1,2,i\}$, are single zeros for
$K(\cdot,p_3,\ldots,\check p_i,\ldots,p_g)$, then the following
statements hold:
\smallskip
\item{a.} For each $3\le j\le g$, there exists a unique $v_j\in C$ such that
$$I(\c_j+\sum_{k\neq
i}p_k-\Delta)\in\Theta_s\ ,$$ for all $i\neq j$, $3\le i\le j$,
where $\c_j:=u_j+v_j$, $3\le j\le g$;\smallskip
\item{b.} The
relations
$$\sum_{k\in I^{1i}_N}\tilde C^{\sigma(i)}_k(\c_j)\sigma\sigma_k=0\ ,$$
$3\le i<j\le g$, considered in Lemma \coefflambda, are linearly independent
and then generate the ideal $I_2$ of quadrics in $\PP_{g-1}$ containing the curve $C$.
\smallskip}

\vskip 6pt

\noindent {\bf Proof.} Since $C$ is trigonal, there exists a unique
(up to a fractional linear transformation) meromorphic function $f$
with three poles. Hence, for each $p\in C$, $f^{-1}(f(p))$ consists of
three (possibly coincident) points; note that, trivially, the sum of
such three points (counting multiplicity) corresponds to the unique
effective divisor of degree three which is special and containing $p$
in its support.

Fix $x_4,\ldots,x_g\in C$, and consider the function
$$F_{x_4,\ldots,x_g}(p):=\prod_{x\in f^{-1}(f(p)))}K(x,x_4,\ldots,x_g)\
,\qquad p\in C\ .
$$ Denote by $[K]_{x_4,\ldots,x_g}\subseteq C$ and $[F]_{x_4,\ldots,x_g}$ the sets of zeros
of $K(\cdot,x_4,\ldots,x_g)$ and $F_{x_4,\ldots,x_g}$, respectively.
Then, one of the following alternatives holds: if
$K(\cdot,x_4,\ldots,x_g)$ is not identically vanishing, then both
$[K]_{x_4,\ldots,x_g}$ and
$$[F]_{x_4,\ldots,x_g}=\bigcup_{x\in [K]_{x_4,\ldots,x_g}}f^{-1}(f(x))\
,$$ are finite sets; otherwise, both $[K]_{x_4,\ldots,x_g}$ and
$[F]_{x_4,\ldots,x_g}$ coincide with $C$.

For each $n$, $1\le n\le g-2$, let
$N^{(n)}_{x_{n+3},\ldots,x_{g}}\subseteq C^n$ denote the set of
$n$-tuples $(p_3,\ldots,p_{n+2})$ such that
$$F^{(n)}_{x_{n+3},\ldots,x_{g}}(p_3,\ldots,p_{n+2}):=\prod_{i=1}^n
F_{p_3,\ldots,\check p_i,\ldots,p_{n+2},x_{n+3},\ldots,x_g}(p_i)\ ,
$$
is not zero. Note that $F^{(1)}\equiv F$ and $N^{(1)}=C\setminus
[F]$.

Now, assume that, for some $m$, $1\le m< g-2$, the set $N^{(n)}$ is
dense in $C^n$ for all $n\le m$. The set
$[F^{(m+1)}]_{x_{m+4},\ldots,x_{g}}$ of zeros of
$$F^{(m+1)}_{x_{m+4},\ldots,x_{g}}(p_3,\ldots,p_{m+2},p)=
F_{p_3,\ldots,p_{m+2},x_{m+4},\ldots,x_g}(p)\prod_{i=1}^m
F_{p_3,\ldots,\check p_i,\ldots,p_{m+2},p,x_{m+4},\ldots,x_g}(p_i)\
,
$$
as a function of $p$, is given by
$$[F^{(m+1)}]_{x_{m+4},\ldots,x_{g}}=\bigcup_{i=1}^m\bigl(\bigcup_{x\in f^{-1}(f(p_i))}
[K]_{p_3,\ldots,\check
p_i,\ldots,p_{m+2},x,x_{m+4},\ldots,x_g}\bigr)\cup
[F]_{p_3,\ldots,p_{m+2},x_{m+4},\ldots,x_g}\ .
$$ If $(p_3,\ldots,p_{m+2})\in N^{(m)}$, then the functions
$$K(\cdot,p_3,\ldots,p_{m+2},x_{m+4},\ldots,x_g)\ ,$$
and
$$K(\cdot,p_3,\ldots,\check
p_i,\ldots,p_{m+2},x,x_{m+4},\ldots,x_g)\ ,$$ for each
$i=1,\ldots,m$, and $x\in f^{-1}(f(p_i))$,
 vanish identically on $C$
(for example, $x_{m+3}$ is not a zero). Hence,
$[F^{(m+1)}]_{x_{m+4},\ldots,x_{g}}\subseteq C$ is a finite set and,
therefore, $N^{(m+1)}_{x_{m+4},\ldots,x_{g}}$ is dense in $C^{m+1}$.
We proved that if $K(\cdot,x_4,\ldots,x_g)$ does not identically
vanish for some $x_4,\ldots,x_g\in C$, then
$N^{n}_{x_{n+3},\ldots,x_g}$ is dense in $C^n$ for all $n$, $1\le
n\le g-2$. It follows that $N^{g-2}$, which does not depend on
$x_4,\ldots,x_g$, is dense in $C^{g-2}$. Also note that the subset
of $C^{g-2}$ for which
$$\bigcup_{i=3}^gf^{-1}(f(p_i))\ ,$$ consists of
pairwise distinct points is dense $C_{g-2}$. Hence, its intersection
with $N^{(g-2)}$ is not empty. Let us choose $(p_3,\ldots,p_g)$ in
such an intersection and fix $u_i\in f^{-1}(f(p_i))$, $u_i\neq p_i$,
for all $i\in I_g\setminus \{1,2\}$. Then, the points $p_3,\ldots,p_g,u_3,\ldots,u_g$ are pairwise
distinct and satisfy the condition
$$K(u_i,p_3,\ldots,\check p_j,\ldots,p_g)=0\ \Leftrightarrow\ i\neq j\ ,$$
for all $i,j\in I_g\setminus\{1,2\}$. Furthermore, if $u_i$, $i\in I_g\setminus\{1,2\}$, is a single zero of $K(\cdot,p_3,\ldots,\check
p_j,\ldots,p_g)$, for all $j\in I_g\setminus\{1,2,i\}$, then there exists a unique point $v_{ij}$
such that $v_{ij}+u_i+\sum_{k\neq j}p_k$ is special. Such a point satisfies necessarily $f^{-1}(f(p_i))=\{p_i,u_i,v_{ij}\}$,
so that it is independent of $j$, and the statement {\it a.} follows.

Finally, note that
$A(\c_j)=p_j$ and $\bar A(\c_j)=\{p_i\mid 3\le i\le g, i\neq j\}$.
Hence, by Theorem \coefflambda, for each $k$, $N<k\le M$ the
coefficients $C^{\sigma(\1_k)}_{l}(\c_{\2_k})$, $l\in I_M$, are given by
$$
C^{\sigma(\1_k)}_{l}(\c_{\2_k})=\Lambda^{(\1_k)}_{k}(\c_{\2_k})C^{\sigma}_{kl}
\ ,$$
where $\Lambda^{(\1_k)}_{k}(\c_{\2_k})\neq 0$.
Linear independence of the $\tilde C^{\sigma(i)}(\c_j)$'s, $3\le i<j\le g$, follows by linear independence of the
$C^{\sigma}_{k}$'s.\hfill$\square$

\vskip 6pt

\subsec{The case of genus $4$}

Consider the case of a non-hyperelliptic curve $C$ of genus $4$. The identity \Kkappa\ reduces to
$$K(p_3,p_4):=-c_{4,2}{\kdue[v]\over\kuno[\sigma]^{5}E(p_3,p_4)^2
\sigma(p_3)\sigma(p_4)}\ ,
$$
where $c_{4,2}=1008$, and can be used to express Eq.\corol\ in terms of the function $K$. For $g=4$ Eq.\corol\
reduces to a unique relation
$$\sum_{i=1}^{10}C^\sigma_i\sigma\sigma_i=0\ .$$
It can be derived from the identity
$${\det_{i,j\in I_{10}} \sigma\sigma_i(x_j)\over \det_{i,j\in I_9} v_i(x_j)}=0\ ,$$
by expanding the determinant at the numerator with respect to the column corresponding to $x_{10}\equiv
z$. One obtains
$$\sum_{i\in I_{10}}(-)^i{\det_{^{j\in I_{10}\setminus\{i\}}_{k\in I_9}}\sigma\sigma_j(x_k)\over
\det_{j,k\in I_9} v_j(x_k)}\sigma\sigma_i(z)=0\ ,$$
where the ratios of determinants do not depend on $x_1,\ldots,x_9$ and correspond to
$${\det_{^{j\in I_{10}\setminus\{i\}}_{k\in I_9}}\sigma\sigma_j(x_k)\over
\det_{j,k\in I_9}
v_j(x_k)}={\kappa[\sigma\sigma_1,\ldots,\check{\sigma\sigma_i},\ldots,\sigma\sigma_{10}]\over
\kappa[v]}\ .$$
Now, note that for $\1_i=\2_i$,
$\kappa[\sigma\sigma_1,\ldots,\check{\sigma\sigma_i},\ldots,\sigma\sigma_{10}]=0$. This can be
checked by observing that all the elements in $\{\sigma\sigma_j\}_{j\in I_{10}\setminus\{i\}}$ vanish
at $p_i$, so that it cannot be a basis of $H^0(K^2_C)$. Hence, we can restrict the summation over all
the $i\in I_{10}$ with $\1_i\neq \2_i$. By a re-labeling of the points $p_1,\ldots,p_4$, the
relation between $\kappa[v]$ and $K$ at genus four is
$$K(p_{\1_i},p_{\2_i})=(-)^{i+1}c_{4,2}{\kappa[\sigma\sigma_1,\ldots,\check{\sigma\sigma_i},\ldots,\sigma\sigma_{10}]\over
\kappa[\sigma]^{5}E(p_{\1_i},p_{\2_i})^2\sigma(p_{\1_i})\sigma(p_{\2_i})}\
,$$
for all $i$, $5\le i\le 10$. Hence,
\eqn\coeffKFour{C_i^\sigma={K(p_{\1_i},p_{\2_i})E(p_{\1_i},p_{\2_i})^2\sigma(p_{\1_i})\sigma(p_{\2_i})\over
K(p_3,p_4)E(p_3,p_4)^2\sigma(p_3)\sigma(p_4)}={k(p_{\1_i},p_{\2_i})\over k(p_3,p_4)}\ ,}
$5\le i\le 10$, with $k$ defined in Eq.\kpiccola, whereas $C_i^\sigma=0$ for $i\le 4$. Since
$\sigma\sigma_i=\sum_{j=1}^{10}X_{ji}\omega\omega_j$, it follows that
$$C^\omega_i=\sum_{j=5}^{10}X^\omega_{ij}C^\sigma_j\ ,$$ $i\in
I_{10}$, and we obtain
\eqn\ComegaKFour{C^\omega_i=\chi_i^{-1}\sum_{k,l=1}^4{k(p_k,p_l)\over
k(p_3,p_4)}{\theta_{\Delta,\1_i}(\a_k)
\theta_{\Delta,\2_i}(\a_l)\over
\sum_{m,n}\theta_{\Delta,m}(\a_k)
\theta_{\Delta,n}(\a_l)\omega_m(p_k)\omega_n(p_l)}\ ,}
$i\in I_{10}$.
Note that $\hat C^\omega_i:=k(p_3,p_4)C^\omega_i$ is symmetric
under any permutation of $p_1,\ldots,p_4$. On the other hand, Corollary \threla\ shows
that $C_i^\omega$, and therefore also $\hat C^\omega_i$, are independent of $p_1,p_2$. We conclude that
$\hat C^\omega_i$, whose explicit form is
$$\eqalignno{\hat C^\omega_i&=-{\chi_i^{-1}\over S(\a)^2 \prod_{1}^{7}\sigma(x_i)}\sum_{_{k\neq
l}^{k,l=1}}^4{\theta_{\Delta,\1_i}(\a_k)
\theta_{\Delta,\2_i}(\a_l)\over
\deltadiv\bigl(p_k+p_l+\sum_{1}^{7}x_i\bigr)\sigma(p_k)\sigma(p_l)\prod_{i\neq
k,l}\bigl(E(p_k,p_i)E(p_l,p_i)\bigr)
}\cr\cr &\cdot \sum_{s\in\perm_{7}}
{S\bigl(\sum_{i=1}^4x_{s_i}\bigr)
S\bigl(\sum_{i=4}^{7}x_{s_i}\bigr) \over
E(x_{s_g},p_k)E(x_{s_g},p_l)}\prod_{i=1}^{3}{
S(x_{s_i}+x_{s_{i+4}}+p_k+p_l) \over
\prod_{_{j\neq
i}^{j=1}}^{\ss 3}E(x_{s_i}, x_{s_{j+4}})} \ ,}
$$
does not depend on $p_1,\ldots,p_4$, for all $i\in I_{10}$.

\vskip 6pt

Note that, at genus $4$, the equivalent relations
$$\sum_{i\in I_M}C^\omega_i\omega\omega_i=0 ,$$
and
$$\sum_{i\in I_M}\hat C^\omega_i\omega\omega_i=0\ ,$$ must be proportional to Eq.\RiemRel, with $e$ one of the
two points in $\Theta_s$; in other words, $C^\omega_i$ and $\hat C^\omega_i$ must be proportional
to $\chi_i^{-1}\theta_{\1_i\2_i}(e)$. In the following proposition, such a proportionality is precisely derived.

\vskip 6pt

\newth\FaKraFour{Proposition}{Let $C$ be non-hyperelliptic of genus $4$ and fix $(p_1,\ldots,p_4)\in
C^4\setminus\B$. The coefficients $C^\sigma_{i}$, $i\in I_M$, correspond to
\eqn\coeffRieFour{C^\sigma_{i}={\sum_{a,b\in I_g}\theta_{ab}(e)\omega_a(p_{\1_i})\omega_b(p_{\2_i})\over
\sum_{a,b\in I_g}\theta_{ab}(e)\omega_a(p_3)\omega_b(p_4)}\ ,
}
for all $i\in I_M$,
and
\eqn\thetaKFour{k(p,q)\equiv K(p,q)E(p,q)^2\sigma(p)\sigma(q)=
A\sum_{a,b\in I_g}\theta_{ab}(e)\omega_a(p)\omega_b(q)\ ,} for all $p,q\in C$, where $e\in\Theta_s$
and $A$ is a complex constant (depending on the moduli). Furthermore,
\eqn\ComegathetaFour{C^\omega_i=\chi_i^{-1}{\theta_{\1_i\2_i}(e)\over\sum_{a,b\in I_g}\theta_{ab}(e)\omega_a(p_3)\omega_b(p_4)} \ ,
}
for all $i\in I_M$.
}

\vskip 6pt

\noindent {\bf Proof.} Lemma \FaKraRel\ gives
$$\sum_{a,b\in I_g}\theta_{ab}(e)\omega_a\cdot\omega_b=
2\sum_{a,b\in I_g}\theta_{ab}(e)\omega_a(p_3)\omega_b(p_4)\sum_{i\in
I_M}C^\sigma_i\sigma\cdot\sigma_i\ ,$$
and, by applying $\chi_i^{-1}(p\cdot p)_i$, $i\in I_M$, to both sides, one obtains Eq.\coeffRieFour\
(note that $\chi_i\neq 1$ if and only if $p_{\1_i}=p_{\2_i}$, and in this case both sides of
\coeffRieFour\ vanish). By comparing Eq.\coeffRieFour\ and Eq.\coeffKFour, we obtain Eq.\thetaKFour.
Finally, by using Eq.\ComegaKFour\ and Eq.\thetaKFour, we have
$$C^\omega_i=\chi_i^{-1}\sum_{c,d\in I_g}\biggl({\theta_{cd}(e)\over
\sum_{a,b\in I_g}\theta_{ab}(e)\omega_a(p_3)\omega_b(p_4)}
\sum_{k=1}^4{\theta_{\Delta,\1_i}(\a_k)\omega_c(p_k)
\over\sum_{m}\theta_{\Delta,m}(\a_k)\omega_m(p_k)}\sum_{l=1}^4{\theta_{\Delta,\2_i}(\a_l)\omega_d(p_l)\over
\sum_n\theta_{\Delta,n}(\a_l)\omega_n(p_l)}\biggr)
$$
and, by Corollary \surprise, we obtain
Eq.\ComegathetaFour.\hfill$\square$

\vskip 6pt

\subsec{Relations among holomorphic cubic differentials}

According to Petri's Theorem, in the most general case the ideal of
a canonical curve $C$ is generated by its ideals of quadrics
together with the ideal of cubics. As discussed in the introduction
of this section, such cubics correspond to linear relations among
holomorphic $3$-differentials on $C$; a generalization of the
previous construction is necessary in order to explicitly determine
such relations.

Fix $p_1,\ldots,p_g\in C$ satisfying the conditions {\it i}), {\it
ii}) and {\it iii}) of Proposition \thtrebase\ with respect to some
fixed $i$, $3\le i\le g$, and let $\{\varphi_j\}_{j\in
I_{N_3-1}}\cup \{\varphi_{i+5g-8}\}$ be the corresponding basis of
$H^0(K_C^3)$. The kernel of the canonical epimorphism from ${\rm
Sym^3}H^0(K_C)$ onto $H^0(K_C^3)$ has dimension
$(g-3)(g^2+6g-10)/6$, and each element corresponds to a linear
combination of the following relations
\eqn\trereli{\sigma_j\sigma_k\sigma_l=\sum_{m\in
I_{N_3-1}}B_{jkl,m}\varphi_m+B_{jkl,i+5g-8}\sigma_2\sigma_i^2\ ,}
$3\le j,k,l\le g$, $j\neq k$, and
\eqn\trerelii{\sigma_2\sigma_j^2=\sum_{m\in
I_{N_3-1}}B_{2jj,m}\varphi_m+ B_{2jj,i+5g-8}\sigma_2\sigma_i^2\ ,}
$3\le j\le g$, $j\neq i$, where $B_{jkl,m},B_{2jj,m}\in\CC$, are
suitable coefficients. On the other hand, a trivial computation
shows that the relations \trereli\ are generated by \trerelii\ and
by the relations among holomorphic quadratic differentials,
\eqn\duerel{\sum_{j=1}^MC^\sigma_{kj}\sigma\sigma_j=0\ ,}
$k=N+1,\ldots,M$. Therefore, relations among holomorphic
$3$-differentials, modulo relations among holomorphic quadratic
differentials, provide at most $g-3$ independent conditions on
products of elements of $H^0(K_C)$.

\vskip 6pt

The relations \trerelii\ can be restated in terms of an arbitrary
basis $\{\eta_j\}_{j\in I_g}$ of $H^0(K_C)$. Let $Y^\eta$ be the
automorphism of $\CC^{M_3}$, determined by
\eqn\YYYY{Y^\eta_{kj}:=\chi_k^{-1}([\eta]^{-1}[\eta]^{-1}[\eta]^{-1})_{jk}\
,} $j,k\in I_{M_3}$, so that
$$\varphi_j=\sum_{k=1}^{M_3}Y^\eta_{kj}\eta\eta\eta_k\ ,$$ $j\in I_{M_3}$.
Consider the following determinants of $d$-dimensional submatrices
of $Y^\eta$
$$
\Ymin{\eta}{\ss j_1 \ldots j_d\cr\ss i_1 \ldots i_d}:=
\det\left(\matrix{Y^\eta_{i_1j_1} & \ldots & Y^\eta_{i_1j_d}\cr
\vdots & \ddots & \vdots \cr Y^\eta_{i_dj_1} & \ldots &
Y^\eta_{i_dj_d}}\right)\ ,$$ $i_1,\ldots,i_d,j_1,\ldots,j_d\in
I_{M_3}$, $d\in I_{M_3}$.

\vskip 6pt

\newth\thcoroltre{Proposition}{\eqn\coroltre{\sum_{j=1}^{M_3}D^\eta_{kj}\eta\eta\eta_j=0\ ,}
$N_3\leq k\le N_3+g-3$, $k\neq i$ where
\eqn\leCtre{D^\eta_{kj}:=\sum_{k_1,\ldots,k_{N_3}=1}^{M_3}
\Ymin{\eta}{\ss 1 \hfill \ldots \hfill (N_3-1)\, i \, k\cr \ss
k_1\hfill \ldots \hfill\  k_{N_3}\, j\hfill}\;
\R_\varphi[\eta\eta\eta_{k_1},\ldots,\eta\eta\eta_{k_{N_3}}]\ ,}
$j\in I_{M_3}$, are $g-3$ independent linear relations among
holomorphic $3$-differentials.}

\vskip 6pt

\noindent{\bf Proof.} Without loss of generality, we can assume
$i=N_3$; such an assumption can always be satisfied after a
re-ordering of the points $p_3,\ldots,p_g$. Fix $N_3+1$ arbitrary
points $x_1,\ldots,x_{N_3},x_{N_3+1}\equiv z\in C$ and consider the
singular matrix $[\varphi_l(x_m)]_{^{l\in I}_{m\in I_{N_3+1}}}$ with
$I:=I_{N_3}\cup \{k\}$, with $N_3< k\le N_3+g-3$. By expressing the
determinant with respect to the column $(\varphi_l(z))_{l\in I}$,
the identity $\det\varphi_l(x_m)=0$, $l\in I$, $m\in I_{N_3+1}$,
yields
$$\sum_{m=1}^{M_3}\biggl[\sum_{l=1}^{N_3}(-)^{l+1}
\R_\varphi[\varphi_1,\ldots,\check{\varphi}_l,\ldots,
\varphi_{N_3},\varphi_k]Y^\eta_{ml}-Y^\eta_{mk}\biggr]\eta\eta\eta_m=0\
.$$ The proposition follows by combinatorial identities analogous to
the proof of Theorem \thcorol.\hfill$\square$

\vskip 6pt

Whereas for $g=4$ the relations \coroltre\ are independent of the
relation among holomorphic quadratic differentials, for $g\ge 5$, \coroltre\ are generated by
\duerel\ in all but some particular curves. Set
$\tilde\psi_{\1_i\2_i,\1_j\2_j}:=\tilde\psi_{ij}$ and
$C^\sigma_{\1_i\2_i,\1_j\2_j}:=C^\sigma_{ij}$, $N+1\le i\le M$,
$j\in I_M$. Consider the $3$-differentials
$\sigma_i\sigma_j\sigma_k$ with $3\le i<j<k\le g$ ($g\ge 5$). By
Eq.\duerel\ and by $C^\sigma_{ij}=\tilde\psi_{ij}-\delta_{ij}$,
$N+1\le i\le M$, $j\in I_M$,
$$\sigma_i\sigma_j\sigma_k=\sum_{m=1}^2\sum_{n=3}^g\tilde\psi_{ij,mn}\sigma_m\sigma_n\sigma_k+
\tilde\psi_{ij,12}\sigma_1\sigma_2\sigma_j \ ,$$
 so that
 $$\eqalignno{&\sum_{m,p=1}^2\sum_{q=3}^g\bigl(\sum_{^{n=3}_{n\neq j}}^g\tilde\psi_{ik,mn}
 \tilde\psi_{nj,pq}\bigr)\sigma_m\sigma_p\sigma_q+
\tilde\psi_{ik,12}\sigma_1\sigma_2\sigma_k+\sum_{m=1}^2\tilde\psi_{ik,mj}\sigma_m\sigma_j^2\cr\cr
=&\sum_{m,p=1}^2\sum_{q=3}^g\bigl(\sum_{^{n=3}_{n\neq
i}}^g\tilde\psi_{jk,mn}
 \tilde\psi_{ni,pq}\bigr)\sigma_m\sigma_p\sigma_q+
\tilde\psi_{jk,12}\sigma_1\sigma_2\sigma_k+\sum_{m=1}^2\tilde\psi_{jk,mi}\sigma_m\sigma_i^2\
.}$$ The above equation yields
$$\eqalignno{C^\sigma_{ik,2j}\sigma_2\sigma_j^2=&
\sum_{m,p=1}^2\sum_{q=3}^g \bigl(\sum_{^{n=3}_{n\neq
i}}^gC^\sigma_{jk,mn}
 C^\sigma_{ni,pq}-\sum_{^{n=3}_{n\neq j}}^gC^\sigma_{ik,mn}
 C^\sigma_{nj,pq}\bigr)\sigma_m\sigma_p\sigma_q\cr\cr
+&C^\sigma_{jk,12}\sigma_1\sigma_2\sigma_k-C^\sigma_{ik,12}\sigma_1\sigma_2\sigma_k+C^\sigma_{jk,1i}\sigma_1\sigma_i^2
-C^\sigma_{ik,1j}\sigma_1\sigma_j^2+C^\sigma_{jk,2i}\sigma_2\sigma_i^2\
.}
$$
If $C^\sigma_{ik,2j}\neq 0$ for some $k$, the above identity shows
that the relation \trerelii\ is generated by Eqs.\duerel. On the
other hand, it can be proved \ottimo\ that if $C^\sigma_{ik,2j}=0$ for all
$3\le k\le g$, $k\neq i,j$, the relation \trerelii\ is independent
of the relations among holomorphic quadratic differentials. This case occurs if
and only if the curve $C$ is trigonal or a smooth quintic.

\vskip 6pt

\newth\thletrerel{Proposition}{Fix $g$ points $p_1,\ldots,p_g\in C$ satisfying
the conditions of theorem \thtrebase. The coefficients
$Y^\omega_{ij}$, defined in Eq.\YYYY\ with $\eta\equiv\omega$, are
given by \eqn\leYtre{Y^\omega_{ij}=
{(1+\delta_{\1_j\2_j}+\delta_{\2_j\3_j})(1+\delta_{\1_j\3_j})\over
6\prod_{{\frak m}\in\{\1,\2,\3\}}\sum_l \theta_{\Delta,l}(\a_{{\frak
m}_i})\omega_l(p_{{\frak
m}_i})}\sum_{s\in\perm_3}\Bigl(\smprod_{{\frak m}\in\{\1,\2,\3\}}
\theta_{\Delta,s({\frak m})_j}(\a_{{\frak m}_i})\Bigr) \ ,} $i,j\in
I_M$.}

\vskip 6pt

\noindent {\bf Proof.} The proposition follows immediately by the
definition \YYYY\ and by Eq.\leomega.\hfill$\square$

\vskip 6pt

\newsec{The canonical basis of $H^0(K_C)$ and $\kuno[\omega]$ in terms of theta functions}\seclab\thetkuno

We saw that the bases for holomorphic differentials we introduced
satisfy several properties. Here we first show that the use of such
bases leads to a straightforward derivation of the Fay's trisecant
identity. We then consider a basic problem in the study of Riemann
surfaces. This arises, for example, in investigating the Schottky
problem or in constructing modular forms, where basic quantities, in
spite of being constants, e.g. the Mumford form, have an expression
that needs the use of points on $C$. We introduce a new general
strategy which is based on the idea of identifying the divisors with
the ones defining spin structures. In doing this one has to consider
several intermediate problems, such as expressing the determinants
of the canonical holomorphic abelian differentials in terms of
theta functions only. Furthermore, another problem concerning
$\det\omega_i(p_j)$ is to consider the $g$-points $p_1,\ldots,p_g$
as defining spin structures which are associate to divisors of
degree $g-1$. Such a question is strictly related to the problem of
expressing $\det\omega_i(p_j)$ without the use of the
$g/2$-differential $\sigma$ and of the constant $\kuno[\omega]$. We
will see that there exists an elegant and natural solution leading
to the explicit expression of basic quantities in terms of divisors
defining spin structures. We also express
the abelian holomorphic differentials in terms of theta functions
only.
This also implies the expression for products of the basic constants
$\kuno_{\nu_k}[\omega]$, corresponding to the main building block of the
Mumford form, in terms of theta functions with spin structures whose
arguments involve the difference of points belonging to the divisors
of such spin structures.

\vskip 6pt

\subsec{Determinants and Fay's identity}

In this section, we will use the bases introduced in section
\primecostr\ to derive
 a combinatorial proof of the
Fay's trisecant identity.

\vskip 6pt

\newth\Faytris{Theorem}{The following are equivalent
\smallskip
\item{\it a)} Proposition \thdettheta\ holds;
\smallskip
\item{\it b)} The Fay's trisecant identity \jfayy\
\eqn\Fay{{\theta(w+\sum_{1}^m(x_i-y_i))\prod_{i<j}E(x_i,x_j)E(y_i,y_j)\over
\theta(w)\prod_{i,j}E(x_i,y_j)}=(-)^{{m(m-1)\over 2}}\det\nolimits_{ij}{\theta(w+x_i-y_j)\over
\theta(w)E(x_i,y_j)}\ ,}
$m\ge 2$, holds for all $x_1,\ldots,x_m,y_1,\ldots,y_m\in C$, $w\in J_0(C)$.
}

\vskip 6pt

\noindent {\bf Proof.} ($a \Rightarrow b$) Fix $x_1,\ldots,x_m,y_1,\ldots,y_m\in C$ and
$w\in J_0(C)$, with $\theta(w)\neq 0$. Choose
$y_1,\ldots,y_m$ distinct, otherwise the identity
is trivial. Set $p_i:= y_i$, $i\in I_m$, and fix $n\in \NN_+$,
with $d:=N_n-m\ge g$, and $p_{m+1},\ldots,p_{N_n}\in C$, in such a way that
$$I(\smsum_1^{N_n}p_i-(2n-1)\Delta)=w\ .$$
By Jacobi Inversion Theorem, such a choice is always possible.
Note that the set of divisors
$p_{m+1}+\ldots+p_{N_n}$, such that $p_i=p_j$ for some $i\neq j\in
I_{N_n}$, is the set of points of a subvariety in the space of
positive divisors of degree $d$. Then the image of such a variety
under the Jacobi map, which is analytic, corresponds to a proper subvariety
$W$ of $J_0(C)$. Hence, the conditions
$\theta(w)\neq 0$ and $$w-I(\smsum_{1}^my_i+(2n-1)\Delta)\in
J_0(C)\setminus W\ ,$$ are satisfied for $w$ a dense subset in $J_0(C)$. It
is therefore sufficient to prove Eq.\Fay\ on such a subset
and the theorem follows by continuity arguments.

Let us then choose the points $p_{m+1},\ldots,p_{N_n}$ to be
pairwise distinct and distinct from $y_1,\ldots,y_m$ and fix a basis
$\{\phi^n_i\}_{i\in I_{N_n}}$ of $H^0(K_C^n)$.
Since $p_1,\ldots,p_{N_n}$ are pairwise distinct and
$$\theta_\Delta(\smsum_1^{N_n}p_i)=\theta(w)\neq 0\ ,$$
it follows by Eq.\detthetaii\
that $\det\phi^n_i(p_j)\neq 0$.
Therefore, by Proposition \thnewbasis, one can define the basis
$\{\gamma^n_i\}_{i\in I_{N_n}}$ of $H^0(K_C^n)$ with the property
$\gamma^n_i(p_j)=\delta_{ij}$, $i,j\in I_{N_n}$. On the other hand, note that
$$\det\gamma^n(x_1,\ldots,x_m,p_{m+1},\ldots,p_{N_n})= \det_{ij\in
I_m}\gamma^n_i(x_j)\ ,$$ can be expressed either
by means of Eq.\glindiff\
$$
\prod_{i=1}^{m}\sigma(x_i,y_i)^{2n-1}\prod_{j=m+1}^{N_n}{E(x_i,p_j)\over
E(y_i,p_j)} {\prod_{^{i,j=1}}^{m}E(x_i,y_j)\over\prod_{^{i,j=1}_{i\neq
j}}^{m}E(y_i,y_j)}\det_{ij}{\theta(w+x_i-y_j)\over
\theta(w)E(x_i,y_j)}\ ,
$$
or by means of \detthetaii\ and \knphi\
$$
\prod_{i=1}^{m}\sigma(x_i,y_i)^{2n-1}\prod_{j=m+1}^{N_n}{E(x_i,p_j)\over
E(y_i,p_j)} {\theta\bigl(w+\sum_{1}^m(x_i-y_i)\bigr)\prod_{i<j}^{m}E(x_i,x_j)\over
\theta(w)\prod_{^{i,j=1}_{i<j}}^{m}E(y_i,y_j)}\ .$$
Eq.\Fay\ then follows by observing that \eqn\perdopo{\prod_{i,j=1\atop i\neq
j}^{m}E(y_i,y_j)=(-)^{m(m-1)/2}\prod_{i,j=1\atop i<
j}^{m}E(y_i,y_j)^2\ .}

\vskip 3pt

\noindent($b\Rightarrow a$) Fix $p_1,\ldots,p_{N_n}\in C$, $n\geq2$, in such a way that the
hypothesis of Proposition \thnewbasis\ is satisfied. Let $\{\gamma^n_i\}_{i\in I_{N_n}}$ be the
corresponding basis of $H^0(K_C^n)$ satisfying \centralle.
$\det\gamma^n_i(z_j)$ can be evaluated, for arbitrary $z_1,\ldots,z_{N_n}\in C$, by expressing $\gamma^n_i(z_j)$ by means of \glindiff.
In particular, by using \Fay\ with $m=N_n$,
$x_i=z_i$, $y_i=p_i$, $i\in I_{N_n}$, and $w=I(\sum_1^{N_n}p_i-(2n-1)\Delta)$, after a
computation analogous to the previous one, \detthetaii\ follows, with $\kenne[\gamma^n]$ given by
Eq.\knphi. Therefore, \detthetaii\ holds for an arbitrary basis $\{\phi_i^n\}_{i\in I_{N_n}}$ of $H^0(K_C^n)$, with
$\kenne[\phi^n]=\kenne[\gamma^n]\det\phi^n_i(p_j)$. The same result holds for \dettheta\
by using \Fay\ with $w=I(\sum_1^gp_i-y-\Delta)$.
\hfill$\square$

\vskip 6pt

\subsec{$\{\omega_i\}_{i\in I_g}$ and $\kuno_\nu[\omega]$: from
divisors on $C^g$ to spin structures}\subseclab\secdivspin

Integrating \ildet\ along the $\alpha$-cycles of $C$ leads
to expressions of the minors of $\omega_i(p_j)$ in terms of
theta functions. In particular, denoting by $\hat\omega_k(p_1,\ldots,\check
p_i,\ldots,p_g)$ the
cofactor of $\omega_k(z)$ in $\det\omega(z,p_1,\ldots,\check
p_i,\ldots,p_g)$, we have

\vskip 6pt

\newth\minors{Proposition}{
\eqn\cofactor{\hat\omega_l(p_1,\ldots,\check
p_i,\ldots,p_g)=\oint_{\alpha_{l,z}}\!\!\det\omega(z,p_1,\ldots,\check
p_i,\ldots,p_g)=\kuno[\omega]\theta_{\Delta,l}(\a_i)\prod_{{j,k\neq i \atop j<k}}E(p_j,p_k)\prod_{j\neq
i}\sigma(p_j)\ , } $l\in I_g$, where $\oint_{\alpha_{l,z}}$ denotes
integration in $z$ along $\alpha_l$. Furthermore,
\eqn\proprioomega{\omega_l(z)=(-)^{l+1}\oint_{\alpha_{1,p_1}}\cdots
\oint_{\alpha_{l-1,p_{l-1}}}\oint_{\alpha_{l+1,p_l}}\cdots\oint_{\alpha_{g-1,p_{g-2}}}
\hat\omega_g(z,p_1,\ldots,p_{g-2})\ ,} $l\in I_g$.}

\vskip 6pt

For $g=2$ the
divisors of $\omega_1$ and $\omega_2$ coincide with the ones
of $\theta_{\Delta,2}$ and $\theta_{\Delta,1}$, respectively.
In particular
\eqn\esplicittto{\omega_1(z)=\kuno[\omega]\theta_{\Delta,2}(z)\sigma(z)
\ ,\qquad\qquad
\omega_2(z)=-\kuno[\omega]\theta_{\Delta,1}(z)\sigma(z) \ .}
This is a particular case of more general relations considered in the following theorem.

\vskip 6pt

\newth\beh{Theorem}{Fix $(p_1,\ldots,p_g)\in C^g$ and set $\a_i$ as in Definition \defAB. Then
$$\det
\theta_{\Delta,i}\big(\a_j\big)\neq0\ ,$$
if and only
if $(p_1,\ldots,p_g)\notin\A$. In this case
\eqn\oomegai{\omega_i(z)=
\sum_{j=1}^g\xi_{ij}^{-1}{\theta_\Delta(\a_j+z-y_j)\theta_\Delta(\a-z)\over
\theta_\Delta(\a-y_j)}{E(p_j,y_j)\over E(p_j,z) E(y_j,z)}\
,} $i\in I_g$, for all $y_i,z\in C$, $i\in I_g$, where
$\xi_{ij}:=\theta_{\Delta,j}(\a_i)$.
Furthermore, if $\alpha_i$, $i\in I_g$, are elements of $\Theta$ such that
$$
\det \theta_i(\alpha_j)\neq0 \ ,
$$
and $\nu_i:=(\nu_i',\nu_i'')$, $i\in I_g$, are odd theta
characteristics satisfying $$\det \theta_i[\nu_j](0)\neq0\ ,$$ then
\eqn\oomegaiuuu{\eqalign{\omega_i(z)&={E(p,q)\over E(p,z) E(q,z)}
\sum_{j=1}^g\lambda_{ij}^{-1}{\theta(\alpha_j+z-q)\theta(\alpha_j+p-z)\over\theta(\alpha_j+p-q)}
\cr\cr &= {E(p,q)\over E(p,z) E(q,z)}
\sum_{j=1}^g\mu_{ij}^{-1}{\theta[\nu_j](z-q)\theta[\nu_j](p-z)\over\theta[\nu_j](p-q)}\
,}} $i\in I_g$, for all $p,q,z\in C$, where
$\lambda_{ij}:=\theta_j(\alpha_i)$ and
$\mu_{ij}:=\theta_j[\nu_i](0)$. The three conditions, and in
particular $\det \theta_i[\nu_j](0)\neq0$, can always be satisfied
for any $C$.}

\vskip 6pt

\noindent{\bf Proof.} Fix $i\in I_g$ and consider the limit $y\to
p_i$ of \essesimm\
\eqn\giragira{\sum_j\theta_{\Delta,j}(\a_i)\omega_j(p_i)=S(\a)
\sigma(p_i)\prod_{j\neq i}E(p_i,p_j)\ .} Independence of
$\xi_{ij}:=\theta_{\Delta,j}(\a_i)$ on $p_i$ implies that the above
relation also holds whether $p_i$ is replaced by an arbitrary point
$q\in C$, and therefore replacing also $\a$ by $\a-p_i+q$. In
particular, since $\prod_{j\neq i}E(p_i,p_j)$ vanishes at $p_i=p_j$,
$j\in I_g$, $j\neq i$, it follows that $\sum_k\xi_{ik}\omega_k(p_j)$
is a diagonal matrix, so that
\eqn\detdd{\det\xi_{ij}\det\omega_i(p_j)=S(\a)^g
\prod_1^g\sigma(p_i)\prod_{_{i\neq j}^{i,j=1}}^gE(p_i,p_j)\ .}
Expressing $\det\omega_i(p_j)$ by Eq.\dettheta, and using \perdopo\
we obtain
\eqn\rigira{\det\xi_{ij}=(-)^{g(g-1)/2}\kappa[\omega]^{-1}S(\a)^{g-1}\prod_{i<j}E(p_i,p_j)\
,} implying that $\det\xi_{ij}\neq 0$ if and only if $S(\a)\neq 0$
and $p_1,\ldots,p_g$ are pairwise distinct. Now, \giragira\ with
$z\equiv p_i$ and \equationn\ give
$$\sum_{j=1}^g\xi_{ij}\omega_j(z)=
{\theta_\Delta(\a_i+z-y)\theta_\Delta(\a-z)\over
\theta_\Delta(\a-y)}{E(p_i,y)\over E(p_i,z) E(y,z)}\ ,
$$
$i\in I_g$, that, if $\det\xi_{ij}\neq 0$, coincides with
\oomegai.

Furthermore, by Fay's trisecant identity it immediately follows that
if $\alpha$ is a non-singular element of $\Theta$, then for any four points
$x_i\in C$, $i\in I_4$, the ``cross ratio''
$$
{\theta(\alpha+x_1-x_2)\theta(\alpha+x_3-x_4)\over \theta(\alpha+x_1-x_4)\theta(\alpha+x_3-x_2)}\ ,
$$
is independent of $\alpha$. In particular, by \thetaconc\
$$
{\theta(\alpha_i+p-q)\theta(\alpha_i+x-z)\over \theta(\alpha_i+p-z)\theta(\alpha_i+x-q)}={\theta[\nu](p-q)\theta[\nu](x-z)\over\theta[\nu](p-z)\theta[\nu](x-q)}\ ,
$$
holds for all $p,q,x,z\in C$ and $\nu$, non-singular odd theta characteristics and $\alpha_i$, $i\in I_g$, a non-singular element of
$\Theta$. In the limit $x\to z$, we have
\eqn\ghghgt{
{\theta(\alpha_i+p-q)
\over \theta(\alpha_i+p-z)\theta(\alpha_i+z-q)}\sum_{j=1}^g\theta_j(\alpha_i)\omega_j(z)
={\theta[\nu](p-q)\over\theta[\nu](p-z)\theta[\nu](z-q)}\sum_{j=1}^g\theta_j[\nu](0)\omega_j(z)\ ,
}
and the first equation in \oomegaiuuu\ follows by noticing that the right hand side of \ghghgt\ coincides with $E(q,p)/E(z,p)E(q,z)$.
The second equation follows by replacing in the above equations $\alpha_i$ by a non-singular theta characteristic $\nu_i$. Each theta function gets a phase factor
that however cancels in the cross ratio. Alternatively, the second equation in \oomegaiuuu\ can be obtained by comparing the expression
for $E(q,p)/E(z,p)E(q,z)$, given by the right hand side of \ghghgt, with the analogous expression with $\nu$ replaced by $\nu_i$. Incidentally,
this shows that independence of the prime form on the specific non-singular odd spin structure is a property related to the Fay trisecant
identity.
It is clear that the conditions $\det\xi_{ij}\neq0$ and $\det\lambda_{ij}\neq0$
can be always satisfied. Existence of
odd theta characteristics such that $\det \mu_{ij}\neq0$, follows by Lefschetz Theorem for abelian varieties which implies
that the rank of the rectangular matrix $\theta_i[\nu_j](0)$, $i\in I_g$,
with $\{\nu_j\}_{j\in I_{2^{g-1}(2^g-1)}}$
the full set of odd characteristics, is always of rank $g$
\ref\SalvatiManni{R.~Salvati Manni, On the nonidentically zero Nullwerte of Jacobians of
theta functions with odd characteristics, {\it Adv. in Math.} {\bf 47} (1983), 88-104.}.\hfill$\square$

\vskip 6pt

An immediate consequence of the previous theorem is that the determinant of the holomorphic 1-differentials can be expressed without the use
of any constant, such as $\kuno[\omega]$ and of the $\sigma$-differential.

\newth\definiT{Corollary}{Let $\nu_i$, $i\in I_g$, be odd theta characteristics
satisfying $\det \theta_i[\nu_j](0)\neq0$, then \eqn\quessstyy{\det
\omega_i(p_j)={1\over \det \theta_i[\nu_j](0)}
{E(p,q)^g\det(\theta[\nu_i](p_j-q)\theta[\nu_i](p-p_j))\over
\prod_{i=1}^g E(p,p_i)E(q,p_i)\theta[\nu_i](p-q)}\ , } holds for all
$p,q,p_i\in C$, $i\in I_g$.}

\vskip 6pt

\noindent {\bf Proof.} Immediate by \oomegaiuuu. $\hfill\square$

\vskip 6pt

A basic implication of the above relations is that also
$\kuno[\omega]$ can be expressed without the use of
$\det\omega_i(z_j)$.

\vskip 6pt

\newth\spettacolare{Corollary}{
\eqn\kappita{\kuno[\omega]=(-)^{g(g-1)/2}\Bigl({\theta_\Delta(\a-y)\over
\sigma(y)\prod_kE(y,p_k)}\Bigr)^{g-1}{\prod_{i<j}E(p_i,p_j)\over
\det\theta_{\Delta,j}(a_i)}\ .} }

\vskip 6pt

\noindent {\bf Proof.} Immediate by \rigira. $\hfill\square$

\vskip 6pt

\newrem\progresso{Remark}{Such a relation was previously known only for $g=1$ \FayMAM, a consequence of the triviality
of $\det\omega_i(z_j)$ on the torus.}

\vskip 6pt

Let $\nu_1,\ldots,\nu_g$ be $g$ odd spin structures and denote by
$p_i^j$, $j\in I_{g-1}$, the corresponding points on $C$, that is
\eqn\strutture{ \nu_i''+\tau\nu_i'=I(\smsum_{j=1}^{g-1}p_i^j-\Delta)
\ , } $i\in I_g$. Set $h_i\equiv h_{\nu_i}$, and recall that
$(h_i)=\sum_{j=1}^{g-1}p_i^j$. Define the symbol
\eqn\definethesymbol{ \{p_i^j\}_n:=\kuno[\omega]^{g(1-g)}e^{\pi
i(g-1)\sum_1^g\tp\nu_i'\tau\nu_i''+2\pi
i(g-1)\sum_1^g\tp\nu_i'\nu_i''}\prod_{l=1}^{g-1}
\prod_{k=1}^g\det_{ij}\omega_i(p_k^j)_{|_{p_k^g\equiv p_{k+n}^l}}\ ,
} $n\in I_{g-1}$, $p_{g+l}^j\equiv p_l^j$, $l\in I_n$, $j\in
I_{g-1}$,  which is a holomorphic $g$-differential in each one of
the $p_i^j$'s. By \ildet\ \eqn\primomodo{
\det_{ij}\omega_i(p_k^j)=(-)^{g-1}e^{-\pi i \tp \nu_k'\tau\nu_k'
-2\pi i \tp
\nu_k'\nu_k''}\kuno[\omega]h_k^2(p_k^g)\prod_{i<j=1}^{g-1}E(p_k^i,p_k^j)\prod_{i=1}^{g-1}\sigma(p_k^i)
\ . } It follows that the expression of \definethesymbol\ for, say
$n=1$, is equivalent to \eqn\outstand{\eqalign{ \{p_i^j\}_1& =
\prod_{i=1}^{g-1}h_1^2(p_2^i)h_2^2(p_3^i)\cdots
h_g^2(p_1^i)\prod_{k=1}^g\prod_{i<j=1}^{g-1}
E(p_k^i,p_k^j)^{g-1}\prod_{i=1}^{g-1}\sigma(p_k^i)^{g-1}\cr\cr & =
\prod_{k=1}^g\prod_{i=1}^{g-1}\sigma(p_k^i)^{g-1}
\prod_{i<j=1}^{g-1} {\theta[\nu_k](p_{k+1}^j-p_{k+1}^i)^{g-1}\over
h_k(p_{k+1}^j)^{g(g-3)}}
 \ .}}
Such an expression contains only points belonging to the divisors
defining non-singular odd spin structures. A more symmetric expression is obtained by
multiplying over $n$
\eqn\bellissimissima{
\prod_{n=1}^{g-1}\{p_i^j\}_n =
\prod_{k=1}^g\prod_{i=1}^{g-1}\sigma(p_k^i)^{(g-1)^2}
\prod_{n=1}^{g-1}\prod_{i<j=1}^{g-1}
{\theta[\nu_k](p_{k+n}^j-p_{k+n}^i)^{g-1}\over h_k(p_{k+n}^j)^{g(g-3)}}
 \ .}

According to Remark \lasigmanu\ Eqs.\dettheta\ and \detthetaii\
define new constants $\kenne_\nu[\phi^n]$, once the
$g/2$-differential $\sigma$ in \dettheta\ and \detthetaii\ is
replaced by $\sigma_\nu$. Explicitly, there are constants
$\kenne_\nu[\phi^n]$ depending only on the marking of $C$ and on
$\{\phi_i^n\}_{i\in I_{N_n}}$ such that
\eqn\detthetaAAA{\kuno_\nu[\phi^1]={
\det\phi_i^1(p_j)\sigma_\nu(y)\prod_1^gE(y,p_i)\over
\deltadiv\bigl(\sum_{1}^gp_i-y\bigr)\prod_1^g\sigma_\nu(p_i)
\prod_{i<j}^gE(p_i,p_j) }\ ,} and
\eqn\detthetaiiAAA{\kenne_\nu[\phi^n]={\det\phi_i^{n}(p_j)\over
\theta_\Delta\bigl(\sum_{1}^{N_n}
p_i\bigr)\prod_{1}^{N_n}\sigma_\nu(p_i)^{2n-1}\prod_{i<j}^{N_n}
E(p_i,p_j)}\ ,}  for $n\geq2$, for all $y,p_1,\ldots,p_{N_n}\in C$.
As we noticed, a property of $\sigma_\nu$ is that, unlike $\sigma$,
it is directly defined on $C$. This allows us to express
$\kuno_\nu[\omega]$ in terms of theta functions only. Most
importantly, products of $\kuno_{\nu_k}[\omega]$ have expressions
where only points belonging to divisors defining odd spin structures
appear.

\vskip 6pt

\newth\basilareeAA{Theorem}{The relation
\eqn\nuooovaa{\kuno_\nu[\omega]=e^{-\pi i \tp \nu'\tau \nu'}
{\prod_{i<j}^g\theta[\nu](p_i-p_j)\exp{\sum_{j=1}^{g-1}({2\pi i\over
g-1} \tp\nu'\K^{y_j}+\sum_{i=1}^g\oint_{\alpha_i} \omega_i(w)\log
\theta[\nu](w-y_j))}\over\prod_{i,j=1}^{g-1}\theta[\nu](p_i-y_j)\det_{ij}\theta_j[\nu](p_g-p_i)}
\ , } holds for all the non-singular odd spin structures
$$
\nu''+\tau\nu'=I(\smsum_{i=1}^{g-1}p_i-\Delta) \ ,
$$
and $p_g,y_i$, $i\in I_{g-1}$, arbitrary points of $C$.
In particular, if
$$\nu_k''+\tau\nu_k'=I(\smsum_{i=1}^{g-1}p_k^i-\Delta) \ ,
$$
$k\in I_n$, $n>2$, are non-singular odd spin structures, then
\eqn\ogep{\prod_{k=1}^n \kuno_{\nu_k}[\omega] \ , } where
$y_k^i\equiv p_{k+1}^i$, $p_k^g\equiv p_{k+2}^1$, $k\in I_n$, $i\in
I_{g-1}$, and $p^i_{n+j}\equiv p^i_n$, $j=1,2$, depends only on the
divisors associated to $\{\nu_k\}_{k\in I_n}$. }

\vskip 6pt

\noindent {\bf Proof.} First note that
$$\kuno_\nu[\omega]\sigma_\nu(y)^{g-1}=\kuno[\omega]\sigma(y)^{g-1}\ ,$$
so that Corollary \spettacolare\ still holds whether $\kuno[\omega]$
and $\sigma$ are replaced by $\kuno_\nu[\omega]$, defined in
\detthetaAAA, and $\sigma_\nu$, defined in \nuovasigma,
respectively. Next, identify $a$ with the divisor $\sum_{i=1}^g
p_i$, with the first $g-1$ points, defining, as in \strutture, the
non-singular odd spin structure $\nu$. It then follows that
$$
I(a-p_g-\Delta)=I(a_i+p_i-p_g-\Delta)=\nu''+\tau\nu' \ ,
$$
so that, recalling the expression of the prime form \pojdlk,
\eqn\alllg{
{\theta_\Delta(a-y)\over E(y,p_g)}=e^{-\pi i \tp \nu'\tau \nu'-2\pi i \tp \nu'(p_g-y+\nu'')}h(p)h(y) \ ,
}
and
$$
\det \theta_{\Delta,j}(a_i)=e^{-\pi ig \tp \nu'\tau \nu'-2\pi i \tp \nu'(gp_g-\sum_{i=1}^gp_i+g\nu'')}\det_{ij}\theta_j[\nu](p_g-p_i) \ .
$$
On the other hand, by \lakkapppaa\
$$I(p_g-a+(g-1)p_g)=-\nu''-\tau\nu' + I((g-1)p_g-\Delta)=
-\nu''-\tau\nu'+\K^{p_g}\ ,
$$
so that, if $a=\sum_{i=1}^{g}p_i$
and $\nu''+\tau\nu'=I(\sum_{j=1}^{g-1}p_j-\Delta)$, then $\det \theta_{\Delta,j}(a_i)$ admits the following nice expression
\eqn\niceee{
\det \theta_{\Delta,j}(a_i)=e^{-\pi i (g-2)\tp \nu'\tau \nu'-2\pi i(g-1)\tp \nu'\nu''-2\pi i\tp\nu'\K^{p_g}}\det_{ij}\theta_j[\nu](p_g-p_i) \ .
}
Replacing \alllg\ and \niceee\ in \kappita, with $\kuno[\omega]$ and $\sigma$ replaced by
$\kuno_\nu[\omega]$ and $\sigma_\nu$, respectively, and then observing that
$$(-)^{g(g-1)/2}\prod_{i<j}^{g-1}\theta[\nu](p_j-p_i)\prod_{i=1}^{g-1}\theta[\nu](p_g-p_i)
=(-)^{g(g-1)}\prod_{i<j}^g\theta[\nu](p_i-p_j)\ ,$$ and that
$(-)^{g(g-1)}=1$, yields \eqn\interrr{\kuno_\nu[\omega]=e^{-\pi i
\tp \nu'\tau \nu' +2\pi i \tp\nu'\K^{y}}
{\prod_{i<j}^g\theta[\nu](p_i-p_j)e^{(g-1)\sum_{i=1}^g\oint_{\alpha_i}
\omega_i(w)\log \theta[\nu](w-y)} \over
\prod_{i=1}^{g-1}\theta[\nu](p_i-y)^{g-1}\det_{ij}\theta_j[\nu](p_g-p_i)}
\ . } Next, independence of $\kuno_\nu[\omega]$ on the points,
implies the identity \eqn\identita{ {e^{2\pi i
\tp\nu'\K^{y}+(g-1)\sum_{i=1}^g\oint_{\alpha_i} \omega_i(w)\log
\theta[\nu](w-y)}\over\prod_{i=1}^{g-1}\theta[\nu](p_i-y)^{g-1}}=
{e^{\sum_{j=1}^{g-1}({2\pi i\over g-1}
\tp\nu'\K^{y_j}+\sum_{i=1}^g\oint_{\alpha_i} \omega_i(w)\log
\theta[\nu](w-y_j))}\over\prod_{i,j=1}^{g-1}\theta[\nu](p_i-y_j)}\
,} that, together with \interrr, reproduces \nuooovaa.

\vskip 6pt

\newrem\reproducing{Remark}{Theorem \basilareeAA\ also provides an immediate derivation of the only previously known solution, namely
the one of the elliptic
case. On the other hand, $g=1$ corresponds to the only case where $\kuno[\omega]$ can be defined independently of $\sigma$, and so it does not need
the definition on the universal covering of $C$. In the elliptic case we in fact have $\kuno[\omega]=\kuno_\nu[\omega]$. By \scriviamolo\
$$\K={1+\tau\over2}\ ,$$
and using $\theta_1\left[^{1/2}_{1/2}\right](0)=-2\pi \eta^3(\tau)$,
Eq.\interrr\ yields \eqn\niceverifica{ \kuno_\nu[\omega]={e^{-\pi i
\tp \nu'\tau \nu' +2\pi i
\tp\nu'\K}\over\theta_1\left[^{1/2}_{1/2}\right](0)}=-{e^{-{\pi
i\over2}(1-{\tau\over2})}\over 2\pi \eta^3(\tau)} \ , } that
coincides with the derivation in \FayMAM\ (see pg.21 there and note
that $\kuno[\omega]$ corresponds to $\kuno_0^{-1}$). We also note
that \ogep\ admits alternative versions depending on the choice of
the identification among the points $p_k^i$ and $y_k^i$.}

\vskip 6pt

An analysis similar to the one for $\kuno_\nu[\omega]$ can be done
also for $\kenne_\nu[\phi^n]$, $n\geq2$. However, in this case the
degree $N_n$ of the effective divisor appearing in \detthetaiiAAA\
is a multiple of $g-1$, so that $I(\sum^{N_n}_1p_i-(2n-1)\Delta)$
can be directly identified with the summation of $2n-1$ non-singular
odd spin structures. It is convenient to change the notation and to
denote the set of points $\{p_i\}_{i\in I_{N_n}}$ as
$\{p_i^j\}_{i\in I_{2n-1},j\in I_{g-1}}$, so that
\eqn\dddf{\nu_k''+\tau\nu_k'=I(\smsum_{i=1}^{g-1}p_k^i-\Delta) \ , }
$k\in I_{2n-1}$, $n\geq2$. The previous investigation naturally
leads to introduce the following constants.

\vskip 6pt

\newrem\thdetthetaBBBABBO{Definition}{Fix $n\ge 2$ such that a set $\{\nu\}\equiv
\{\nu_k\}_{k\in I_{2n-1}}$ of non-singular odd spin structures, as
in \dddf, is defined and set
\eqn\dettheaaaBV{\kenne_{\{\nu\}}[\phi^n]:={\det\phi_i^n(p_j)
\exp[(2n-1){\sum_{k=1}^{2n-1}\sum_{j=1}^{g-1}\sum_{i=1}^g\oint_{\alpha_i}}\omega_i(w)\log
\theta[\nu_k](w-p_k^j)]\over \theta\bigl(\sum_{1}^{2n-1}
\nu_k\bigr)\prod_{k<l=1}^{2n-1} \prod_{i,j=1}^{g-1}
E(p_k^i,p_l^j)\prod_{k=1}^{2n-1}\prod_{i=1}^{g-1}\bigl(h_{\nu_k}(p_k^i)^{(2n-1)g}
\prod_{j>i}^{g-1} E(p_k^i,p_k^j)\bigr)}\ ,} where the points $p_i$,
$i\in I_{N_n}$, in the determinant are given by $p_{i+j(g-1)}\equiv
p_j^i$, $i\in I_{g-1}$, $j\in I_{2n-1}$.}

\vskip 6pt

\newsec{Siegel metric induced on $\hat{\M}_g$ and Bergman reproducing kernel}\seclab\Siegsec

In this section we derive the explicit expression of the metric $ds^2_{|\hat\M_g}$
on the moduli space $\hat\M_g$ of genus $g$ canonical curves induced by the Siegel metric. This
was previously known only for the trivial cases $g=2$ and $g=3$. By Wirtinger Theorem the explicit expression
for the volume form on $\hat\M_g$ is also obtained. A remarkable property of $ds^2_{|\hat\M_g}$ is that
it is given by the Kodaira-Spencer map of the square of the Bergman reproducing kernel (times $4\pi^2$).
This is one of the basic properties of the Bergman reproducing kernel
derived in this section. Such an approach will led to a
notable relation satisfied by the determinant of powers of the Bergman reproducing kernel.
The results are a natural consequence of the present approach,
which also uses, as for the present derivation of $ds^2_{|\hat\M_g}$,
the isomorphisms introduced in section \notation.

The Torelli space $\T_g$ of smooth algebraic curves of genus $g$ can
be embedded in $\H_g$ by the period mapping, which assigns to a
curve $C$, with a fixed basis of $H_1(C,\ZZ)$, representing a point
in $\T_g$, the corresponding period matrix. The period mapping has
maximal rank $3g-3$ on the subspace $\hat\T_g$ of non-hyperelliptic
curves and therefore a metric on $\H_g$ induces the pull-back metric
on $\hat\T_g$. It is therefore natural to consider the Siegel metric
on $\H_g$ \ref\siegelajm{C.~L. Siegel, Symplectic Geometry, {\it
Amer.\ J.\ Math.\ } {\bf 65} (1943), 1-86} \eqn\siegg{ds^2:=\Tr\,
(Y^{-1}dZ Y^{-1}d\bar Z) \ ,} where $Y:=\im Z$, $Z\in \H_g$. Such a
metric is ${\rm Sp}(2g,\RR)$ invariant, and since
$\hat\M_g\cong\hat\T_g/\Gamma_g$, it also induces a metric on
$\hat\M_g$. The Siegel volume form is \siegelajm\ \eqn\sieggvv{
d\nu={i^M\over 2^g}{\bigwedge_{i\le j}^g ({\rm d} Z_{ij}\wedge{\rm
d}\bar Z_{ij})\over(\det Y)^{g+1} }\ . }

The explicit expression of the volume forms on $\hat{\M}_g$ induced
by the Siegel metric, which coincides with \sieggvv\ for $g=2$ and
$g=3$ non-hyperelliptic curves, is given in Theorem \futuref{volumeSieg}. It
is simply written in terms of the Riemann period matrix $\tau_{ij}$ and of the
basis $\{d\tau_{ij}\}$ of $T^*\hat\T_g$.

The Laplacian associated to the Siegel's symplectic metric were
derived, ten years after Siegel's paper \siegelajm, by H. Maass
\ref\maass{H. Maass, Die Differentialgleichungen in der Theorie der
Siegelschen Modulfuktionen, {\it Math. Ann.} {\bf 126} (1953),
44-68.} \eqn\Maass{\Delta=4\Tr
\Big(Y\,{}^{{}^t}\Big(Y{\partial\over\partial\bar Z}\Big){\partial
\over\partial Z}\Big) \ . } As we will see, as a byproduct of the
present approach, and of the formalism developed in section
\notation\ in particular, both \sieggvv\ and \Maass\ are straightforwardly derived.

\vskip 6pt

\subsec{Derivation of the volume form and the Laplace-Beltrami
operator on $\H_g$}

\newth\xyomegasss{Proposition}{The Siegel metric \siegg\ can be equivalently
expressed in the form
\eqn\nvbh{ds^2=\sum_{i,j=1}^Mg^S_{ij}dZ_id\bar
Z_j\ ,} where \eqn\pulis{g^S_{ij}(Z,\bar
Z):=2{Y^{-1}_{\1_i\1_j}Y^{-1}_{\2_i\2_j}+
Y^{-1}_{\1_i\2_j}Y^{-1}_{\2_i\1_j}\over
(1+\delta_{\1_i\2_i})(1+\delta_{\1_j\2_j})}=2\chi_i^{-1}\chi_j^{-1}
(Y^{-1}Y^{-1})_{ij}\ ,} $i,j\in I_M$.}

\vskip 6pt

\noindent {\bf Proof.} For $n=2$ the identity \identittddd\ reads
$$
\sum_{i,j=1}^gf(i,j)=\sum_{k=1}^M(2-\delta_{\1_k\2_k})f(\1_k,\2_k)\ ,$$
where we used the identity
$$2-\delta_{ij}={2\over 1+\delta_{ij}}\ .$$ Hence
\eqn\great{\eqalign{ds^2&=\sum_{i,j,k,l=1}^gY^{-1}_{ij}d Z_{jk}Y^{-1}_{kl}d\bar Z_{li}\cr
&=\sum_{i,j=1}^gd\bar Z_{ji}\sum_{m=1}^M{Y^{-1}_{i\1_m}Y^{-1}_{j\2_m}+
Y^{-1}_{i\2_m}Y^{-1}_{j\1_m}\over 1+\delta_{\1_m\2_m}}d Z_{\1_m\2_m}\cr
&=\sum_{m,n=1}^M(2-\delta_{\1_n\2_n})d\bar Z_{\1_n\2_n}{Y^{-1}_{\1_n\1_m}
Y^{-1}_{\2_n\2_m}+
Y^{-1}_{\1_n\2_m}Y^{-1}_{\2_n\1_m}\over 1+\delta_{\1_m\2_m}}d Z_{\1_m\2_m}\cr
&=\sum_{m,n=1}^M2\chi_m^{-1}\chi_n^{-1}(Y^{-1}Y^{-1})_{nm}d Z_md\bar Z_n\ .}}
\hfill$\square$

\vskip 6pt

\noindent Let \eqn\formomega{\omega:={i\over
2}\sum_{i,j=1}^Mg^S_{ij}dZ_i\wedge d\bar Z_j\ ,} be the $(1,1)$-form
associated to the Siegel metric on $\H_g$, so that the volume form
on $\H_g$ is
$${1\over M!}\omega^M=\Bigl({i\over 2}\Bigr)^M
\det g^S_{ij}\bigwedge_{i\le j}^g ({\rm d} Z_{ij}\wedge{\rm d}\bar
Z_{ij})\ .$$

\newth\volumeHg{Proposition}{
$$\det g^S_{ij}={2^{M-g}\over (\det Y)^{g+1}} \ . $$}

\vskip 6pt

\noindent {\bf Proof.} Since $Y$ is symmetric and positive-definite,
we have $PY^{-1}P^{-1}={\rm
diag}\,(\lambda_1,\ldots,\lambda_g)\equiv D$, for some non-singular
$g\times g$ matrix $P$ and some positive
$\lambda_1,\ldots,\lambda_g$. By \pulis\ and \formdet\
$$\eqalignno{\det g^S_{ij}=&2^M\det\nolimits_{ij}\bigl((Y^{-1}Y^{-1})_{ij}\chi^{-1}_i\chi^{-1}_j\bigr)\cr\cr
=&2^M\det\nolimits_{ij}\bigl((PP)_{ij}\chi_j^{-1}\bigr)\det\nolimits_{ij}\bigl((P^{-1}
P^{-1})_{ij}\chi_j^{-1}\bigr)\det\nolimits_{ij}\bigl((Y^{-1}Y^{-1})_{ij}\chi^{-1}_i\chi^{-1}_j\bigr)\
,}
$$
and by \formprod\
$$\det g^S_{ij}=2^M\det\nolimits_{ij}\bigl((DD)_{ij}\chi^{-1}_i\chi^{-1}_j\bigr)=
2^M\det\nolimits_{ij}\bigl(\lambda\lambda_i(\delta\delta)_{ij}\chi^{-1}_i\chi^{-1}_j\bigr)\
.$$ The proposition then follows observing that
$(\delta\delta)_{ij}=\chi_j\delta_{ij}$ and that \formsimprod\
yields
$$\det g^S_{ij}=2^M\prod_{i=1}^M\lambda\lambda_i\chi^{-1}_i=2^{M-g}\Bigl(\prod_{k=1}^g\lambda_k\Bigr)^{g+1}\ .$$
\hfill$\square$

\vskip 6pt

\newth\laplac{Proposition}{The Laplace-Beltrami operator acting on functions on $\H_g$ is
$$\Delta={1\over 2}\sum_{i,j=1}^M(YY)_{ij}{\partial\over\partial Z_i}{\partial\over\partial
\bar Z_j}\ .$$ }

\vskip 6pt

\noindent{\bf Proof.} Just use the definition of $\Delta$ and note
that $g^{S\,ij}=(YY)_{ij}/2$.\hfill$\square$

\vskip 6pt

\subsec{Basis of the fiber of $T^*\hat\T_g$ and Siegel metric on $\hat\M_g$}\subseclab\Kodaira

The following theorem provides a modular invariant basis of the fiber of
$T^*\hat\T_g$ at the point representing $C$.

\newth\thdexi{Theorem}{If $p_3,\ldots,p_g\in C$ are $g-2$ pairwise distinct points
such that $K(p_3,\ldots,p_g)\neq 0$, then
\eqn\dexi{
\Xi_i:= \sum_{j=1}^MX^\omega_{ji}d\tau_j\ ,} $i\in I_N$, with $X^\omega_{ij}$, $i,j\in I_M$, defined
in Eq.\leX, is a modular invariant basis of the fiber of $T^*\hat\T_g$
at the point representing $C$.}

\vskip 6pt

\noindent {\bf Proof.}
Consider the Kodaira-Spencer map $k$ identifying the space
of quadratic differentials on $C$ with the fiber of the cotangent
bundle of $\M_g$ at the point representing $C$. Next, consider a Beltrami differential
$\mu\in\Gamma(\bar K_C\otimes K_C^{-1})$ (see \ref\BonoraCJ{
  L.~Bonora, A.~Lugo, M.~Matone and J.~Russo,
  A global operator formalism on higher genus Riemann surfaces: b-c
  systems,
  { \it Commun.\ Math.\ Phys.\ } {\bf 123} (1989), 329-359.}
for explicit constructions) and recall that it defines a tangent vector at $C$ of $\T_g$.
The derivative of the period map $\tau_{ij}:\T_g\to \CC$ at $C$ in the direction
of $\mu$ is given by Rauch's formula
$$
d_C\tau_{ij}(\mu)=\int_C\mu\omega_i\omega_j \ .$$
It follows that
$$k(\omega_j\omega_k)={1\over2\pi i}d\tau_{jk} \ ,$$ $j,k\in I_g$,
so that, by \vXww,
\eqn\dxiiilll{k(v_j)={1\over2\pi i}\sum_{k=1}^M X_{kj}^\omega d\tau_k \ ,} $j\in I_N$,
where
$$d\tau_i:=d\tau_{\1_i\2_i}\ ,$$ $i\in I_M$.
It follows that the differentials
\eqn\dxiii{\Xi_j:=2\pi i\,k(v_j) \ ,} $j\in I_N$, are linearly
independent. Furthermore, since by construction the basis
$\{v_i\}_{i\in I_N}$ is independent of the choice of a symplectic
basis of $H_1(C,\ZZ)$, such differentials are modular invariant,
i.e. \eqn\vacitata{\Xi_i\;\mapsto\; \tilde\Xi_i=\Xi_i \ ,} $i\in I_N$, under
\modull. \hfill$\square$

\vskip 6pt

Let $ds_{|\hat\M_g}^2$ be the metric on
$\hat{\cal M}_g$ induced by the Siegel metric.
Set
\eqn\giesse{g^\tau_{ij}:=g_{ij}^S(\tau,\bar\tau)=2\chi_i^{-1}\chi_j^{-1}
({\imtau}^{-1}{\imtau}^{-1})_{ij}\ .}

\vskip 6pt

\newth\thdenudue{Corollary}{\eqn\tog{ds_{|\hat\M_g}^2= \sum_{i,j=1}^Ng^\Xi_{ij}\Xi_i\bar\Xi_j
\ ,} where
\eqn\giXi{g^\Xi_{ij}:=\sum_{k,l=1}^Mg^\tau_{kl}B^\omega_{ik}\bar
B^\omega_{jl}\ ,} and $B^\omega$ is the matrix defined in {\rm \leB}
with $\eta_i\equiv\omega_i$, $i\in I_g$. Furthermore, the volume
form on $\hat\M_g$ induced by the Siegel metric is
 \eqn\denudue{d\nu_{|\hat\M_g}=\Bigl({i\over 2}\Bigr)^N\det
g^\Xi\, dw\wedge d\bar w\ ,}
where
\eqn\dewu{dw:=\sum_{i_N> \ldots> i_1=1}^M\Xmin{\omega}{\ss 1 \hfill
\ldots \hfill N \cr \ss i_1\hfill \ldots \hfill
i_N}\,d\tau_{i_1}\wedge\cdots \wedge d\tau_{i_N}\ .}}

\vskip 6pt

\noindent {\bf Proof.}
By \nvbh\ and \pulis\
\eqn\together{ds_{|\hat\M_g}^2=\sum_{k,l=1}^Mg^\tau_{ij}d\tau_{i}d\bar\tau_{j}\
.} Furthermore, by applying the Kodaira-Spencer map to both sides of
Eq.\teor, one obtains
\eqn\domega{d\tau_i=\sum_{j=1}^NB^\omega_{ji}\,\Xi_j \ ,}
$i\in I_M$, and \tog\ follows.
On the other hand, by \tog\
\eqn\ds{d\nu_{|\hat\M_g}=
\Bigl({i\over 2}\Bigr)^N
\det g^\Xi\ {\textstyle\bigwedge_{1}^N}\!(\Xi_i\wedge\bar\Xi_i)\ ,}
and by Theorem \thdexi\ the proof is completed. \hfill$\square$

\vskip 6pt

\noindent Applying the Kodaira-Spencer map to \coroll\ yields the
linear relations satisfied by $d\tau_i$, $i\in I_M$.

\vskip 6pt

\newth\thfinale{Corollary}{The $(g-2)(g-3)/2$ linear relations \eqn\finale{\sum_{j=1}^MC^\omega_{ij}
d\tau_j=0\ ,} $N+1\leq i\le M$, where the matrices $C^\omega$ are defined
in \leCi, hold.}

\vskip 6pt
Set $\imtau:=\im \tau$ and consider the Bergman reproducing kernel
$$B(z,\bar w):=\sum_{i,j=1}^g \omega_i(z)(\imtau^{-1})_{ij}\bar\omega_j(w)\ ,$$ for all $z,w\in C$., and
Set ${\cal K}(\phi\bar\psi):=k(\phi)\bar k(\psi)$, $\bar
k(\bar\psi)=\overline{k(\psi)}$, for all $\phi,\psi\in H^0(K_C^2)$,
where $k$ is the Kodaira-Spencer map.

\vskip 6pt

\newth\volumeSieg{Theorem}{\eqn\daiie{ds^2_{|\hat\M_g}=4\pi^2{\cal K}(B^2) \ .}
Furthermore, the volume form on $\hat\M_g$ induced by the Siegel metric is
\eqn\basicc{d\nu_{|\hat\M_g}=i^N\sum_{{i_N>\ldots>i_1=1\atop
j_N>\ldots>j_1=1}}^M{\left|\imtau^{-1}\imtau^{-1}\right|^{i_1\ldots
i_N}_{j_1\ldots j_N}\over
\prod_{k=1}^N(1+\delta_{\1_{i_k}\2_{i_k}})(1+\delta_{\1_{j_k}\2_{j_k}})}\bigwedge_{l=1}^N(d\tau_{i_l}\wedge
d\bar\tau_{j_l})\ . }}

\vskip 6pt

\noindent {\bf Proof.} Eq.\daiie\ is an immediate consequence of Proposition \xyomegasss\ and of the
application of the Kodaira-Spencer map to the identity
\eqn\idennn{\sum_{i,j=1}^M\omega\omega_i(z)g^\tau_{ij}\bar\omega\bar\omega_j(w)=B^2(z,\bar w)\ .}
Consider the $(1,1)$-form $\omega$ defined in Eq.\formomega. By Wirtinger's Theorem
\ref\GH{P.A.~Griffiths and J.~Harris, {\it Principles of Algebraic Geometry}, Wiley, 1978.}, the
volume form on a $d$-dimensional complex submanifold $S$ is
$${1\over d!}\omega^d\ ,$$
so that the volume of $S$ is expressed as the integral over $S$ of a globally defined differential form on $\H_g$.
Note that
$$\eqalignno{d\nu_{|\hat\M_g}=&{i^N\over 2^NN!}\sum_{{i_1,\ldots,i_N=1\atop j_1,\ldots,j_N=1}}^M
\prod_{k=1}^Ng^\tau_{i_kj_k}\bigwedge_{l=1}^N(d\tau_{i_l}\wedge
d\bar \tau_{j_l})\cr\cr
=&{i^N\over 2^NN!}\sum_{{i_N>\ldots>i_1=1\atop j_N>\ldots>j_1=1}}^M\sum_{r,s\in\perm_N}
\epsilon(r)\epsilon(s)\prod_{k=1}^Ng^\tau_{i_{r(k)}j_{s(k)}}\bigwedge_{l=1}^N(d\tau_{i_l}\wedge
d\bar \tau_{j_l})\ ,
}$$
and Eq.\basicc\ follows by the identity
$$\sum_{r,s\in\perm_N}
\epsilon(r)\epsilon(s)\prod_{k=1}^Ng^\tau_{i_{r(ks)}j_{s(k)}}=N!\left|g^\tau\right|^{i_1\ldots
i_N}_{j_1\ldots j_N}\ .$$
\hfill$\square$

\vskip 6pt

Fix the points $z_1,\ldots,z_N\in C$ satisfying the conditions of
Proposition \thnewbasis. The basis $\{\gamma_i\}_{i\in I_N}$ of
$H^0(K_C^2)$, with $\gamma_i\equiv\gamma^2_i$, $i\in I_N$, defined
by Eq.\basendiff\ in the case $n=2$, satisfies the relations
$$\omega\omega_i=\sum_{j=1}^N\omega\omega_i(z_j)\gamma_j\ ,\qquad
v_i=\sum_{j=1}^Nv_i(z_j)\gamma_j\ ,$$ $i\in I_M$. Set
$\Gamma_i:=(2\pi i)^{-1}k(\gamma_i)$ and $[v]_{ij}:=v_i(z_j)$, $i,j\in I_N$.

\vskip 6pt

\newth\Bergman{Corollary}{Fix the points $z_1,\ldots,z_N\in C$ in such a way that
$\det\phi_i(z_j)\neq 0$, for any arbitrary basis $\{\phi_i\}_{i\in
I_N}$ of $H^0_C(K^2)$. The metric on $\hat\M_g$ induced by the
Siegel metric is \eqn\siegelbergman{
ds^2_{|\hat\M_g}=\sum_{i,j=1}^NB^2(z_i,\bar
z_j)\Gamma_i\bar\Gamma_j\ , } and the volume form is
\eqn\volumee{d\nu_{|\hat\M_g}=\Bigl({i\over 2}\Bigr)^N\det\,
B^2(z_i,\bar z_j)\ {\textstyle\bigwedge_1^N}\!(\Gamma_i\wedge \bar
\Gamma_i) =\Bigl({i\over 2}\Bigr)^N\,{\det B^2(z_i,\bar z_j)\over
|\det v_i(z_j)|^2}\, dw\wedge d\bar w\ ,} where $\{v_i\}_{i\in I_N}$
is the basis of $H^0(K_C^2)$ defined in Proposition {\rm
\futuref{thlev}} and $dw$ is defined in Eq.{\rm \futuref{dewu}}. }

\vskip 6pt

\noindent {\bf Proof.} Eq.\siegelbergman, and therefore the first equality in
Eq.\volumee, follows substituting
$$d\tau_i=\sum_{j=1}^N\omega\omega_i(z_j)\Gamma_j\ ,$$ $i\in I_M$,
in \together\ and then using the identity \idennn. Next, note that
comparing \siegelbergman\ and \tog, and by
$\Xi_i=\sum_{j=1}^N[v]_{ij}\Gamma_j$, $i\in I_N$, yields
$$\sum_{k,l=1}^N[v]_{ki}g^\Xi_{kl}[\bar v]_{lj}=B^2(z_i,\bar z_j)\
,$$ which also follows by the definition \giXi\ of $g^\Xi$ and by
Eq.\teor, with $\eta_i\equiv\omega_i$, $i\in I_g$, and Eq.\idennn.
Hence \eqn\ildettt{\det\, g^\Xi = {\det B^2(z_i,\bar z_j)\over|\det
v_i(z_j)|^2}\ ,} which also follows by $\det\gamma_i(z_j)=1$ and
$$
\Xi_1\wedge\cdots\wedge\Xi_N=\det\, v_i(z_j)\Gamma_1\wedge\cdots\wedge\Gamma_N \ ,
$$
and the second equality in Eq.\volumee\ follows. \hfill$\square$

\vskip 6pt

\subsec{Determinants of powers of the Bergman reproducing kernel}

Corollary \Bergman, in particular Eq.\volumee, implies that the
ratio $\det B^2(z_i,\bar z_j)/|\det v_i(z_j)|^2$ does not depend on
$z_i$, $i\in I_N$, and therefore $\det B^2(z_i,\bar z_j)$ factorizes
into a product of a holomorphic times an antiholomorphic function of
$z_1,\ldots,z_N$. This is a special case of a more general theorem.

\vskip 6pt

\newth\fattorizza{Theorem}{Fix $n\in\NN_+$ and set
$$B_A(z,\bar w):=\sum_{i,j=1}^g\omega_i(z)A_{ij}\bar\omega_j(w)\ ,
$$ where $A$ is a complex $g\times g$ matrix. Then, for all $z_i,w_i\in C$, $i\in I_{N_n}$,
\eqn\fattor{\det B_A^n(z_i,\bar
w_j)=\bigl|\kenne[\phi^n]\bigr|^{-2}\det\phi^n(z_1,\ldots,z_{N_n})
\det\bar\phi^n(w_1,\ldots,w_{N_n}){K}_n(A)\ ,} where
$\{\phi^n_i\}_{i\in I_{N_n}}$ is an arbitrary basis of $H^0(K_C^n)$
and
\eqn\effennea{\eqalign{{K}_n(A)=\sum_{{i_{N_n}>\ldots>
i_1=1\atop j_{N_n}>\ldots
> j_{1}=1}}^{M_n}&\kenne[\omega\cdots\omega_{i_1},\ldots,\omega\cdots\omega_{i_{N_n}}]\cr
&\qquad\cdot {|A\ldots A|^{i_1\ldots i_{N_n}}_{j_1\ldots
j_{N_n}}\over
\prod_{k=1}^{{N_n}}\chi_{i_k}\chi_{j_k}}\,\bar\kenne[\omega\cdots\omega_{j_1},\ldots,
\omega\cdots\omega_{j_{N_n}}]\ .}}
 Furthermore, for $n\geq2$
\eqn\fattorbb{\det B_A^n(z_i,\bar
z_j)=\Bigl|\theta_\Delta\big(\smsum_{1}^{N_n}z_i\big)\prod_{i<j}^{N_n}E(z_i,z_j)\prod_{1}^{N_n}\sigma(z_i)^{2n-1}\Bigr|^2{K}_n(A)\ .
}}

\vskip 6pt

\noindent {\bf Proof.}
Observe that
$$\eqalignno{B_A^n(z_i,\bar w_j)=&\sum_{{k_1,\ldots,k_n=1\atop
l_1,\ldots,l_n=1}}^g\omega_{k_1}(z_i)\cdots\omega_{k_n}(z_i)A_{k_1l_1}\cdots
A_{k_nl_n}\bar\omega_{l_1}(w_j)\cdots\bar\omega_{l_n}(w_j)\cr\cr
=&\sum_{k,l=1}^{M_n} \omega\cdots\omega_k(z_i){(A\cdots A)_{kl}\over
\chi_k\chi_l}\bar\omega\cdots\bar\omega_l(w_j)\ , }$$ with the
notation of section \notation. Then
$$\det B_A^n(z_i,\bar w_j)=\sum_{{k_1,\ldots,k_{N_n}=1\atop l_1,\ldots,
l_n=1}}^{M_n}\sum_{s\in\perm_{N_n}}\sgn(s)\prod_{i=1}^{N_n}\omega\cdots\omega_{k_i}(z_i)
\bar\omega\cdots\bar\omega_{l_i}(w_{s_i}){(A\cdots A)_{k_il_i}\over \chi_{k_i}\chi_{l_i}}\ ,$$
and by defining $m_{s_i}:=l_i$, $i\in I_{M_n}$,
$$\eqalignno{\det B_A^n(z_i,\bar w_j)=&\sum_{{k_1,\ldots,k_{N_n}=1\atop
m_1,\ldots,
m_{N_n}=1}}^{M_n}{|A\ldots A|^{k_1\ldots k_{N_n}}_{m_1\ldots m_{N_n}}\over
\prod_{i=1}^{{N_n}}\chi_{k_i}\chi_{m_i}}\prod_{i=1}^{N_n}\omega\cdots\omega_{k_i}(z_i)
\bar\omega\cdots\bar\omega_{m_i}(w_i)\cr\cr
=&\sum_{{k_{N_n}>\ldots> k_1=1\atop m_{N_n}>\ldots
> m_{1}=1}}^{M_n}{|A\ldots A|^{k_1\ldots k_{N_n}}_{m_1\ldots m_{N_n}}\over
\prod_{i=1}^{{N_n}}\chi_{k_i}\chi_{m_i}}\sum_{r,s\in\perm_{N_n}}\sgn(r)\sgn(s)\prod_{i=1}^{N_n}\omega\cdots\omega_{k_{r_i}}(z_i)
\bar\omega\cdots\bar\omega_{m_{s_i}}(w_i)\cr\cr
=&\sum_{{k_{N_n}>\ldots> k_1=1\atop m_{N_n}>\ldots
> m_{1}=1}}^{M_n}{|A\ldots A|^{k_1\ldots k_{N_n}}_{m_1\ldots m_{N_n}}\over
\prod_{i=1}^{{N_n}}\chi_{k_i}\chi_{m_i}}\det_{{i\in \{k_1,\ldots,k_{N_n}\}\atop j\in
I_{N_n}}}\omega\cdots\omega_{i}(z_j)\det_{{i\in \{m_1,\ldots,m_{N_n}\}\atop j\in
I_{N_n}}}\bar\omega\cdots\bar\omega_{i}(w_j)\ .
}$$
By Eq.\ratiodef, for an arbitrary basis $\{\phi_i^n\}_{i\in I_{N_n}}$ of $H^0(K_C^n)$
$$\det_{{i\in \{k_1,\ldots,k_{N_n}\}\atop j\in
I_{N_n}}}\omega\cdots\omega_i(z_j)=\det\phi^n(z_1,\ldots,z_{N_n})
{\kenne[\omega\cdots\omega_{k_1},\ldots,\omega\cdots\omega_{k_{N_n}}]\over\kenne[\phi^n]}\
,$$ leading to \fattor. Eq.\fattorbb\ then follows by
Eq.\detthetaii. \hfill$\square$

\vskip 6pt

\newrem\kditaudued{Remark}{Since by \volumee\ $\det B^2(z_i,\bar z_j)$ is positive definite, it follows that the $K_2(\tau_2^{-1})>0$.
Even if there are stringent arguments suggesting that
$K_n(\tau_2^{-1})>0$ also for $n>2$, a complete proof is still lacking.}

\vskip 6pt

\newth\nuovabasse{Corollary}{Let $\{\phi^n_i\}_{i\in I_{N_n}}$ be an arbitrary basis of $H^0(K_C^n)$ and fix the points
$w_1,\ldots,w_{N_n}\in C$ in such a way that $\det\,
\phi^n_i(w_j)\neq0$. \eqn\bdf{\alpha^n_i(z)=B^n_A(z,\bar w_i)\ ,}
$i\in I_{N_n}$, is a basis of $H^0(K_C^n)$ whenever $K_n(A)\neq0$.}

\vskip 6pt

\noindent {\bf Proof.} Immediate by Theorem \fattorizza.
\hfill$\square$

\vskip 6pt

\newsec{Geometry of $\Theta_s$, $K(p_3,\ldots,p_g)$
and Mumford isomorphism}\seclab\Mumfsec

In this section, we first use the construction of section
\secdivspin\ to derive an expression for the Mumford form which does
not involve any determinant of holomorphic $1$-differentials.
Furthermore, we express the volume form $d\nu_{|_{\hat\M_g}}$ in
terms of the Mumford form, a result previously known for $g=2$ and
$g=3$ only. Next, by means of the Mumford isomorphism we investigate
the modular properties of $K(p_3,\ldots,p_g)$ in order to construct
sections of bundles on $\M_g$. For $g=2$ and $g=3$ such sections
reproduces the building blocks for the Mumford form. In the case of
$g=4$, a modular form on the Jacobian locus is obtained, which is
proportional to the Hessian of the theta function evaluated on
$\Theta_s$. This is a remarkable result in view of
\Farkas\grushsalv, where it is shown that the vanishing of such a
Hessian on the Andreotti-Mayer locus $\N_0=\J_4\cup \theta_{null}$,
where $\theta_{null}\subset\ppavmod_4$ is the locus of the ppav's
with a vanishing theta-null, characterizes the intersection
$\J_4\cap \theta_{null}$. This indicates that the sections on $\M_g$
built in terms of $K$ may be considered as generalizations to $g>4$
of such a Hessian, thus providing a tool for the analysis of the
geometry of $\Theta_s$ and of the Andreotti-Mayer locus $\N_{g-4}$.
We explicit construct such sections for even genus and for the case
$g=5$.

\vskip 6pt

\subsec{Mumford isomorphism and Siegel volume form}

Let ${\cal C}_g \sopra{\longrightarrow}{{}_\pi}{\cal M}_g$ be the
universal curve over ${\cal M}_g$ and $L_n=R\pi_*(K^n_{{\cal
C}_g/{\cal M}_g})$ the vector bundle on ${\cal M}_g$ of rank
$(2n-1)(g-1)+\delta_{n1}$ with fiber $H^0(K_C^n)$ at the point of
${\cal M}_g$ representing $C$. Let $\lambda_n:=\det L_n$ be the
determinant line bundle. According to Mumford
\ref\Mumford{D.~Mumford, Stability of projective varieties, {\it
Enseign.\ Math.\ } {\bf 23} (1977), 39-110.}
$$
\lambda_n\cong\lambda_1^{\otimes c_n}\ ,
$$
where $c_n=6n^2-6n+1$, which corresponds to (minus) the central charge of
the chiral $b-c$ system of conformal weight $n$
\BonoraCJ. The Mumford form $\mu_{g,n}$ is the unique, up to a
constant, holomorphic section of
$\lambda_n\otimes\lambda_1^{-\otimes c_n}$ nowhere vanishing on
${\cal M}_g$.

\vskip 6pt

\newth\nuovamumfordformm{Theorem}{Let $\{\phi^n_i\}_{i\in I_{N_n}}$ be a basis of
$H^0(K_C^n)$, $n\geq2$. For any set $\{\nu\}\equiv\{\nu_i\}_{i\in
I_{2n-1}}$ of non-singular odd spin structures, the Mumford form is,
up to a universal constant,
$$
\mu_{g,n}=F_{g,n}[\phi^n]{\phi^n_1\wedge\cdots\wedge\phi^n_{N_n}\over
(\omega_1\wedge\cdots\wedge\omega_g)^{c_n}}\ , $$ where
\eqn\nuovamumfordd{F_{g,n}[\phi^n]:= {\prod_{k=1}^{2n-1}
\kuno_{\nu_k}^{2n-1}[\omega]\over \kenne_{\{\nu\}}[\phi_n]} \ ,}
where
$$\nu_k''+\tau\nu_k'=I(\smsum_{i=1}^{g-1}p_k^i-\Delta) \ ,
$$
$k\in I_{2n-1}$, where $\kenne_{\{\nu\}}[\phi_n]$ is defined in
\dettheaaaBV\ and $\kuno_{\nu_k}^{2n-1}[\omega]$ in \nuooovaa, with
$y_k^i\equiv p_{k+1}^i$, $p_k^g\equiv p_{k+2}^1$, $k\in I_n$, $i\in
I_{g-1}$, and $p^i_{n+j}\equiv p^i_n$, $j=1,2$. }

\vskip 6pt

\noindent {\bf Proof.} By construction $F_{g,n}$ may depend only on
the marking of $C$. However, the transformation properties of the
theta functions derived in section \Riemanndef, Lemma
\riemannconstmod\ in particular, show that $\mu_{g,n}$ is
independent of the marking of $C$. $\hfill\square$

\vskip 6pt

\newrem\sumumf{Remark}{ Eq.\nuovamumfordd\ expresses, for any $n\geq2$ and
$g\geq2$, the Mumford form in terms of points defining odd spin
structures only. Explicit expressions of the Mumford form were
derived in
\lref\BelavinCY{
  A.~A.~Belavin and V.~G.~Knizhnik,
  Algebraic geometry and the geometry of quantum strings,
  {\it Phys.\ Lett.\ } B {\bf 168} (1986), 201-206.
}
\lref\BeilinsonZW{
  A.~A.~Beilinson and Y.~I.~Manin,
  The Mumford form and the Polyakov measure in string theory,
  {\it Commun.\ Math.\ Phys.\ }  {\bf 107} (1986), 359-376.
}
\lref\AlvarezGaumeVM{
  L.~Alvarez-Gaume, J.~B.~Bost, G.~W.~Moore, P.~C.~Nelson and C.~Vafa,
  Bosonization on higher genus Riemann Surfaces,
 {\it Commun.\ Math.\ Phys.\ } {\bf 112} (1987), 503-552.
}
\lref\VerlindeKW{
  E.~P.~Verlinde and H.~L.~Verlinde,
  Chiral bosonization, determinants and the string partition function,
{\it  Nucl.\ Phys.\ } B {\bf 288} (1987), 357-396.
} \refs{\BelavinCY\BeilinsonZW\AlvarezGaumeVM{-}\VerlindeKW} and
\FayMAM. In particular, Fay in \FayMAM\ provides the following
expression \eqn\mumfordd{
F_{g,n}[\phi^n]:={\deltadiv\bigl(\sum_{1}^{N_n}p_i\bigr)\prod_{i<j}^{N_n}E(p_i,p_j)\over
\det \phi^n_i(p_j)\prod_{1}^{N_n}c(p_i)^{2n-1\over g-1}}=
{\theta((2n-1)\K^p)\over W[\phi^n](p)c(p)^{(2n-1)^2}}\ , } where
\eqn\cppv{c(p):={\deltadiv\bigl(\sum_{1}^gp_i-y\bigr)
\prod_1^g\sigma(p_i,p)\prod_{i<j}E(p_i,p_j)\over
\det\omega_i(p_j)\sigma(y,p)
\prod_1^gE(y,p_i)}={1\over\kuno[\omega]\sigma(p)^{g-1}}\ ,} for all
$p,p_1,\ldots,p_g,y\in C$. By \dettheta\ and \detthetaii, such an
expression corresponds to
\eqn\efffe{F_{g,n}[\phi^n]={\kuno[\omega]\over
\kenne[\phi^n]}^{(2n-1)^2}\ .} Note that while $\det\omega_i(p_j)$
explicitly appears as a building block in Eqs.\mumfordd\ and \efffe,
Theorem \nuovamumfordformm\ provides an expression without such a
determinant.}

\vskip 6pt

The previous investigations allow now a straightforward derivation
of the relation between the Mumford form and the volume form on
$\hat\M_g$, induced by the Siegel metric.

\newth\corolllary{Theorem}{\eqn\nicessimo{|\kdue[\phi^2]|^{-2}|\phi_1^2\wedge\cdots\wedge\phi_N^2|^2=
\Bigl({2\over i}\Bigr)^N{d\nu_{|\hat\M_g}\over K_2(\imtau^{-1})}\ ,}
where $K_2$ is given in Eq.\effennea. In particular,
\eqn\gahsty{|\mu_{g,2}|^2
|\omega_1\wedge\cdots\wedge\omega_g|^{26}=\Bigl({2\over
i}\Bigr)^N{|\kuno[\omega]|^{18}\over
K_2(\imtau^{-1})}d\nu_{|\hat\M_g}\ . }}

\vskip 6pt

\noindent {\bf Proof.} By \ds\ildettt\ and \fattor\
$$
{\textstyle\bigwedge_{1}^N}\!(\Xi_i\wedge\bar\Xi_i)=\Bigl({2\over
i}\Bigr)^N{|\det v_i(z_j)|^2\over\det B^2(z_i,\bar z_j)}d
\nu_{|\hat\M_g}=\Bigl({2\over i}\Bigr)^N|\kappa[v]|^2
K_2^{-1}(\imtau^{-1})d \nu_{|\hat\M_g} \ ,
$$
which reproduces \nicessimo\ after the base change $\{v_i\}_{i\in I_N}\to\{\phi^2_i\}_{i\in I_N}$. \hfill$\square$

\vskip 6pt

It has been shown in
\ref\BelavinTV{
  A.~A.~Belavin, V.~Knizhnik, A.~Morozov and A.~Perelomov,
  Two and three loop amplitudes in the bosonic string theory,
  {\it Phys.\ Lett.\ } B {\bf 177} (1986), 324-328. 
  }
\ref\MorozovDA{
  A.~Morozov,
  Explicit formulae for one, two, three and four loop string amplitudes,
 {\it  Phys.\ Lett.\ } B {\bf 184} (1987), 171-176.
} that \eqn\gduee{F_{2,2}[\omega\omega]={c_{2,2}\over
\Psi_{10}(\tau)}\ , } where $\Psi_{10}$ is the modular form of
weight 10
$$
\Psi_{10}(\tau):=\prod_{a,b\, even}\theta\left[^a_b\right](0)^2\ ,
$$ where the product is over the 10 even characteristics of $g=2$.
It has been shown in \DHokerQP\ that $c_{2,2}=1/\pi^{12}$.
Furthermore, it has been conjectured that
 \BelavinTV\MorozovDA\ \eqn\gtreee{F_{3,2}[\omega\omega]={c_{3,2}\over \Psi_{9}(\tau)}\
, } with $\Psi_9(\tau)^2\equiv\Psi_{18}(\tau)$
$$
\Psi_{18}(\tau):=\prod_{a,b\, even}\theta\left[^a_b\right](0)\ ,
$$ where the product is over the 36 even characteristics of $g=3$
and $c_{3,2}=1/2^6\pi^{18}$
\ref\DHokerCE{
  E.~D'Hoker and D.~H.~Phong,
  Asyzygies, modular forms, and the superstring measure. II,
  {\it Nucl.\ Phys.\ } B {\bf 710} (2005), 83-116
  [arXiv:hep-th/0411182].
}.

\vskip 6pt

\newth\detbarprt{Proposition}{For $g=2$
\eqn\dajjie{\kuno[\omega]^6=-{1\over\pi^{12}}
{\deltadiv(p_1+p_2-p_3)\deltadiv(p_1+p_3-p_2)\deltadiv(p_2+p_3-p_1)\over
\Psi_{10}(\tau)\deltadiv(p_1+p_2+p_3)[E(p_1,p_2)E(p_1,p_3)E(p_2,p_3)
\sigma(p_1)\sigma(p_2)\sigma(p_3)]^2}
 \ .}
For $g=3$ \eqn\dajjieii{\kuno[\omega]^5={
\sum_{s\in\perm_6'}\sgn(s)\prod_{k=1}^4[\deltadiv\bigl(\sum_{i=1}^3p_{d^k_i(s)}-y\bigr)
\prod_{i<j}^3E(p_{d^k_i(s)},p_{d^k_j(s)})]\over 2^6
15\pi^{18}\Psi_9(\tau)\deltadiv\bigl(\sum_{1}^6p_i\bigr)\prod_{i=1}^6
\sigma(p_i)\sigma(y)^4\prod_{i=1}^6E(y,p_i)^2\prod_{i<j}^6E(p_i,p_j)}\
. } Furthermore, for $g=2$ \eqn\dajjieoobbb{(\det
\omega\omega_i(p_j))^2=-{1\over\pi^{12}}
{\prod_{i=1}^3\deltadiv(\sum_{j=1}^3p_j-2p_i)^3\over
\Psi_{10}(\tau)\deltadiv(\sum_1^3p_i)\prod_{i<j}E(p_i,p_j)^4} \ .}
 }

\vskip 6pt

\noindent {\bf Proof.} Eq.\dajjie\ follows by simply replacing the
expression of $\det\omega\omega_i(p_j)$ in \due\ in
$F_{2,2}[\omega\omega]$, then using
$c(p)={\sigma(p)^{1-g}/\kuno[\omega]}$. Similarly, by \tre\
$${\det\omega\omega(p_1,\ldots,p_6)\over \kuno[\omega]^4}=
{\prod_{1}^6\sigma(p_i)^2
\sum_{s\in\perm_6'}\sgn(s)\prod_{k=1}^4[\deltadiv\bigl(\sum_{i=1}^3p_{d^k_i(s)}-y\bigr)
\prod_{i<j}^3E(p_{d^k_i(s)},p_{d^k_j(s)})]\over
15\sigma(y)^4\prod_{i=1}^6E(y,p_i)^2}\ ,$$ and \dajjieii\ follows by
the identity $\smprod_{i=1}^6c(p_i)^{-{3\over
2}}=\kuno[\omega]^9\smprod_{i=1}^6\sigma(p_i)^3$. Eq.\dajjieoobbb\
follows by direct computation. \hfill$\square$

\vskip 6pt

\subsec{Modular properties of $K(p_3,\ldots,p_g)$ and the geometry
of $\Theta_s$}

For each $n\in\ZZ_{>0}$, let us consider the rank $N_n$ vector
bundle $L_n$ on $\M_g$, defined in the previous subsection, whose
fiber at the point corresponding to a curve $C$ is $H^0(K_C^n)$. A
general section $s\in L_n^m$, $i>1$, admits the local expression on
an open set $U\subset\M_g$ \eqn\locsec{s(p)=\sum_{i_1,\ldots,i_m\in
I_{N_n}}s_{i_1\ldots
i_m}(p)\phi_{i_1}\otimes\phi_{i_2}\otimes\cdots\otimes\phi_{i_m}\
,\qquad p\in U\subset \M_g\ , } with respect to a set
$\{\phi_i\}_{i\in I_{N_n}}$ of linearly independent local sections
of $L_n$ on $U$.

For each non-hyperelliptic $C$ of genus $g\ge 3$,
$k(p_3,\ldots,p_g)$ as defined in \kpiccola, is a holomorphic
$(g-3)$-differential in each variable, and is symmetric (for $g$
even) or anti-symmetric (for $g$ odd) in its $g-2$ arguments. Hence,
\eqn\kapsec{k:=\sum_{i_1,\ldots,i_{g-2}\in I_{N_{g-3}}}k_{i_1\ldots
i_{g-2}} \phi_{i_1}\otimes\cdots\otimes\phi_{i_{g-2}}\ , } can be
naturally seen as an element of $E_g$, where
$$E_g:=\left\{\vcenter{\vbox{\halign{\strut\hskip 6pt $ # $ \hfil & \hskip 2cm$ # $ \hfil\cr
\Sym^{g-2}H^0(K_C^{g-3})\ , &g \hbox{ even}\ ,\cr
\bigwedge^{g-2}H^0(K_C^{g-3})\ , &g \hbox{ odd}\ ,\cr}}}\right.
$$
for a fixed basis $\{\phi_i\}_{i\in I_{N_{g-3}}}$ of
$H^0(K_C^{g-3})$. The definition can be extended in a continuous way
to hyperelliptic curves, by setting $k_{i_1\ldots i_{g-2}}\equiv 0$
in this case. At genus $g=3$, $k(p_3)$ is a holomorphic function on
$C$ and therefore is a constant. Furthermore, Eq.\kpiccola\ also
makes sense at genus $g=2$; in this case, $k$ is again a constant.
For $g> 3$, let us define $\EE_g$ by
$$\EE_g:=\left\{\vcenter{\vbox{\halign{\strut\hskip 6pt $ # $ \hfil & \hskip 2cm$ # $ \hfil\cr
\Sym^{g-2}L_{g-3}\ , &g \hbox{ even}\ ,\cr \bigwedge^{g-2}L_{g-3}\ ,
&g \hbox{ odd}\ .\cr}}}\right.
$$
In view of Eqs.\locsec\ and \kapsec, it is natural to seek for a section
$\kk\in\EE_g$ such that, at the point $p_C\in\M_g$ corresponding to the curve $C$, it satisfies
$$E_g\ni k(p_3,\ldots,p_g)\cong \kk(p_C)\in(\EE_g)_{|p_C}\ ,$$ under the identification
$(\EE_g)_{|p_C}\cong E_g$. On the other hand, $k(p_3,\ldots,p_g)$
is not modular invariant, and then it does not correspond to a well-defined element of
$E_g$ for each $p_C\in\M_g$. The correct statement is given by the following
theorem.

\vskip 6pt

\newth\linebun{Theorem}{
$$\kk:=\kuno[\omega]^{g-8}k\otimes(\omega_1\wedge\ldots\wedge\omega_g)^{12-g}\ ,$$
is a holomorphic section of $\lambda_1^{12-g}$ for $g=2,3$ and of $\EE_g\otimes\lambda_1^{12-g}$ for $g>3$,
which vanishes only in the
hyperelliptic locus for $g\geq3$.}

\vskip 6pt

\noindent {\bf Proof.} Let us derive the modular properties of
$$\kuno[\omega]^{g-8}k(p_3,\ldots,p_g)\ .
$$
Eq.\kpiccola\ and the identity $\kuno[\sigma]=\kuno[\omega]/\det\omega_i(p_j)$ yield
$$\kuno[\omega]^{g-8}k(p_3,\ldots,p_g)=(-)^{g+1}
c_{g,2}{\kdue[v]\over\kuno[\omega]^{9}}(\det\omega_i(p_j))^{g+1}
\ .
$$
By Eq.\efffe, it follows that $\kdue[v]/\kuno[\omega]^{9}$ has a
simple modular transformation
$${\kdue[v]\over\kuno[\omega]^{9}}\to {\kdue[v]\over\kuno[\omega]^{9}}(\det(C\tau+D))^{-13}\
,\qquad\quad \left(\matrix{A & B\cr C & D}\right)\in {\rm
Sp}(2g,\ZZ)\ ,
$$
and, by using the modular transformation $\det\omega_i(p_j)\to\det\omega_i(p_j)\det(C\tau+D)$, we
obtain
$$\kuno[\omega]^{g-8}k
\to \kuno[\omega]^{g-8}k (\det(C\tau+D))^{g-12}\  .
$$
Hence,
$\kuno[\omega]^{g-8}k\otimes(\omega_1\wedge\ldots\wedge\omega_g)^{12-g}$
is modular invariant and determines a section of
$\Sym^{g-2}E_{g-3}\otimes\lambda_1^{12-g}$ on $\M_g$. Since
$\kappa[\omega]\neq 0$ for all $C$, $\kk=0$ at the point
corresponding to the $C$ if and only if $k(p_3,\ldots,p_g)=0$ for
all $p_3,\ldots,p_g\in C$, or, equivalently, if and only if $C$ is
hyperelliptic. \hfill$\square$

\vskip 6pt

For $g=2$ the section $\kk$ corresponds to
$$\kk=\kappa[\omega]^6k(\omega_1\wedge\omega_2)^{10}\ ,$$
and for $g=3$
$$
\kk=\kappa[\omega]^5k(\omega_1\wedge\omega_2\wedge\omega_3)^{9}\ .
$$
Note that, for $g=2,3$, Eqs.\detdue\ and \dettre\ lead to the following relations
$${\kdue[v]\over \kuno[\sigma]^{g+1}}={\kdue[\omega\omega]\over \kuno[\omega]^{g+1}}\ ,
$$
and, together with \gduee\ and \gtreee, we obtain the identification
$$\eqalign{k=6\pi^{12}\kappa[\omega]^6\Psi_{10}\ ,\qquad & g=2\ ,\cr
k=15\cdot 2^6\kappa[\omega]^5\pi^{18}\Psi_{9}\ ,\qquad & g=3\ , }$$
recovering the results of Proposition \detbarprt.

\vskip 6pt

Let $C$ be a non-hyperelliptic curve of genus $g=4$. In this case,
$k(p_3,p_4)$ is a holomorphic $1$-differential in both $p_3$ and
$p_4$, symmetric in its arguments. Then,
$$k^{(4)}:={\det k(p_i,p_j)\over (\det\omega_i(p_j))^2}=\det k_{ij}\ ,$$
is a meromorphic function on $C$ in each $p_i$, $i\in I_4$.

\vskip 6pt

\newth\kappaquattro{Proposition}{The function $k^{(4)}$ is a constant on $C$ that depends only on the
choice of the marking. Furthermore, $k^{(4)}=0$ if and only if $C$ is hyperelliptic or if $C$ is
non-hyperelliptic and admits a (necessarily even)
singular spin structure.}

\vskip 6pt

\noindent{\bf Proof.} Let us suppose that, for a suitable
choice of $p_1,\ldots,p_4\in C$, $\{k(p_i,z)\}_{i\in I_4}$ is a basis of $H^0(K_C)$. Then,
the determinant
$${\det k(p_i,z_j)\over \det\omega_i(z_j)}\ ,$$
does not depend on the points $z_1,\ldots,z_4\in C$. Hence, the ratio
$${\det k(p_i,z_j)\over \det\omega_i(z_j)\det\omega_i(p_j)}\ ,$$ is a non-vanishing constant on $C$.
In particular, by taking $p_i=z_i$, it follows that such a constant is $k^{(4)}$. On the contrary, if
for all $p_1,\ldots,p_4\in C$, the holomorphic $1$-differentials $k(p_i,z)$, $i\in I_4$, are
linearly dependent, then $k^{(4)}$ vanishes identically.

Such a construction shows that $k^{(4)}$ vanishes if and only if
$k(p_i,z)$, $i\in I_4$, are linearly dependent for all
$p_1,\ldots,p_4\in C$. If $C$ is hyperelliptic, then $k(p_i,p_j)=0$
for all $p_i,p_j\in C$ and $k^{(4)}=0$. Assume that $C$ admits a
singular spin structure $\alpha$ and let $L_\alpha$ be the
corresponding holomorphic line bundle with $L_\alpha^2\cong K_C$.
This implies that $\Theta_s$ consists of a unique point of order $2$
in the Jacobian torus. For each $p\in C$, the holomorphic
$1$-differential $k(p,z)$ is the square of an element of
$H^0(L_\alpha)$; by varying $p\in C$, such $1$-differentials span
the image of the map $\varphi:{\rm Sym}^2H^0(L_\alpha)\to H^0(K_C)$.
If $\alpha$ is even, then $h^0(L_\alpha)=2$ and ${\rm
Sym}^2H^0(L_\alpha)$ has dimension three, so that $\varphi$ cannot
be surjective and $k^{(4)}=0$. If $\alpha$ is odd, then
$h^0(L_\alpha)=3$ so that, for each point $p\in C$,
$h^0(L_\alpha\otimes\O(-p))\ge 2$; if $h_1,h_2$ span
$H^0(L_\alpha\otimes\O(-p))$, then $h_1/h_2$ is a non-constant
meromorphic function with $2$ poles and $C$ is hyperelliptic.

Suppose that $C$ is
non-hyperelliptic and does not admit a singular spin structure.
Then, $\Theta_{s}$ consists of $2$ distinct points, $e$ and $-e$. Let us first observe that if
there exist two points $p,q\in C$ such that $I(p-q)=2e$, then they are unique. For, if
$I(\tilde p-\tilde q)=2e=I(p-q)$, then $p+\tilde q -\tilde p-q$ is the divisor of a meromorphic
function on $C$. But, since $C$ is non-hyperelliptic, the unique meromorphic function with less that
$3$ poles are the constants and, since $p\neq q$ (because $2e\neq 0$ in $J_0(C)$), it follows that
$\tilde p=p$ and $\tilde q=q$.

Also, observe that $K(z,z)$ is not identically vanishing as a function of $z$; since $C$ is compact,
$K(z,z)$ has only a finite number of
zeros. Fix a point $p_1\in C$ and define
$x_1,x_2,y_1,y_2\in
C$ by
$$I(p_1+x_1+x_2-\Delta)=e\ ,\qquad I(p_1+y_1+y_2-\Delta)=-e\ .$$
Then the divisor of $k(p_1,z)$ with respect to $z$ is
$2p_1+x_1+x_2+y_1+y_2$. Observe that at least one between $x_1$ and
$x_2$ is distinct from $y_1$ and $y_2$, since otherwise we would
have $e=-e$. We choose $p_1$ in such a way that
$p_1,x_1,x_2,y_1,y_2$ are distinct from the zeros of $K(z,z)$ and
from the points $p,q$ such that $I(p-q)=2e$ (if they exist). Note
that such a condition can always be fulfilled, since it is
equivalent to require that $p_1$ is distinct from the zeros of
$k(p,\cdot)$, $k(q,\cdot)$ and $k(w,\cdot)$ for each $w$ such that
$K(w,w)=0$. Then, the points for which such a condition is not
satisfied is a finite set.

Set $p_2:=x_1$
and $p_3:=y_1$. The divisor of $k(p_3,z)$ is $(k(p_3,z))=2p_3+p_1+y_2+z_1+z_2$, where $z_1,z_2$
satisfy
$$I(p_3+z_1+z_2-\Delta)=e\ .$$ Since the condition on the choice of $p_1$ implies
$K(p_3,p_3)\neq 0$, it follows that $z_1$ and $z_2$ are distinct from $p_3$.
Set $p_4:=z_1$, so that $\det k(p_i,p_j)=k(p_1,p_4)^2k(p_2,p_3)^2$. The identities
$$I(p_1+p_2+x_2-p_3-p_4-z_2)=0\ ,\qquad\qquad I(p_4+z_2-p_1-y_2)=2e\ ,$$ imply that $p_4$ and $z_2$ are
distinct from $p_1,p_2,x_2$ (for example, if $p_4=x_2$, then $p_1+p_2-p_3-z_2$ is the divisor of a
meromorphic function and $C$ is hyperelliptic) and  from $y_2$ (if $p_4=y_2$, then $I(z_2-p_1)=2e$,
counter the requirement that $p_1$ is distinct from $q$ and $p$). Therefore, $k(p_1,p_4)k(p_2,p_3)\neq
0$ and then $k^{(4)}\neq 0$.
\hfill$\square$

\vskip 6pt

By Propositions \kappaquattro\ and \linebun, it follows that, for
$g=4$,
$$\kk^{(4)}:=\kappa[\omega]^{-16}\det
k_{ij}(\omega_1\wedge\cdots\wedge\omega_4)^{34}\ ,$$ is a
holomorphic section of $\lambda_1^{34}$ vanishing only on the
hyperelliptic locus, with a zero of order $4[(3g-3)-(2g-1)]=8$, and
on the locus of Riemann surfaces with an even singular spin
structure, with a zero of order $1$. By Eq.\thetaKFour, the
following relation holds
$$k^{(4)}=A^4\det_{ij\in I_4}\theta_{ij}(e)\ ,
$$
where the constant $A$ depends on the moduli. Recently, it has been
shown that the Hessian $\det_{ij\in I_4}\theta_{ij}(e)$ plays a key
role in the analysis of the Andreotti-Mayer loci at genus $4$ and in
the corresponding applications to the Schottky problem
\Farkas\grushsalv. Whereas no natural generalization of such a
Hessian exists at genus $g>4$, the section $\kk^{(4)}$ is the $g=4$
representative of a set of sections $\kk^{(g)}$ of a tensor power of
$\lambda_1$ on $\M_g$, defined for each even $g\ge 4$.

\newrem\geralkappa{Definition}{Let $C$ be a curve of even genus $g\ge 4$. Fix $N_{g-3}=h^0(K_C^{g-3})$
points $p_1,\ldots,p_{N_{g-3}}\in C$ and let $\{\phi_i\}_{i\in
I_{N_{g-3}}}$ be a basis of $H^0(K_C^{g-3})$. Set
$$\kk^{(g)}:={\kenne[\phi]^{g-2}\sum_{s^1,\ldots,s^{g-2}\in\perm_{N_{g-3}}}
\prod_{i=1}^{g-2}\sgn(s^i)\prod_{j=1}^{N_{g-3}}
k(p_{s_j^1},\ldots,p_{s_j^{g-2}})\over
N_{g-3}!\kappa[\omega]^{(2g-7)^2(g-2)+(8-g)N_{g-3}}\bigl(\det\phi(p_1,\ldots,p_{N_{g-3}})\bigr)^{g-2}}
(\omega_1\wedge\cdots\wedge\omega_g)^{d_g}\ ,
$$
where $d_g:=(12-g)N_{g-3}+(g-2)[6(g-3)(g-4)+1]$.
}

\vskip 6pt

\newth\generalkappa{Proposition}{For all the even $g\ge 4$, $\kk^{(g)}$ does not depend on the points
$p_1,\ldots,p_{N_{g-3}}\in C$ and on the basis $\{\phi_i\}_{i\in I_{N_{g-3}}}$ of
$H^0(K_C^{g-3})$ and determines a section of $\lambda_1^{d_g}$ on $\M_g$.}

\vskip 6pt

\noindent {\bf Proof.} Choose $(g-2)N_{g-3}$ points $p_1^i,\ldots,p_{N_{g-3}}^i\in C$, $i\in
I_{g-2}$ and note that
\eqn\permprod{\sum_{s^1,\ldots,s^{g-2}\in\perm_{N_{g-3}}}
\prod_{i=1}^{g-2}\sgn(s^i)\prod_{j=1}^{N_{g-3}}
k(p^i_{s_j^1},\ldots,p^i_{s_j^{g-2}})\ ,}
is a product of $g-3$ differentials in each $p_j^i$, $i\in I_{g-2}, j\in I_{N_{g-3}}$. Such a
product is completely anti-symmetric with respect to the permutations of each $N_{g-3}$-tuple
$(p^i_1,\ldots,p^i_{N_{g-3}})$, for all $i\in I_{g-2}$, so that it must be proportional to the
determinant
$\det\phi(p^i_1,\ldots,p^i_{N_{g-3}})$.
Therefore, the ratio of Eq.\permprod\ and
$$\prod_{i\in I_{g-2}}\det\phi(p^i_1,\ldots,p^i_{N_{g-3}})\ ,
$$
does not depend on the points $p_1^i,\ldots,p_{N_{g-3}}^i\in C$,
$i\in I_{g-2}$; in particular, by choosing, for each $j\in
I_{N_{g-3}}$, $p_j^1\equiv p_j^2\equiv\ldots\equiv p_j^{g-2}\equiv
p_j$, where $p_1,\ldots,p_{N_{g-3}}$ are the points in the
definition \geralkappa, it follows that $\kk^{(g)}$ is a constant as
a function of $C^{N_{g-3}}$. The proposition follows trivially by
Theorem \linebun\ and by the expression \efffe\ of the Mumford form,
with $n=g-3$.\hfill$\square$

\vskip 6pt

Definition \geralkappa\ and Proposition \generalkappa\ makes sense
also at odd genera; however, simple algebraic considerations show
that, in this case, $\kk^{(g)}$ is identically null on $\M_g$. In
general, there exist some non-trivial generalizations of $\kk^{(4)}$
at odd genus, but they are not as simple as the ones at even $g$. An
example at genus $g=5$ is
$$\eqalign{{(\omega_1\wedge\cdots\wedge\omega_5)^{164}\kdue[\phi]^4\over
\kuno[\omega]^{84}\bigl(\det\phi(p_1,\ldots,p_{12})\bigr)^4}
&\sum\nolimits_{i,j,k,l\in\perm_{12}}\sgn(i)\sgn(j)\sgn(k)\sgn(l)\cr
&\qquad \cdot k(p_{i_1},p_{i_2},p_{i_3})k(p_{i_4},p_{i_5},p_{j_{1}})
k(p_{i_6},p_{i_7},p_{k_{1}})k(p_{i_8},p_{i_9},p_{l_{1}})\cr &\qquad
\cdot
k(p_{i_{10}},p_{j_{2}},p_{j_{3}})k(p_{i_{11}},p_{k_{2}},p_{k_{3}})
k(p_{i_{12}},p_{l_{2}},p_{l_{3}})k(p_{j_{4}},p_{j_{5}},p_{j_{6}})\cr
&\qquad \cdot
k(p_{j_{7}},p_{j_{8}},p_{k_{4}})k(p_{j_{9}},p_{j_{10}},p_{l_{4}})
k(p_{j_{11}},p_{k_{5}},p_{k_{6}})k(p_{j_{12}},p_{l_{5}},p_{l_{6}})\cr
&\qquad \cdot
k(p_{k_{7}},p_{k_{8}},p_{k_{9}})k(p_{k_{10}},p_{k_{11}},p_{l_{7}})
k(p_{k_{12}},p_{l_{8}},p_{l_{9}})k(p_{l_{10}},p_{l_{11}},p_{l_{12}})\
,}
$$
which does not depend on the points $p_1,\ldots,p_{12}\in C$ and corresponds to a section of
$\lambda_1^{164}$ on $\M_5$.

\vskip 6pt

 {\bf Acknowledgements}. Work partially supported by the
European Community's Human Potential Programme under contract
MRTN-CT-2004-005104 ``Constituents, Fundamental Forces and
Symmetries of the Universe". R.V. thanks the Max-Planck-Insitut
f\"ur Mathematik, Bonn, for hospitality during completion of this
work.

\listrefs

\end